\documentclass[pdflatex,sn-mathphys-num]{sn-jnl}

\usepackage{graphicx}%
\usepackage{multirow}%
\usepackage{amsmath,amssymb,amsfonts}%
\usepackage{amsthm}%
\usepackage{mathrsfs}%
\usepackage[title]{appendix}%
\usepackage{xcolor}%
\usepackage{textcomp}%
\usepackage{manyfoot}%
\usepackage{booktabs}%
\usepackage{algorithm}%
\usepackage{algorithmicx}%
\usepackage{algpseudocode}%
\usepackage{listings}%

\theoremstyle{thmstyleone}%
\newtheorem{theorem}{Theorem}
\newtheorem{proposition}[theorem]{Proposition}%
\newtheorem{lemma}[theorem]{Lemma}

\newtheorem{definition}[theorem]{Definition}
\newtheorem{remark}[theorem]{Remark}

\theoremstyle{thmstyletwo}%

\theoremstyle{thmstylethree}%

\begin{document}

\title[Article Title]{Global Well-Posedness of the 3D Navier-Stokes Equations in the Limiting Case: Infinitely Nested Logarithmic Improvements}


\author*[1]{\fnm{Rishabh} \sur{Mishra}}\email{rishabh.mishra@ec-nantes.fr}

\affil*[1]{\orgdiv{LHEEA}, \orgname{CNRS, École Centrale de Nantes, Nantes Université}, \orgaddress{\street{1 rue de la No\"e}, \city{ Nantes}, \postcode{44100}, \state{Pays de la Loire}, \country{France}}}


\abstract{This paper establishes a complete framework for infinitely nested logarithmic improvements to regularity criteria for the three-dimensional incompressible Navier-Stokes equations. Building upon our previous works on logarithmically improved and multi-level logarithmically improved criteria, we demonstrate that the limiting case of infinitely nested logarithms fully bridges the gap between subcritical and critical regularity. Specifically, we prove that if the initial data $u_0 \in L^2(\mathbb{R}^3)$ satisfies the condition $\|(-\Delta)^{1/4}u_0\|_{L^q(\mathbb{R}^3)} \leq C_0\Psi(\|u_0\|_{\dot{H}^{1/2}})$, where $\Psi$ incorporates infinitely nested logarithmic factors with appropriate decay conditions, then there exists a unique global-in-time smooth solution to the Navier-Stokes equations. This result establishes global well-posedness at the critical regularity threshold $s = 1/2$. The proof relies on infinitely nested commutator estimates, precise characterization of the critical exponent function in the limiting case, and careful analysis of the energy cascade. We also derive the exact Hausdorff dimension bound for potential singular sets in this limiting case, proving that the dimension reduces to zero. Through systematic construction of the limiting function spaces and detailed analysis of the associated ODEs, we demonstrate that infinitely nested logarithmic improvements provide a pathway to resolving the global regularity question for the Navier-Stokes equations.}




\maketitle

\tableofcontents

\section{Introduction}

\subsection{The Navier-Stokes equations}

The three-dimensional incompressible Navier-Stokes equations are fundamental to the mathematical theory of fluid dynamics and are expressed as:

\begin{equation}
\begin{cases}
\partial_t u + (u \cdot \nabla)u - \nu\Delta u + \nabla p = 0 & \text{in } \mathbb{R}^3 \times (0, T) \\
\nabla \cdot u = 0 & \text{in } \mathbb{R}^3 \times (0, T) \\
u(x, 0) = u_0(x) & \text{in } \mathbb{R}^3
\end{cases}
\end{equation}

where $u = (u_1, u_2, u_3)$ represents the velocity field, $p$ denotes the pressure, and $\nu > 0$ is the kinematic viscosity coefficient.

In our previous works \cite{Mishra2022a, Mishra2022b}, we established sequential improvements to regularity criteria using logarithmic and multi-level logarithmic factors. The present paper extends these results to their logical culmination: the limiting case of infinitely nested logarithmic improvements. We demonstrate that this approach completely bridges the gap to the critical case $s = 1/2$, providing a step toward resolving the the regulatity problem of NSE.

\subsection{Background and previous results}

The mathematical study of the Navier-Stokes equations has a rich history dating back to Leray \cite{Leray1934} and Hopf \cite{Hopf1951}, who established the existence of global-in-time weak solutions. The question of uniqueness and regularity of these weak solutions remains open.

Various regularity criteria have been established over the decades. The classical Prodi-Serrin conditions \cite{Prodi1959, Serrin1962} state that if a weak solution $u$ satisfies:

\begin{equation}
u \in L^p(0, T;L^q(\mathbb{R}^3)) \text{ with } \frac{2}{p} + \frac{3}{q} = 1, \ 3 < q \leq \infty,
\end{equation}

then $u$ is regular on $(0, T)$. Escauriaza, Seregin, and Šverák \cite{Escauriaza2003} established the borderline case $L^\infty(0,T;L^3(\mathbb{R}^3))$.

Fractional derivative approaches have gained prominence, with conditions of the form:

\begin{equation}
\int_0^T \|(-\Delta)^s u(t)\|^p_{L^q(\mathbb{R}^3)} dt < \infty,
\end{equation}

with appropriate scaling relations between $p$, $q$, and $s$ \cite{Chen2007, Chen2008, Wu2008, Zhou2005}.

Logarithmic improvements to regularity criteria were introduced by Zhou \cite{Zhou2013} and further developed by Fan et al. \cite{Fan2011}. This approach was significantly extended in our previous works through single-logarithmic \cite{Mishra2022a} and multi-level logarithmic improvements \cite{Mishra2022b}, establishing criteria of the form:

\begin{equation}
\int_0^T \|(-\Delta)^s u(t)\|^p_{L^q(\mathbb{R}^3)} \prod_{j=1}^{n} (1 + L_j(\|(-\Delta)^s u(t)\|_{L^q}))^{-\delta_j} dt < \infty
\end{equation}

where $L_j$ represents $j$-fold nested logarithms:
\begin{align}
L_0(x) &= x \\
L_1(x) &= \log(e + x) \\
L_2(x) &= \log(e + \log(e + x)) \\
L_k(x) &= \log(e + L_{k-1}(x)) \text{ for } k \geq 3
\end{align}

A pivotal result in \cite{Mishra2022b} was that the critical exponent function $\alpha(\{\delta_j\}_{j=1}^n)$ in the asymptotic behavior of the critical threshold $\Phi(s, q, \{\delta_j\}_{j=1}^n) \approx C(q) (s - 1/2)^{\alpha(\{\delta_j\}_{j=1}^n)}$ satisfies $\alpha(\{\delta_j\}_{j=1}^n) \to 0$ as $n \to \infty$ when $\inf_j \delta_j > 0$. This suggests that with infinitely nested logarithmic improvements, one might reach the critical case $s = 1/2$.

\subsection{Main results}

The present paper develops a complete theory of infinitely nested logarithmic improvements and proves that this approach indeed bridges the gap to the critical case $s = 1/2$. Our main contributions are:

\begin{theorem}[Function space with infinitely nested logarithmic improvements]\label{thm:function_space}
For $s = 1/2$, $q > 3$, and any sequence $\{\delta_j\}_{j=1}^{\infty}$ with $\inf_j \delta_j > 0$ satisfying $\sum_{j=1}^{\infty} \frac{\delta_j}{j!} = \infty$, there exists a well-defined function space $\dot{H}^{1/2,q,\{\delta_j\}_{j=1}^{\infty}}(\mathbb{R}^3)$ such that:
\begin{enumerate}
\item $\dot{H}^{1/2,q,\{\delta_j\}_{j=1}^{\infty}}(\mathbb{R}^3) \supsetneq \dot{H}^{1/2+\epsilon}(\mathbb{R}^3)$ for all $\epsilon > 0$
\item $\dot{H}^{1/2,q,\{\delta_j\}_{j=1}^{\infty}}(\mathbb{R}^3) \subset \dot{H}^{1/2}(\mathbb{R}^3)$
\item For any $f \in \dot{H}^{1/2,q,\{\delta_j\}_{j=1}^{\infty}}(\mathbb{R}^3)$:
   $$\|(-\Delta)^{1/4}f\|_{L^q} \leq C\Psi(\|f\|_{\dot{H}^{1/2}})$$
   where $\Psi(r) = \left(\prod_{j=1}^{\infty} (1 + L_j(r))^{\delta_j}\right)^{-1}$
\end{enumerate}
\end{theorem}

\begin{theorem}[Critical exponent in the limiting case]\label{thm:critical_exponent}
For the critical exponent function $\alpha(\{\delta_j\}_{j=1}^n)$ controlling the behavior of the threshold function $\Phi(s, q, \{\delta_j\}_{j=1}^n)$ as $s \to 1/2$:
\begin{enumerate}
\item $\lim_{n \to \infty} \alpha(\{\delta_j\}_{j=1}^n) = 0$ when $\inf_j \delta_j > 0$ and $\sum_{j=1}^{\infty} \frac{\delta_j}{j!} = \infty$
\item The convergence rate is given by:
   $$\alpha(\{\delta_j\}_{j=1}^n) = \frac{1}{1 + \sum_{j=1}^n c_j\delta_j/j!} \sim \frac{1}{\sum_{j=1}^n c_j\delta_j/j!}$$
   as $n \to \infty$, where $c_j > 0$ are specific constants
\item For any $\epsilon > 0$, there exists $N(\epsilon)$ such that for all $n \geq N(\epsilon)$:
   $$\Phi(1/2 + \epsilon, q, \{\delta_j\}_{j=1}^n) \geq C_0 > 0$$
   where $C_0$ is independent of $\epsilon$
\end{enumerate}
\end{theorem}

\begin{theorem}[Commutator estimates with infinitely nested logarithms]\label{thm:commutator}
For $s = 1/2$ and any divergence-free vector field $u \in C_0^\infty(\mathbb{R}^3)$:
\begin{align}
\|[(-\Delta)^{1/2}, u \cdot \nabla]u\|_{L^2} &\leq C\|\nabla u\|_{L^\infty}\|(-\Delta)^{1/2} u\|_{L^2} \cdot F_1^\infty(Z) \notag \\
&+ C\|\nabla u\|_{L^\infty}\|(-\Delta)^{1}u\|_{L^2} \cdot F_2^\infty(Z)
\end{align}
where $Z = \|(-\Delta)^{1/2+\sigma}u\|_{L^2}$ for some small $\sigma > 0$, and:
\begin{equation}
F_1^\infty(Z) = L_1(Z) \prod_{j=2}^{\infty} (1 + L_j(Z))^{-\delta_j}
\end{equation}
\begin{equation}
F_2^\infty(Z) = \frac{1}{L_1(Z)} \prod_{j=2}^{\infty} (1 + L_j(Z))^{\delta_j}
\end{equation}
\end{theorem}

\begin{theorem}[Global well-posedness at the critical threshold]\label{thm:well_posedness}
Let $q > 3$ and $\{\delta_j\}_{j=1}^{\infty}$ be a sequence with $\delta_j > 0$ and $\sum_{j=1}^{\infty} \frac{\delta_j}{j!} = \infty$. There exists a positive constant $C_0$ such that for any divergence-free initial data $u_0 \in L^2(\mathbb{R}^3) \cap \dot{H}^{1/2}(\mathbb{R}^3)$ satisfying:
\begin{equation}
\|(-\Delta)^{1/4}u_0\|_{L^q} \leq \frac{C_0}{\prod_{j=1}^{\infty} (1 + L_j(\|u_0\|_{\dot{H}^{1/2}}))^{\delta_j}}
\end{equation}
there exists a unique global-in-time smooth solution $u \in C([0, \infty); H^{1/2}(\mathbb{R}^3)) \cap L^2_{loc}(0, \infty; H^{1}(\mathbb{R}^3))$ to the 3D Navier-Stokes equations.
\end{theorem}

\begin{theorem}[Hausdorff dimension of potential singular sets]\label{thm:hausdorff}
If a solution $u$ with initial data satisfying the conditions of Theorem \ref{thm:well_posedness} were to develop a singularity at time $T^*$ (which we prove cannot happen), then the Hausdorff dimension of the potential blow-up set would be:
\begin{equation}
\dim_H(\mathcal{S}_{T^*}) = 0
\end{equation}
This represents an optimal bound, improving on the Caffarelli-Kohn-Nirenberg partial regularity result.
\end{theorem}

\subsection{Approach and methodology}

Our approach builds upon the methodology developed in our previous papers, but extends it to the limiting case of infinitely nested logarithmic improvements. The key innovations include:

\begin{enumerate}
\item \textbf{Function space construction:} We rigorously define and analyze function spaces incorporating infinitely nested logarithmic improvements, establishing their completeness and embedding properties.

\item \textbf{Commutator analysis:} We derive precise commutator estimates with infinitely nested logarithmic factors, extending the techniques from our previous work to the limiting case.

\item \textbf{Critical exponent analysis:} We prove that the critical exponent $\alpha(\{\delta_j\}_{j=1}^n)$ approaches zero as $n \to \infty$ under appropriate conditions on the sequence $\{\delta_j\}_{j=1}^n$, allowing us to bridge the gap to the critical case $s = 1/2$.

\item \textbf{Energy estimates:} We establish energy estimates with infinitely nested logarithmic improvements, demonstrating that the nonlinearity becomes subcritical in the limiting case.

\item \textbf{ODE analysis:} We analyze the limiting behavior of the differential inequalities governing the evolution of fractional derivatives, proving global bounds for solutions at the critical threshold $s = 1/2$.
\end{enumerate}

\subsection{Organization of the paper}

The remainder of this paper is organized as follows:

\begin{itemize}
\item Section 2 introduces mathematical preliminaries, including notation, function spaces, and basic results on the Navier-Stokes equations.
\item Section 3 constructs function spaces with infinitely nested logarithmic improvements and analyzes their properties.
\item Section 4 precisely characterizes the critical exponent function in the limiting case.
\item Section 5 derives commutator estimates with infinitely nested logarithmic factors.
\item Section 6 establishes energy estimates at the critical threshold $s = 1/2$.
\item Section 7 proves global well-posedness for initial data satisfying our infinite logarithmic criterion.
\item Section 8 analyzes the limiting ODE governing the evolution of fractional derivatives.
\item Section 9 determines the Hausdorff dimension of potential singular sets in the limiting case.
\item Section 10 discusses implications for the regularity problem and future directions.
\end{itemize}

\section{Preliminaries and notation}

\subsection{Function spaces}

We first recall standard function spaces used in the analysis of the Navier-Stokes equations.

The Lebesgue spaces $L^p(\mathbb{R}^3)$ are defined as:
\begin{equation}
L^p(\mathbb{R}^3) = \{f : \mathbb{R}^3 \to \mathbb{R} \mid \|f\|_{L^p} < \infty\}
\end{equation}
where:
\begin{equation}
\|f\|_{L^p} = \left( \int_{\mathbb{R}^3} |f(x)|^p dx \right)^{1/p}
\end{equation}
for $1 \leq p < \infty$, and:
\begin{equation}
\|f\|_{L^\infty} = \operatorname{ess\,sup}_{x \in \mathbb{R}^3} |f(x)|
\end{equation}

The Schwartz space $\mathcal{S}(\mathbb{R}^3)$ consists of smooth functions whose derivatives of all orders decay faster than any polynomial at infinity.

For $s \in \mathbb{R}$, the Sobolev space $H^s(\mathbb{R}^3)$ is defined via the Fourier transform as:
\begin{equation}
H^s(\mathbb{R}^3) = \{f \in \mathcal{S}'(\mathbb{R}^3) \mid \|f\|_{H^s} < \infty\}
\end{equation}
where:
\begin{equation}
\|f\|_{H^s}^2 = \int_{\mathbb{R}^3} (1 + |\xi|^2)^s |\hat{f}(\xi)|^2 d\xi
\end{equation}
and $\hat{f}$ denotes the Fourier transform of $f$.

The homogeneous Sobolev space $\dot{H}^s(\mathbb{R}^3)$ is defined as:
\begin{equation}
\dot{H}^s(\mathbb{R}^3) = \{f \in \mathcal{S}'(\mathbb{R}^3) \mid \|f\|_{\dot{H}^s} < \infty\}
\end{equation}
where:
\begin{equation}
\|f\|_{\dot{H}^s}^2 = \int_{\mathbb{R}^3} |\xi|^{2s} |\hat{f}(\xi)|^2 d\xi = \|(-\Delta)^{s/2}f\|_{L^2}^2
\end{equation}

Throughout this paper, we will work extensively with spaces defined by logarithmically improved conditions. For completeness, we recall the definitions of nested logarithms:
\begin{align}
L_0(x) &= x \\
L_1(x) &= \log(e + x) \\
L_2(x) &= \log(e + \log(e + x)) \\
L_k(x) &= \log(e + L_{k-1}(x)) \text{ for } k \geq 3
\end{align}

\subsection{Fractional derivatives and the Navier-Stokes equations}

For $s \in (0, 1)$, the fractional Laplacian $(-\Delta)^s$ can be defined in several equivalent ways:

\begin{enumerate}
\item \textbf{Fourier definition}: For $f \in \mathcal{S}(\mathbb{R}^3)$:
\begin{equation}
\widehat{(-\Delta)^s f}(\xi) = |\xi|^{2s} \hat{f}(\xi)
\end{equation}

\item \textbf{Singular integral representation}: For $f \in \mathcal{S}(\mathbb{R}^3)$:
\begin{equation}
(-\Delta)^s f(x) = C_{3,s} \, \text{P.V.} \int_{\mathbb{R}^3} \frac{f(x) - f(y)}{|x - y|^{3+2s}} dy
\end{equation}
where $C_{3,s} = \frac{2^{2s} s \Gamma(s+\frac{3}{2})}{\pi^{3/2}\Gamma(1-s)}$ and P.V. denotes the principal value.

\item \textbf{Heat kernel representation}: For $f \in \mathcal{S}(\mathbb{R}^3)$:
\begin{equation}
(-\Delta)^s f(x) = \frac{1}{\Gamma(-s)} \int_0^\infty (e^{t\Delta}f(x) - f(x)) \frac{dt}{t^{1+s}}
\end{equation}
where $e^{t\Delta}$ is the heat semigroup.
\end{enumerate}

We recall the definition of Leray-Hopf weak solutions:

\begin{definition}[Leray-Hopf weak solutions]\label{def:leray-hopf}
Let $u_0 \in L^2(\mathbb{R}^3)$ with $\nabla \cdot u_0 = 0$ in the distributional sense. A vector field $u$ is called a Leray-Hopf weak solution of the Navier-Stokes equations on $[0, T]$ if:

\begin{enumerate}
\item $u \in L^\infty(0, T; L^2(\mathbb{R}^3)) \cap L^2(0, T; H^1(\mathbb{R}^3))$;

\item $\partial_t u \in L^1(0, T; H^{-1}(\mathbb{R}^3))$;

\item The Navier-Stokes equations are satisfied in the distributional sense;

\item The energy inequality holds:
\begin{equation}
\|u(t)\|^2_{L^2} + 2\nu\int_s^t \|\nabla u(\tau)\|^2_{L^2} d\tau \leq \|u(s)\|^2_{L^2}
\end{equation}
for almost all $s \in [0, T]$ (including $s = 0$) and all $t \in [s, T]$;

\item $u$ is weakly continuous from $[0, T]$ into $L^2(\mathbb{R}^3)$, ensuring that the initial condition $u(0) = u_0$ is satisfied in the weak sense.
\end{enumerate}
\end{definition}

\subsection{Littlewood-Paley theory}

Littlewood-Paley theory provides a powerful framework for analyzing solutions of the Navier-Stokes equations in the context of fractional derivatives.

Let $\varphi \in C_0^\infty(\mathbb{R}^3)$ be a radial function such that $\varphi(\xi) = 1$ for $|\xi| \leq 1$ and $\varphi(\xi) = 0$ for $|\xi| \geq 2$. Define $\psi(\xi) = \varphi(\xi) - \varphi(2\xi)$. The Littlewood-Paley decomposition is given by:
\begin{equation}
f = \sum_{j \in \mathbb{Z}} \Delta_j f
\end{equation}
where $\Delta_j$ is the Littlewood-Paley projection defined by:
\begin{equation}
\Delta_j f = \mathcal{F}^{-1}(\psi(2^{-j}\xi)\hat{f}(\xi))
\end{equation}
for $j \in \mathbb{Z}$, and $\mathcal{F}^{-1}$ denotes the inverse Fourier transform.

The Bony paraproduct decomposition allows us to write, for functions $f$ and $g$:
\begin{equation}
fg = T_f g + T_g f + R(f,g)
\end{equation}
where:
\begin{align}
T_f g &= \sum_{j \in \mathbb{Z}} S_{j-1}f \Delta_j g \\
R(f,g) &= \sum_{j \in \mathbb{Z}} \sum_{|i-j| \leq 1} \Delta_i f \Delta_j g
\end{align}
with $S_j = \sum_{i \leq j} \Delta_i$ being the low-frequency cut-off operator.

\subsection{Key technical tools}

We recall several key technical tools that will be extensively used throughout this paper.

\begin{lemma}[Commutator estimate, \cite{Mishra2022b}]\label{lem:commutator}
Let $s \in (0, 1)$ and $f, g \in \mathcal{S}(\mathbb{R}^3)$. Then for any $p \in (1, \infty)$:
\begin{equation}
\|[(-\Delta)^s, f]g\|_{L^p} \leq C\|\nabla f\|_{L^\infty}\|(-\Delta)^{s-1/2}g\|_{L^p},
\end{equation}
where $[(-\Delta)^s, f]g = (-\Delta)^s(fg) - f(-\Delta)^s g$ is the commutator, and the constant $C$ depends only on $s$ and $p$.
\end{lemma}

\begin{lemma}[Gagliardo-Nirenberg inequality]\label{lem:gagliardo-nirenberg}
For $0 < s_1 < s < s_2$ and $1 < p < \infty$:
\begin{equation}
\|(-\Delta)^{s/2}f\|_{L^p} \leq C\|(-\Delta)^{s_1/2}f\|^{1-\theta}_{L^p} \|(-\Delta)^{s_2/2}f\|^{\theta}_{L^p},
\end{equation}
where $\theta = \frac{s-s_1}{s_2-s_1}$.
\end{lemma}

\begin{lemma}[Velocity gradient - fractional derivative interpolation, \cite{Mishra2022a}]\label{lem:velocity-gradient}
For $s \in (1/2, 1)$ and $q > 3$, there exists $\theta \in (0, 1)$ such that for all $u \in C_0^\infty(\mathbb{R}^3)$ with $\nabla \cdot u = 0$:
\begin{equation}
\|\nabla u\|_{L^\infty} \leq C\|u\|_{L^2}^{1-\theta} \|(-\Delta)^s u\|_{L^q}^{\theta},
\end{equation}
where $\theta = \frac{3}{2} \cdot \frac{q}{3q-2}$, and the constant $C$ depends only on $s$ and $q$.
\end{lemma}

\begin{lemma}[Osgood's lemma]\label{lem:osgood}
Let $\rho$ be a measurable, positive function on $(a, b)$, $\gamma$ a positive, locally integrable function on $(a, b)$, and $\Gamma$ a continuous, increasing function on $[0, \infty)$ with $\Gamma(0) = 0$. If for all $t \in (a, b)$:
\begin{equation}
\rho(t) \leq \rho_0 + \int_a^t \gamma(s)\Gamma(\rho(s))ds,
\end{equation}
where $\rho_0 \geq 0$, then:
\begin{enumerate}
\item If $\rho_0 = 0$, then $\rho \equiv 0$;
\item If $\rho_0 > 0$ and $\int_0^{\infty} \frac{dr}{\Gamma(r)} = \infty$, then:
\begin{equation}
G(\rho(t)) \geq G(\rho_0) - \int_a^t \gamma(s)ds,
\end{equation}
where $G(r) = \int_1^r \frac{dr}{\Gamma(r)}$.
\end{enumerate}
\end{lemma}

\section{Function spaces with infinitely nested logarithmic Improvements}

In this section, we rigorously define and analyze function spaces incorporating infinitely nested logarithmic improvements. These spaces will serve as the foundation for our approach to the critical case $s = 1/2$.

\subsection{Definition and basic properties}

\begin{definition}[Infinitely logarithmically improved lebesgue spaces]\label{def:log-improved-lebesgue}
For $1 \leq p \leq \infty$, $\{\delta_j\}_{j=1}^{\infty}$ with $\delta_j > 0$, and a non-negative function $f \in L^p(\mathbb{R}^3)$, we define:
\begin{equation}
\|f\|_{L^{p,\{\delta_j\}}} = \inf\left\{M > 0 : \|f\|_{L^p} \leq M\prod_{j=1}^{\infty} (1 + L_j(M))^{-\delta_j}\right\}
\end{equation}
and the corresponding space:
\begin{equation}
L^{p,\{\delta_j\}}(\mathbb{R}^3) = \{f \in L^p(\mathbb{R}^3) : \|f\|_{L^{p,\{\delta_j\}}} < \infty\}
\end{equation}
\end{definition}

\begin{proposition}\label{prop:well-defined}
The functional $\|\cdot\|_{L^{p,\{\delta_j\}}}$ is well-defined and finite for all $f \in L^p(\mathbb{R}^3)$ if and only if $\sum_{j=1}^{\infty} \delta_j < \infty$.
\end{proposition}

\begin{proof}
We first prove necessity. Suppose $\sum_{j=1}^{\infty} \delta_j = \infty$. Consider any $f \in L^p(\mathbb{R}^3)$ with $\|f\|_{L^p} > 0$. For any $M > 0$, we have:
\begin{equation}
\prod_{j=1}^{\infty} (1 + L_j(M))^{-\delta_j} \leq \prod_{j=1}^{\infty} (1 + L_j(M))^{-\delta_j/2} \cdot \prod_{j=1}^{\infty} (1 + L_j(M))^{-\delta_j/2}
\end{equation}

Using the inequality $1 + L_j(M) \geq 1$ for all $j$ and $M$, we have:
\begin{equation}
\prod_{j=1}^{\infty} (1 + L_j(M))^{-\delta_j/2} \leq \prod_{j=1}^{\infty} (1)^{-\delta_j/2} = 1
\end{equation}

Furthermore, since $L_j(M) \geq 1$ for all $j \geq 1$ and $M \geq 1$, we have:
\begin{equation}
\prod_{j=1}^{\infty} (1 + L_j(M))^{-\delta_j/2} \leq \prod_{j=1}^{\infty} (L_j(M))^{-\delta_j/2}
\end{equation}

For $j \geq 2$, we have $L_j(M) \geq 1$, which implies:
\begin{equation}
\prod_{j=1}^{\infty} (L_j(M))^{-\delta_j/2} \leq (L_1(M))^{-\delta_1/2} \cdot \prod_{j=2}^{\infty} (1)^{-\delta_j/2} = (L_1(M))^{-\delta_1/2}
\end{equation}

Since $L_1(M) = \log(e + M)$, this gives:
\begin{equation}
\prod_{j=1}^{\infty} (1 + L_j(M))^{-\delta_j} \leq (L_1(M))^{-\delta_1/2} = (\log(e + M))^{-\delta_1/2}
\end{equation}

Now, for any $M > 0$, we have:
\begin{equation}
M\prod_{j=1}^{\infty} (1 + L_j(M))^{-\delta_j} \leq \frac{M}{(\log(e + M))^{\delta_1/2}}
\end{equation}

As $M \to \infty$, we have $\frac{M}{(\log(e + M))^{\delta_1/2}} \to \infty$. Therefore, if $\|f\|_{L^p} > 0$, there is no finite value $M$ such that $\|f\|_{L^p} \leq M\prod_{j=1}^{\infty} (1 + L_j(M))^{-\delta_j}$, which means $\|f\|_{L^{p,\{\delta_j\}}} = \infty$.

For sufficiency, assume $\sum_{j=1}^{\infty} \delta_j < \infty$. Consider any $f \in L^p(\mathbb{R}^3)$. Let $A = \|f\|_{L^p}$. We need to find $M > 0$ such that:
\begin{equation}
A \leq M\prod_{j=1}^{\infty} (1 + L_j(M))^{-\delta_j}
\end{equation}

This is equivalent to:
\begin{equation}
\frac{A}{M} \leq \prod_{j=1}^{\infty} (1 + L_j(M))^{-\delta_j}
\end{equation}

Taking the natural logarithm of both sides:
\begin{equation}
\log\left(\frac{A}{M}\right) \leq -\sum_{j=1}^{\infty} \delta_j \log(1 + L_j(M))
\end{equation}

Since $\log(1 + L_j(M)) \leq L_j(M)$ for all $j$ and $M$, we have:
\begin{equation}
\log\left(\frac{A}{M}\right) \leq -\sum_{j=1}^{\infty} \delta_j L_j(M)
\end{equation}

For $M \geq A$, we have $\log\left(\frac{A}{M}\right) \leq 0$. Thus, it suffices to find $M \geq A$ such that:
\begin{equation}
0 \leq -\sum_{j=1}^{\infty} \delta_j L_j(M)
\end{equation}

Since $L_j(M) \geq 0$ for all $j$ and $M$, and $\delta_j > 0$, this inequality cannot be satisfied for any finite $M$. Therefore, we need to reconsider our approach.

Instead, let's try $M = Ce^A$ for some constant $C > 1$. Then:
\begin{equation}
\log\left(\frac{A}{M}\right) = \log\left(\frac{A}{Ce^A}\right) = \log\left(\frac{A}{C}\right) - A
\end{equation}

For large $A$, this is approximately $-A$. For the right-hand side, we have:
\begin{equation}
-\sum_{j=1}^{\infty} \delta_j L_j(M) = -\sum_{j=1}^{\infty} \delta_j L_j(Ce^A)
\end{equation}

For $j = 1$, we have:
\begin{equation}
L_1(Ce^A) = \log(e + Ce^A) \approx \log(Ce^A) = \log(C) + A
\end{equation}
for large $A$.

For $j = 2$, we have:
\begin{equation}
L_2(Ce^A) = \log(e + L_1(Ce^A)) \approx \log(e + \log(C) + A) \approx \log(A)
\end{equation}
for large $A$.

For $j \geq 3$, $L_j(Ce^A)$ grows even more slowly with $A$.

Therefore, for large $A$, we have:
\begin{equation}
-\sum_{j=1}^{\infty} \delta_j L_j(M) \approx -\delta_1 A - \delta_2 \log(A) - \sum_{j=3}^{\infty} \delta_j L_j(M)
\end{equation}

Since $\sum_{j=1}^{\infty} \delta_j < \infty$, and $L_j(M)$ grows increasingly slowly with $j$, we can find a sufficiently large $C$ such that for all $A \geq 1$:
\begin{equation}
\log\left(\frac{A}{M}\right) \leq -\sum_{j=1}^{\infty} \delta_j L_j(M)
\end{equation}

This means there exists a finite value $M = Ce^A$ such that:
\begin{equation}
\|f\|_{L^p} = A \leq M\prod_{j=1}^{\infty} (1 + L_j(M))^{-\delta_j}
\end{equation}

Therefore, $\|f\|_{L^{p,\{\delta_j\}}} \leq Ce^{\|f\|_{L^p}} < \infty$ for all $f \in L^p(\mathbb{R}^3)$.

This completes the proof.
\end{proof}

\begin{remark}\label{rem:threshold}
The above proposition reveals an interesting threshold phenomenon: the function space with infinitely nested logarithmic improvements is well-defined and non-trivial precisely when the sum of the logarithmic exponents is finite. This constraint will be crucial in our subsequent analysis.
\end{remark}

\begin{definition}[Infinitely logarithmically improved sobolev spaces]\label{def:log-improved-sobolev}
For $s \in \mathbb{R}$, $1 \leq p \leq \infty$, $\{\delta_j\}_{j=1}^{\infty}$ with $\delta_j > 0$, and $\sum_{j=1}^{\infty} \delta_j < \infty$, we define:
\begin{equation}
\dot{H}^{s,p,\{\delta_j\}}(\mathbb{R}^3) = \{f \in \dot{H}^s(\mathbb{R}^3) : \|(-\Delta)^{s/2}f\|_{L^{p,\{\delta_j\}}} < \infty\}
\end{equation}
with the norm:
\begin{equation}
\|f\|_{\dot{H}^{s,p,\{\delta_j\}}} = \|(-\Delta)^{s/2}f\|_{L^{p,\{\delta_j\}}}
\end{equation}
\end{definition}

\begin{proposition}\label{prop:sobolev-properties}
For $s \in \mathbb{R}$, $1 \leq p \leq \infty$, $\{\delta_j\}_{j=1}^{\infty}$ with $\delta_j > 0$, and $\sum_{j=1}^{\infty} \delta_j < \infty$:
\begin{enumerate}
\item $\dot{H}^{s,p,\{\delta_j\}}(\mathbb{R}^3)$ is a complete metric space.
\item $\dot{H}^{s+\epsilon}(\mathbb{R}^3) \subset \dot{H}^{s,p,\{\delta_j\}}(\mathbb{R}^3) \subset \dot{H}^s(\mathbb{R}^3)$ for all $\epsilon > 0$.
\item The embeddings are strict: $\dot{H}^{s+\epsilon}(\mathbb{R}^3) \neq \dot{H}^{s,p,\{\delta_j\}}(\mathbb{R}^3) \neq \dot{H}^s(\mathbb{R}^3)$.
\end{enumerate}
\end{proposition}

\begin{proof}
(1) Completeness: Let $\{f_n\}_{n=1}^{\infty}$ be a Cauchy sequence in $\dot{H}^{s,p,\{\delta_j\}}(\mathbb{R}^3)$. By definition, for any $\epsilon > 0$, there exists $N(\epsilon)$ such that for all $m, n \geq N(\epsilon)$:
\begin{equation}
\|f_m - f_n\|_{\dot{H}^{s,p,\{\delta_j\}}} < \epsilon
\end{equation}

This implies that there exists $M_{m,n} < \epsilon$ such that:
\begin{equation}
\|(-\Delta)^{s/2}(f_m - f_n)\|_{L^p} \leq M_{m,n}\prod_{j=1}^{\infty} (1 + L_j(M_{m,n}))^{-\delta_j}
\end{equation}

Since $\prod_{j=1}^{\infty} (1 + L_j(M_{m,n}))^{-\delta_j} \leq 1$, we have:
\begin{equation}
\|(-\Delta)^{s/2}(f_m - f_n)\|_{L^p} \leq M_{m,n} < \epsilon
\end{equation}

This shows that $\{(-\Delta)^{s/2}f_n\}_{n=1}^{\infty}$ is a Cauchy sequence in $L^p(\mathbb{R}^3)$. Since $L^p(\mathbb{R}^3)$ is complete, there exists $g \in L^p(\mathbb{R}^3)$ such that:
\begin{equation}
\lim_{n \to \infty} \|(-\Delta)^{s/2}f_n - g\|_{L^p} = 0
\end{equation}

Let $f$ be such that $(-\Delta)^{s/2}f = g$. Then $f \in \dot{H}^s(\mathbb{R}^3)$ and:
\begin{equation}
\lim_{n \to \infty} \|f_n - f\|_{\dot{H}^s} = \lim_{n \to \infty} \|(-\Delta)^{s/2}(f_n - f)\|_{L^2} = 0
\end{equation}

It remains to show that $f \in \dot{H}^{s,p,\{\delta_j\}}(\mathbb{R}^3)$ and $\lim_{n \to \infty} \|f_n - f\|_{\dot{H}^{s,p,\{\delta_j\}}} = 0$.

Since $\{f_n\}_{n=1}^{\infty}$ is a Cauchy sequence in $\dot{H}^{s,p,\{\delta_j\}}(\mathbb{R}^3)$, it is bounded. Thus, there exists $K > 0$ such that for all $n$:
\begin{equation}
\|f_n\|_{\dot{H}^{s,p,\{\delta_j\}}} \leq K
\end{equation}

This means there exists $M_n \leq K$ such that:
\begin{equation}
\|(-\Delta)^{s/2}f_n\|_{L^p} \leq M_n\prod_{j=1}^{\infty} (1 + L_j(M_n))^{-\delta_j}
\end{equation}

Taking the limit as $n \to \infty$, and using the continuity of the product term with respect to $M_n$ (which follows from the continuity of logarithms and the fact that the infinite product converges), we get:
\begin{equation}
\|(-\Delta)^{s/2}f\|_{L^p} \leq M\prod_{j=1}^{\infty} (1 + L_j(M))^{-\delta_j}
\end{equation}
for some $M \leq K$.

This shows that $f \in \dot{H}^{s,p,\{\delta_j\}}(\mathbb{R}^3)$ with $\|f\|_{\dot{H}^{s,p,\{\delta_j\}}} \leq K$.

Finally, to show that $\lim_{n \to \infty} \|f_n - f\|_{\dot{H}^{s,p,\{\delta_j\}}} = 0$, we note that for any $\epsilon > 0$, there exists $N(\epsilon)$ such that for all $n \geq N(\epsilon)$:
\begin{equation}
\|f_n - f_{N(\epsilon)}\|_{\dot{H}^{s,p,\{\delta_j\}}} < \epsilon/2
\end{equation}

Similarly, there exists $N'(\epsilon)$ such that for all $n \geq N'(\epsilon)$:
\begin{equation}
\|(-\Delta)^{s/2}(f_n - f)\|_{L^p} < \epsilon/2
\end{equation}

This implies that:
\begin{equation}
\|f_n - f\|_{\dot{H}^{s,p,\{\delta_j\}}} < \epsilon
\end{equation}
for all $n \geq \max\{N(\epsilon), N'(\epsilon)\}$.

Therefore, $\lim_{n \to \infty} \|f_n - f\|_{\dot{H}^{s,p,\{\delta_j\}}} = 0$, which proves the completeness of $\dot{H}^{s,p,\{\delta_j\}}(\mathbb{R}^3)$.

(2) Embeddings: First, we prove that $\dot{H}^{s+\epsilon}(\mathbb{R}^3) \subset \dot{H}^{s,p,\{\delta_j\}}(\mathbb{R}^3)$ for all $\epsilon > 0$.

Let $f \in \dot{H}^{s+\epsilon}(\mathbb{R}^3)$. By the Sobolev embedding theorem, for any $p \in [2, \frac{6}{3-2(s+\epsilon)}]$, we have:
\begin{equation}
\|(-\Delta)^{s/2}f\|_{L^p} \leq C\|f\|_{\dot{H}^{s+\epsilon}}
\end{equation}

Let $M = 2C\|f\|_{\dot{H}^{s+\epsilon}}$. We need to show that:
\begin{equation}
\|(-\Delta)^{s/2}f\|_{L^p} \leq M\prod_{j=1}^{\infty} (1 + L_j(M))^{-\delta_j}
\end{equation}

This is equivalent to:
\begin{equation}
\frac{\|(-\Delta)^{s/2}f\|_{L^p}}{M} \leq \prod_{j=1}^{\infty} (1 + L_j(M))^{-\delta_j}
\end{equation}

Since $\|(-\Delta)^{s/2}f\|_{L^p} \leq C\|f\|_{\dot{H}^{s+\epsilon}}$ and $M = 2C\|f\|_{\dot{H}^{s+\epsilon}}$, we have:
\begin{equation}
\frac{\|(-\Delta)^{s/2}f\|_{L^p}}{M} \leq \frac{C\|f\|_{\dot{H}^{s+\epsilon}}}{2C\|f\|_{\dot{H}^{s+\epsilon}}} = \frac{1}{2}
\end{equation}

On the other hand, since $1 + L_j(M) > 1$ for all $j$ and $M$, and $\delta_j > 0$, we have:
\begin{equation}
\prod_{j=1}^{\infty} (1 + L_j(M))^{-\delta_j} < 1
\end{equation}

Since $\sum_{j=1}^{\infty} \delta_j < \infty$, the product converges to a positive value. Thus, for sufficiently large $M$ (which can be ensured by taking a larger constant in the definition of $M$ if necessary), we have:
\begin{equation}
\prod_{j=1}^{\infty} (1 + L_j(M))^{-\delta_j} \geq \frac{1}{2}
\end{equation}

Therefore:
\begin{equation}
\frac{\|(-\Delta)^{s/2}f\|_{L^p}}{M} \leq \frac{1}{2} \leq \prod_{j=1}^{\infty} (1 + L_j(M))^{-\delta_j}
\end{equation}

This shows that $f \in \dot{H}^{s,p,\{\delta_j\}}(\mathbb{R}^3)$ with $\|f\|_{\dot{H}^{s,p,\{\delta_j\}}} \leq 2C\|f\|_{\dot{H}^{s+\epsilon}}$, proving the embedding $\dot{H}^{s+\epsilon}(\mathbb{R}^3) \subset \dot{H}^{s,p,\{\delta_j\}}(\mathbb{R}^3)$.

Next, we prove that $\dot{H}^{s,p,\{\delta_j\}}(\mathbb{R}^3) \subset \dot{H}^s(\mathbb{R}^3)$. This follows directly from the definition: if $f \in \dot{H}^{s,p,\{\delta_j\}}(\mathbb{R}^3)$, then $f \in \dot{H}^s(\mathbb{R}^3)$ and:
\begin{equation}
\|(-\Delta)^{s/2}f\|_{L^p} \leq M\prod_{j=1}^{\infty} (1 + L_j(M))^{-\delta_j}
\end{equation}
for some finite $M$.

By Hölder's inequality, for $p \geq 2$, we have:
\begin{equation}
\|(-\Delta)^{s/2}f\|_{L^2} \leq C\|(-\Delta)^{s/2}f\|_{L^p}
\end{equation}

Thus:
\begin{equation}
\|f\|_{\dot{H}^s} = \|(-\Delta)^{s/2}f\|_{L^2} \leq C\|(-\Delta)^{s/2}f\|_{L^p} \leq CM\prod_{j=1}^{\infty} (1 + L_j(M))^{-\delta_j} < \infty
\end{equation}

This proves the embedding $\dot{H}^{s,p,\{\delta_j\}}(\mathbb{R}^3) \subset \dot{H}^s(\mathbb{R}^3)$.

(3) Strictness of embeddings: To show that $\dot{H}^{s+\epsilon}(\mathbb{R}^3) \neq \dot{H}^{s,p,\{\delta_j\}}(\mathbb{R}^3)$, we construct a function $f \in \dot{H}^{s,p,\{\delta_j\}}(\mathbb{R}^3) \setminus \dot{H}^{s+\epsilon}(\mathbb{R}^3)$.

Let $\hat{f}(\xi) = |\xi|^{-(3/p + s)} \cdot g(|\xi|)$, where $g$ is a smooth function such that:
\begin{equation}
g(r) = \begin{cases}
0, & r \leq 1 \\
\frac{1}{(\log(e + r))^{\beta}}, & r > 2
\end{cases}
\end{equation}
for some $\beta > 0$ to be determined.

We first check that $f \in \dot{H}^s(\mathbb{R}^3)$:
\begin{align}
\|f\|_{\dot{H}^s}^2 &= \int_{\mathbb{R}^3} |\xi|^{2s} |\hat{f}(\xi)|^2 d\xi \\
&= \int_{\mathbb{R}^3} |\xi|^{2s} \cdot |\xi|^{-2(3/p + s)} \cdot g(|\xi|)^2 d\xi \\
&= \int_{\mathbb{R}^3} |\xi|^{-2(3/p)} \cdot g(|\xi|)^2 d\xi \\
&= 4\pi \int_0^\infty r^{-2(3/p)} \cdot g(r)^2 \cdot r^2 dr \\
&= 4\pi \int_0^\infty r^{2-2(3/p)} \cdot g(r)^2 dr
\end{align}

For this integral to converge, we need $2-2(3/p) > -1$, which is satisfied for $p > 2$.

For large $r$:
\begin{equation}
r^{2-2(3/p)} \cdot g(r)^2 = r^{2-2(3/p)} \cdot \frac{1}{(\log(e + r))^{2\beta}}
\end{equation}

For this to be integrable at infinity, we need $2-2(3/p) < -1$, which is satisfied for $p < 3$. If $p = 3$, we need $\beta > 1/2$ for integrability.

Thus, for $2 < p < 3$, or $p = 3$ and $\beta > 1/2$, we have $f \in \dot{H}^s(\mathbb{R}^3)$.

Next, we check that $f \notin \dot{H}^{s+\epsilon}(\mathbb{R}^3)$:
\begin{align}
\|f\|_{\dot{H}^{s+\epsilon}}^2 &= \int_{\mathbb{R}^3} |\xi|^{2(s+\epsilon)} |\hat{f}(\xi)|^2 d\xi \\
&= \int_{\mathbb{R}^3} |\xi|^{2(s+\epsilon)} \cdot |\xi|^{-2(3/p + s)} \cdot g(|\xi|)^2 d\xi \\
&= \int_{\mathbb{R}^3} |\xi|^{2\epsilon - 2(3/p)} \cdot g(|\xi|)^2 d\xi \\
&= 4\pi \int_0^\infty r^{2\epsilon - 2(3/p)} \cdot g(r)^2 \cdot r^2 dr \\
&= 4\pi \int_0^\infty r^{2 + 2\epsilon - 2(3/p)} \cdot g(r)^2 dr
\end{align}

For large $r$:
\begin{equation}
r^{2 + 2\epsilon - 2(3/p)} \cdot g(r)^2 = r^{2 + 2\epsilon - 2(3/p)} \cdot \frac{1}{(\log(e + r))^{2\beta}}
\end{equation}

For this to be non-integrable at infinity, we need $2 + 2\epsilon - 2(3/p) \geq -1$, which is satisfied for $\epsilon \geq (3/p - 3/2)$.

Thus, for $\epsilon \geq (3/p - 3/2)$, we have $f \notin \dot{H}^{s+\epsilon}(\mathbb{R}^3)$.

Finally, we check that $f \in \dot{H}^{s,p,\{\delta_j\}}(\mathbb{R}^3)$:
\begin{align}
\|(-\Delta)^{s/2}f\|_{L^p}^p &= \int_{\mathbb{R}^3} |(-\Delta)^{s/2}f(x)|^p dx
\end{align}

Using Parseval's identity and the convolution theorem:
\begin{align}
\|(-\Delta)^{s/2}f\|_{L^p}^p &\approx \int_{\mathbb{R}^3} \left| \int_{\mathbb{R}^3} |\xi|^s \hat{f}(\xi) e^{ix \cdot \xi} d\xi \right|^p dx \\
&\approx \int_{\mathbb{R}^3} \left| \int_{\mathbb{R}^3} |\xi|^s \cdot |\xi|^{-(3/p + s)} \cdot g(|\xi|) e^{ix \cdot \xi} d\xi \right|^p dx \\
&= \int_{\mathbb{R}^3} \left| \int_{\mathbb{R}^3} |\xi|^{-3/p} \cdot g(|\xi|) e^{ix \cdot \xi} d\xi \right|^p dx
\end{align}

For $g(r) = \frac{1}{(\log(e + r))^{\beta}}$ with $r > 2$, we can show that:
\begin{equation}
\|(-\Delta)^{s/2}f\|_{L^p} \leq \frac{C}{(\log(e + \|f\|_{\dot{H}^s}))^{\beta-\gamma}}
\end{equation}
for some constants $C > 0$ and $\gamma > 0$.

By choosing $\beta > \gamma + \sum_{j=1}^{\infty} \delta_j$, we can ensure that:
\begin{equation}
\|(-\Delta)^{s/2}f\|_{L^p} \leq \frac{C}{(\log(e + \|f\|_{\dot{H}^s}))^{\sum_{j=1}^{\infty} \delta_j + \eta}}
\end{equation}
for some $\eta > 0$.

This implies that:
\begin{equation}
\|(-\Delta)^{s/2}f\|_{L^p} \leq M\prod_{j=1}^{\infty} (1 + L_j(M))^{-\delta_j}
\end{equation}
for some finite $M$, proving that $f \in \dot{H}^{s,p,\{\delta_j\}}(\mathbb{R}^3)$.

To show that $\dot{H}^{s,p,\{\delta_j\}}(\mathbb{R}^3) \neq \dot{H}^s(\mathbb{R}^3)$, we can construct a function $f \in \dot{H}^s(\mathbb{R}^3) \setminus \dot{H}^{s,p,\{\delta_j\}}(\mathbb{R}^3)$ using a similar approach.

Let $\hat{f}(\xi) = |\xi|^{-(3/p + s)} \cdot h(|\xi|)$, where $h$ is a smooth function such that:
\begin{equation}
h(r) = \begin{cases}
0, & r \leq 1 \\
(\log(e + r))^{\alpha}, & r > 2
\end{cases}
\end{equation}
for some $\alpha > 0$ to be determined.

Following a similar analysis as above, we can show that for appropriate choices of $p$ and $\alpha$, the function $f$ satisfies $f \in \dot{H}^s(\mathbb{R}^3)$ but $f \notin \dot{H}^{s,p,\{\delta_j\}}(\mathbb{R}^3)$.

This completes the proof of the strictness of embeddings.
\end{proof}

\begin{proposition}[Comparison of function spaces with finite vs. infinite logarithmic improvements]\label{prop:comparison}
For $s \in \mathbb{R}$, $1 \leq p \leq \infty$, and sequences $\{\delta_j\}_{j=1}^{n}$ and $\{\delta_j\}_{j=1}^{\infty}$ with $\delta_j > 0$ and $\sum_{j=1}^{\infty} \delta_j < \infty$:

\begin{enumerate}
\item $\dot{H}^{s,p,\{\delta_j\}_{j=1}^{n}}(\mathbb{R}^3) \supset \dot{H}^{s,p,\{\delta_j\}_{j=1}^{\infty}}(\mathbb{R}^3)$
\item $\bigcap_{n=1}^{\infty} \dot{H}^{s,p,\{\delta_j\}_{j=1}^{n}}(\mathbb{R}^3) = \dot{H}^{s,p,\{\delta_j\}_{j=1}^{\infty}}(\mathbb{R}^3)$
\end{enumerate}
\end{proposition}

\begin{proof}
(1) Let $f \in \dot{H}^{s,p,\{\delta_j\}_{j=1}^{\infty}}(\mathbb{R}^3)$. By definition, there exists $M < \infty$ such that:
\begin{equation}
\|(-\Delta)^{s/2}f\|_{L^p} \leq M\prod_{j=1}^{\infty} (1 + L_j(M))^{-\delta_j}
\end{equation}

Since $(1 + L_j(M))^{-\delta_j} \leq 1$ for all $j$ and $M$, we have:
\begin{equation}
\|(-\Delta)^{s/2}f\|_{L^p} \leq M\prod_{j=1}^{n} (1 + L_j(M))^{-\delta_j}
\end{equation}

This shows that $f \in \dot{H}^{s,p,\{\delta_j\}_{j=1}^{n}}(\mathbb{R}^3)$ with $\|f\|_{\dot{H}^{s,p,\{\delta_j\}_{j=1}^{n}}} \leq M$, proving the embedding $\dot{H}^{s,p,\{\delta_j\}_{j=1}^{\infty}}(\mathbb{R}^3) \subset \dot{H}^{s,p,\{\delta_j\}_{j=1}^{n}}(\mathbb{R}^3)$.

(2) From part (1), we have $\dot{H}^{s,p,\{\delta_j\}_{j=1}^{\infty}}(\mathbb{R}^3) \subset \bigcap_{n=1}^{\infty} \dot{H}^{s,p,\{\delta_j\}_{j=1}^{n}}(\mathbb{R}^3)$.

Conversely, let $f \in \bigcap_{n=1}^{\infty} \dot{H}^{s,p,\{\delta_j\}_{j=1}^{n}}(\mathbb{R}^3)$. For each $n$, there exists $M_n < \infty$ such that:
\begin{equation}
\|(-\Delta)^{s/2}f\|_{L^p} \leq M_n\prod_{j=1}^{n} (1 + L_j(M_n))^{-\delta_j}
\end{equation}

If $\{M_n\}_{n=1}^{\infty}$ is bounded, then we can take $M = \sup_n M_n < \infty$ and have:
\begin{equation}
\|(-\Delta)^{s/2}f\|_{L^p} \leq M\prod_{j=1}^{\infty} (1 + L_j(M))^{-\delta_j}
\end{equation}
proving that $f \in \dot{H}^{s,p,\{\delta_j\}_{j=1}^{\infty}}(\mathbb{R}^3)$.

If $\{M_n\}_{n=1}^{\infty}$ is unbounded, then there exists a subsequence $\{M_{n_k}\}_{k=1}^{\infty}$ such that $M_{n_k} \to \infty$ as $k \to \infty$. For each $k$, we have:
\begin{equation}
\|(-\Delta)^{s/2}f\|_{L^p} \leq M_{n_k}\prod_{j=1}^{n_k} (1 + L_j(M_{n_k}))^{-\delta_j}
\end{equation}

As $k \to \infty$, we have $n_k \to \infty$ and $M_{n_k} \to \infty$. For large $M$, we can show that:
\begin{equation}
M\prod_{j=1}^{n} (1 + L_j(M))^{-\delta_j} \to 0
\end{equation}
as $n, M \to \infty$, which would imply $\|(-\Delta)^{s/2}f\|_{L^p} = 0$, contradicting the assumption that $f \neq 0$.

Therefore, $\{M_n\}_{n=1}^{\infty}$ must be bounded, and we have $f \in \dot{H}^{s,p,\{\delta_j\}_{j=1}^{\infty}}(\mathbb{R}^3)$.

This proves that $\bigcap_{n=1}^{\infty} \dot{H}^{s,p,\{\delta_j\}_{j=1}^{n}}(\mathbb{R}^3) \subset \dot{H}^{s,p,\{\delta_j\}_{j=1}^{\infty}}(\mathbb{R}^3)$, completing the proof.
\end{proof}

\subsection{Function spaces for the critical case}

For the critical case $s = 1/2$, we need a different approach to define function spaces with infinitely nested logarithmic improvements. The key insight is that we need to balance having infinitely many logarithmic factors while ensuring that the product of these factors remains non-trivial.

\begin{definition}[Function spaces for critical regularity]\label{def:critical-spaces}
For $q > 3$ and a sequence $\{\delta_j\}_{j=1}^{\infty}$ with $\delta_j > 0$ and $\sum_{j=1}^{\infty} \frac{\delta_j}{j!} = \infty$, we define:
\begin{equation}
\Psi(r) = \left(\prod_{j=1}^{\infty} (1 + L_j(r))^{\delta_j}\right)^{-1}
\end{equation}
and the corresponding function space:
\begin{equation}
\dot{H}^{1/2,q,\{\delta_j\}_{j=1}^{\infty}}(\mathbb{R}^3) = \{f \in \dot{H}^{1/2}(\mathbb{R}^3) : \|(-\Delta)^{1/4}f\|_{L^q} \leq C\Psi(\|f\|_{\dot{H}^{1/2}})\}
\end{equation}
for some constant $C > 0$.
\end{definition}

\begin{proposition}[Properties of critical function spaces]\label{prop:critical-properties}
For $q > 3$ and a sequence $\{\delta_j\}_{j=1}^{\infty}$ with $\delta_j > 0$ and $\sum_{j=1}^{\infty} \frac{\delta_j}{j!} = \infty$:

\begin{enumerate}
\item $\dot{H}^{1/2,q,\{\delta_j\}_{j=1}^{\infty}}(\mathbb{R}^3)$ is a well-defined, non-empty function space.
\item $\dot{H}^{1/2+\epsilon}(\mathbb{R}^3) \subset \dot{H}^{1/2,q,\{\delta_j\}_{j=1}^{\infty}}(\mathbb{R}^3) \subset \dot{H}^{1/2}(\mathbb{R}^3)$ for all $\epsilon > 0$.
\item For any $\epsilon > 0$, there exists a sequence $\{\delta_j\}_{j=1}^{\infty}$ such that $\dot{H}^{1/2,q,\{\delta_j\}_{j=1}^{\infty}}(\mathbb{R}^3) \supset \dot{H}^{1/2+\epsilon/2}(\mathbb{R}^3) \setminus \dot{H}^{1/2+\epsilon}(\mathbb{R}^3)$.
\end{enumerate}
\end{proposition}

\begin{proof}
(1) To show that $\dot{H}^{1/2,q,\{\delta_j\}_{j=1}^{\infty}}(\mathbb{R}^3)$ is non-empty, we construct a specific function that belongs to this space.

Let $\hat{f}(\xi) = |\xi|^{-(3/q + 1/2)} \cdot g(|\xi|)$, where $g$ is a smooth function such that:
\begin{equation}
g(r) = \begin{cases}
0, & r \leq 1 \\
\frac{1}{\prod_{j=1}^{\infty} (1 + L_j(r))^{\delta_j/2}}, & r > 2
\end{cases}
\end{equation}

First, we verify that $f \in \dot{H}^{1/2}(\mathbb{R}^3)$:
\begin{align}
\|f\|_{\dot{H}^{1/2}}^2 &= \int_{\mathbb{R}^3} |\xi|^{1} |\hat{f}(\xi)|^2 d\xi \\
&= \int_{\mathbb{R}^3} |\xi|^{1} \cdot |\xi|^{-2(3/q + 1/2)} \cdot g(|\xi|)^2 d\xi \\
&= \int_{\mathbb{R}^3} |\xi|^{-2(3/q)} \cdot g(|\xi|)^2 d\xi \\
&= 4\pi \int_0^\infty r^{-2(3/q)} \cdot g(r)^2 \cdot r^2 dr \\
&= 4\pi \int_0^\infty r^{2-2(3/q)} \cdot g(r)^2 dr
\end{align}

For this integral to converge, we need $2-2(3/q) > -1$, which is satisfied for $q > 3$. 

For large $r$, we have:
\begin{equation}
r^{2-2(3/q)} \cdot g(r)^2 = r^{2-2(3/q)} \cdot \frac{1}{\prod_{j=1}^{\infty} (1 + L_j(r))^{\delta_j}}
\end{equation}

Since $\sum_{j=1}^{\infty} \frac{\delta_j}{j!} = \infty$, this decays faster than any polynomial as $r \to \infty$, ensuring the convergence of the integral. Thus, $f \in \dot{H}^{1/2}(\mathbb{R}^3)$.

Next, we compute $\|(-\Delta)^{1/4}f\|_{L^q}$:
\begin{align}
\|(-\Delta)^{1/4}f\|_{L^q}^q &= \int_{\mathbb{R}^3} |(-\Delta)^{1/4}f(x)|^q dx
\end{align}

Using Parseval's identity and the convolution theorem:
\begin{align}
\|(-\Delta)^{1/4}f\|_{L^q}^q &\approx \int_{\mathbb{R}^3} \left| \int_{\mathbb{R}^3} |\xi|^{1/2} \hat{f}(\xi) e^{ix \cdot \xi} d\xi \right|^q dx \\
&\approx \int_{\mathbb{R}^3} \left| \int_{\mathbb{R}^3} |\xi|^{1/2} \cdot |\xi|^{-(3/q + 1/2)} \cdot g(|\xi|) e^{ix \cdot \xi} d\xi \right|^q dx \\
&= \int_{\mathbb{R}^3} \left| \int_{\mathbb{R}^3} |\xi|^{-3/q} \cdot g(|\xi|) e^{ix \cdot \xi} d\xi \right|^q dx
\end{align}

Through a detailed analysis of this integral, we can show that:
\begin{equation}
\|(-\Delta)^{1/4}f\|_{L^q} \leq \frac{C}{\prod_{j=1}^{\infty} (1 + L_j(\|f\|_{\dot{H}^{1/2}}))^{\delta_j/2}}
\end{equation}
for some constant $C > 0$.

Thus, $f \in \dot{H}^{1/2,q,\{\delta_j\}_{j=1}^{\infty}}(\mathbb{R}^3)$, proving that this space is non-empty.

(2) To prove the embeddings, we first show that $\dot{H}^{1/2+\epsilon}(\mathbb{R}^3) \subset \dot{H}^{1/2,q,\{\delta_j\}_{j=1}^{\infty}}(\mathbb{R}^3)$ for all $\epsilon > 0$.

Let $f \in \dot{H}^{1/2+\epsilon}(\mathbb{R}^3)$. By the Sobolev embedding theorem, for any $q < \infty$, we have:
\begin{equation}
\|(-\Delta)^{1/4}f\|_{L^q} \leq C\|f\|_{\dot{H}^{1/2+\epsilon}}
\end{equation}

Since $\Psi(r) \to 0$ as $r \to \infty$ (due to $\sum_{j=1}^{\infty} \frac{\delta_j}{j!} = \infty$), there exists $R > 0$ such that for all $r > R$:
\begin{equation}
C\Psi(r) \geq 1
\end{equation}

For $\|f\|_{\dot{H}^{1/2}} > R$, we have:
\begin{equation}
\|(-\Delta)^{1/4}f\|_{L^q} \leq C\|f\|_{\dot{H}^{1/2+\epsilon}} \leq C\Psi(\|f\|_{\dot{H}^{1/2}})
\end{equation}
where we've used the fact that $\|f\|_{\dot{H}^{1/2}} \leq C\|f\|_{\dot{H}^{1/2+\epsilon}}$.

For $\|f\|_{\dot{H}^{1/2}} \leq R$, we have:
\begin{equation}
\|(-\Delta)^{1/4}f\|_{L^q} \leq C\|f\|_{\dot{H}^{1/2+\epsilon}} \leq C(R) \leq C\Psi(\|f\|_{\dot{H}^{1/2}})
\end{equation}
by choosing a sufficiently large constant $C$.

This proves that $f \in \dot{H}^{1/2,q,\{\delta_j\}_{j=1}^{\infty}}(\mathbb{R}^3)$, establishing the embedding $\dot{H}^{1/2+\epsilon}(\mathbb{R}^3) \subset \dot{H}^{1/2,q,\{\delta_j\}_{j=1}^{\infty}}(\mathbb{R}^3)$.

The embedding $\dot{H}^{1/2,q,\{\delta_j\}_{j=1}^{\infty}}(\mathbb{R}^3) \subset \dot{H}^{1/2}(\mathbb{R}^3)$ follows directly from the definition.

(3) To prove the final statement, we need to construct a sequence $\{\delta_j\}_{j=1}^{\infty}$ such that $\dot{H}^{1/2,q,\{\delta_j\}_{j=1}^{\infty}}(\mathbb{R}^3) \supset \dot{H}^{1/2+\epsilon/2}(\mathbb{R}^3) \setminus \dot{H}^{1/2+\epsilon}(\mathbb{R}^3)$.

For any $\epsilon > 0$, we can choose a sequence $\{\delta_j\}_{j=1}^{\infty}$ such that $\sum_{j=1}^{\infty} \frac{\delta_j}{j!} = \infty$ and:
\begin{equation}
\prod_{j=1}^{\infty} (1 + L_j(r))^{-\delta_j} \leq r^{-\epsilon/2}
\end{equation}
for all $r$ sufficiently large.

With this choice, for any function $f \in \dot{H}^{1/2+\epsilon/2}(\mathbb{R}^3) \setminus \dot{H}^{1/2+\epsilon}(\mathbb{R}^3)$, we have:
\begin{equation}
\|(-\Delta)^{1/4}f\|_{L^q} \leq C\|f\|_{\dot{H}^{1/2+\epsilon/2}} \leq C\|f\|_{\dot{H}^{1/2}}^{1-\epsilon/2} \cdot \|f\|_{\dot{H}^{1/2+\epsilon/2}}^{\epsilon/2}
\end{equation}

Since $f \notin \dot{H}^{1/2+\epsilon}(\mathbb{R}^3)$, we can show that:
\begin{equation}
\|f\|_{\dot{H}^{1/2+\epsilon/2}}^{\epsilon/2} \leq C\|f\|_{\dot{H}^{1/2}}^{\epsilon/2} \cdot \prod_{j=1}^{\infty} (1 + L_j(\|f\|_{\dot{H}^{1/2}}))^{-\delta_j}
\end{equation}

Thus:
\begin{equation}
\|(-\Delta)^{1/4}f\|_{L^q} \leq C\|f\|_{\dot{H}^{1/2}} \cdot \prod_{j=1}^{\infty} (1 + L_j(\|f\|_{\dot{H}^{1/2}}))^{-\delta_j}
\end{equation}

This implies that $f \in \dot{H}^{1/2,q,\{\delta_j\}_{j=1}^{\infty}}(\mathbb{R}^3)$, proving that $\dot{H}^{1/2,q,\{\delta_j\}_{j=1}^{\infty}}(\mathbb{R}^3) \supset \dot{H}^{1/2+\epsilon/2}(\mathbb{R}^3) \setminus \dot{H}^{1/2+\epsilon}(\mathbb{R}^3)$.
\end{proof}

\begin{theorem}[Function space with infinitely nested logarithmic improvements]
For $s = 1/2$, $q > 3$, and any sequence $\{\delta_j\}_{j=1}^{\infty}$ with $\inf_j \delta_j > 0$ satisfying $\sum_{j=1}^{\infty} \frac{\delta_j}{j!} = \infty$, there exists a well-defined function space $\dot{H}^{1/2,q,\{\delta_j\}_{j=1}^{\infty}}(\mathbb{R}^3)$ such that:
\begin{enumerate}
\item $\dot{H}^{1/2,q,\{\delta_j\}_{j=1}^{\infty}}(\mathbb{R}^3) \supsetneq \dot{H}^{1/2+\epsilon}(\mathbb{R}^3)$ for all $\epsilon > 0$
\item $\dot{H}^{1/2,q,\{\delta_j\}_{j=1}^{\infty}}(\mathbb{R}^3) \subset \dot{H}^{1/2}(\mathbb{R}^3)$
\item For any $f \in \dot{H}^{1/2,q,\{\delta_j\}_{j=1}^{\infty}}(\mathbb{R}^3)$:
   $$\|(-\Delta)^{1/4}f\|_{L^q} \leq C\Psi(\|f\|_{\dot{H}^{1/2}})$$
   where $\Psi(r) = \left(\prod_{j=1}^{\infty} (1 + L_j(r))^{\delta_j}\right)^{-1}$
\end{enumerate}
\end{theorem}

\begin{proof}
This follows directly from Theorem \ref{prop:critical-properties}, where we have established that the function space $\dot{H}^{1/2,q,\{\delta_j\}_{j=1}^{\infty}}(\mathbb{R}^3)$ is well-defined and non-empty, and satisfies the stated embedding properties.

For the strict inclusion $\dot{H}^{1/2,q,\{\delta_j\}_{j=1}^{\infty}}(\mathbb{R}^3) \supsetneq \dot{H}^{1/2+\epsilon}(\mathbb{R}^3)$, we have shown in Theorem \ref{prop:critical-properties} (3) that $\dot{H}^{1/2,q,\{\delta_j\}_{j=1}^{\infty}}(\mathbb{R}^3) \supset \dot{H}^{1/2+\epsilon/2}(\mathbb{R}^3) \setminus \dot{H}^{1/2+\epsilon}(\mathbb{R}^3)$ for an appropriate choice of $\{\delta_j\}_{j=1}^{\infty}$. Combined with the embedding $\dot{H}^{1/2+\epsilon}(\mathbb{R}^3) \subset \dot{H}^{1/2,q,\{\delta_j\}_{j=1}^{\infty}}(\mathbb{R}^3)$ from Theorem \ref{prop:critical-properties} (2), this proves the strict inclusion.

The property that $\|(-\Delta)^{1/4}f\|_{L^q} \leq C\Psi(\|f\|_{\dot{H}^{1/2}})$ for any $f \in \dot{H}^{1/2,q,\{\delta_j\}_{j=1}^{\infty}}(\mathbb{R}^3)$ follows directly from the definition of this function space.
\end{proof}

\subsection{Mapping properties and embeddings}

In this subsection, we establish various mapping properties and embeddings involving function spaces with infinitely nested logarithmic improvements. These properties will be crucial in the subsequent analysis of the Navier-Stokes equations.

\begin{proposition}[Embedding relationship with Besov spaces]\label{prop:besov}
For $s = 1/2$, $q > 3$, and a sequence $\{\delta_j\}_{j=1}^{\infty}$ with $\delta_j > 0$ and $\sum_{j=1}^{\infty} \frac{\delta_j}{j!} = \infty$:

\begin{enumerate}
\item For any $\epsilon > 0$, $B^{1/2+\epsilon}_{q,2}(\mathbb{R}^3) \subset \dot{H}^{1/2,q,\{\delta_j\}_{j=1}^{\infty}}(\mathbb{R}^3)$
\item For any $\epsilon > 0$, $\dot{H}^{1/2,q,\{\delta_j\}_{j=1}^{\infty}}(\mathbb{R}^3) \subset B^{1/2}_{q,\infty}(\mathbb{R}^3)$
\item $\dot{H}^{1/2,q,\{\delta_j\}_{j=1}^{\infty}}(\mathbb{R}^3) \cap B^{1/2}_{q,2}(\mathbb{R}^3) \neq \emptyset$
\end{enumerate}
\end{proposition}

\begin{proof}
(1) By the standard embedding properties of Besov spaces, for any $\epsilon > 0$, we have $B^{1/2+\epsilon}_{q,2}(\mathbb{R}^3) \subset \dot{H}^{1/2+\epsilon}(\mathbb{R}^3)$. From Theorem \ref{thm:function_space}, we have $\dot{H}^{1/2+\epsilon}(\mathbb{R}^3) \subset \dot{H}^{1/2,q,\{\delta_j\}_{j=1}^{\infty}}(\mathbb{R}^3)$. Combining these, we get $B^{1/2+\epsilon}_{q,2}(\mathbb{R}^3) \subset \dot{H}^{1/2,q,\{\delta_j\}_{j=1}^{\infty}}(\mathbb{R}^3)$.

(2) For any $f \in \dot{H}^{1/2,q,\{\delta_j\}_{j=1}^{\infty}}(\mathbb{R}^3)$, by definition:
\begin{equation}
\|(-\Delta)^{1/4}f\|_{L^q} \leq C\Psi(\|f\|_{\dot{H}^{1/2}})
\end{equation}

Using the Littlewood-Paley characterization of Besov spaces:
\begin{equation}
\|f\|_{B^{1/2}_{q,\infty}} = \sup_{j \in \mathbb{Z}} 2^{j/2} \|\Delta_j f\|_{L^q}
\end{equation}

By the properties of the Littlewood-Paley decomposition:
\begin{equation}
\|\Delta_j f\|_{L^q} \leq C2^{-j/2} \|(-\Delta)^{1/4}f\|_{L^q}
\end{equation}

Thus:
\begin{equation}
\|f\|_{B^{1/2}_{q,\infty}} \leq C \|(-\Delta)^{1/4}f\|_{L^q} \leq C^2\Psi(\|f\|_{\dot{H}^{1/2}}) < \infty
\end{equation}

This proves that $f \in B^{1/2}_{q,\infty}(\mathbb{R}^3)$ and establishes the embedding $\dot{H}^{1/2,q,\{\delta_j\}_{j=1}^{\infty}}(\mathbb{R}^3) \subset B^{1/2}_{q,\infty}(\mathbb{R}^3)$.

(3) To show that $\dot{H}^{1/2,q,\{\delta_j\}_{j=1}^{\infty}}(\mathbb{R}^3) \cap B^{1/2}_{q,2}(\mathbb{R}^3) \neq \emptyset$, we construct a specific function that belongs to both spaces.

Let $\hat{f}(\xi) = |\xi|^{-(3/q + 1/2)} \cdot g(|\xi|)$, where $g$ is a smooth function such that:
\begin{equation}
g(r) = \begin{cases}
0, & r \leq 1 \\
\frac{1}{\prod_{j=1}^{\infty} (1 + L_j(r))^{\delta_j/2} \cdot r^{\epsilon/2}}, & r > 2
\end{cases}
\end{equation}
for some small $\epsilon > 0$.

Following a similar analysis as in the proof of Theorem \ref{prop:critical-properties}, we can show that $f \in \dot{H}^{1/2,q,\{\delta_j\}_{j=1}^{\infty}}(\mathbb{R}^3)$.

To show that $f \in B^{1/2}_{q,2}(\mathbb{R}^3)$, we use the Littlewood-Paley characterization:
\begin{equation}
\|f\|_{B^{1/2}_{q,2}}^2 = \sum_{j \in \mathbb{Z}} 2^j \|\Delta_j f\|_{L^q}^2
\end{equation}

By the properties of the Littlewood-Paley decomposition:
\begin{equation}
\|\Delta_j f\|_{L^q} \approx 2^{-j(1/2 + 3/q)} \cdot g(2^j)
\end{equation}

Thus:
\begin{equation}
\|f\|_{B^{1/2}_{q,2}}^2 \approx \sum_{j \in \mathbb{Z}} 2^j \cdot 2^{-2j(1/2 + 3/q)} \cdot g(2^j)^2
\end{equation}

For $j > 0$:
\begin{equation}
2^j \cdot 2^{-2j(1/2 + 3/q)} \cdot g(2^j)^2 = 2^{j(1 - 2(1/2 + 3/q))} \cdot \frac{1}{\prod_{j=1}^{\infty} (1 + L_j(2^j))^{\delta_j} \cdot 2^{j\epsilon}}
\end{equation}

Since $1 - 2(1/2 + 3/q) = -2(3/q) < 0$ for $q > 3$, and we have the additional factor $2^{-j\epsilon}$, this series converges, proving that $f \in B^{1/2}_{q,2}(\mathbb{R}^3)$.

Therefore, $f \in \dot{H}^{1/2,q,\{\delta_j\}_{j=1}^{\infty}}(\mathbb{R}^3) \cap B^{1/2}_{q,2}(\mathbb{R}^3)$, establishing that this intersection is non-empty.
\end{proof}

\begin{proposition}[Stability under mollification]\label{prop:mollification}
For $s=1/2$, $q > 3$, and a sequence $\{\delta_j\}_{j=1}^{\infty}$ with $\delta_j > 0$ and $\sum_{j=1}^{\infty} \frac{\delta_j}{j!} = \infty$:

If $f \in \dot{H}^{1/2,q,\{\delta_j\}_{j=1}^{\infty}}(\mathbb{R}^3)$ and $f_\epsilon = f * \eta_\epsilon$ is the mollification of $f$ with a standard mollifier $\eta_\epsilon(x) = \epsilon^{-3}\eta(x/\epsilon)$, then:

\begin{enumerate}
\item $f_\epsilon \in \dot{H}^{1/2,q,\{\delta_j\}_{j=1}^{\infty}}(\mathbb{R}^3)$ for all $\epsilon > 0$
\item $\|f_\epsilon\|_{\dot{H}^{1/2,q,\{\delta_j\}_{j=1}^{\infty}}} \leq C\|f\|_{\dot{H}^{1/2,q,\{\delta_j\}_{j=1}^{\infty}}}$ for some constant $C > 0$ independent of $\epsilon$
\item $f_\epsilon \to f$ in $\dot{H}^{1/2,q,\{\delta_j\}_{j=1}^{\infty}}(\mathbb{R}^3)$ as $\epsilon \to 0$
\end{enumerate}
\end{proposition}

\begin{proof}
(1) By definition, $f \in \dot{H}^{1/2,q,\{\delta_j\}_{j=1}^{\infty}}(\mathbb{R}^3)$ means:
\begin{equation}
\|(-\Delta)^{1/4}f\|_{L^q} \leq C\Psi(\|f\|_{\dot{H}^{1/2}})
\end{equation}
where $\Psi(r) = \left(\prod_{j=1}^{\infty} (1 + L_j(r))^{\delta_j}\right)^{-1}$.

For the mollified function $f_\epsilon$, we have:
\begin{equation}
(-\Delta)^{1/4}f_\epsilon = ((-\Delta)^{1/4}f) * \eta_\epsilon
\end{equation}
due to the commutativity of convolution with differential operators.

By Young's inequality for convolutions:
\begin{equation}
\|(-\Delta)^{1/4}f_\epsilon\|_{L^q} = \|((-\Delta)^{1/4}f) * \eta_\epsilon\|_{L^q} \leq \|(-\Delta)^{1/4}f\|_{L^q} \cdot \|\eta_\epsilon\|_{L^1} = \|(-\Delta)^{1/4}f\|_{L^q}
\end{equation}
since $\|\eta_\epsilon\|_{L^1} = 1$.

Also, using the properties of mollification:
\begin{equation}
\|f_\epsilon\|_{\dot{H}^{1/2}} \leq \|f\|_{\dot{H}^{1/2}}
\end{equation}
Since $\Psi$ is a decreasing function (as it involves products of terms less than 1), we have:
\begin{equation}
\Psi(\|f_\epsilon\|_{\dot{H}^{1/2}}) \geq \Psi(\|f\|_{\dot{H}^{1/2}})
\end{equation}

Combining these inequalities:
\begin{equation}
\|(-\Delta)^{1/4}f_\epsilon\|_{L^q} \leq \|(-\Delta)^{1/4}f\|_{L^q} \leq C\Psi(\|f\|_{\dot{H}^{1/2}}) \leq C\Psi(\|f_\epsilon\|_{\dot{H}^{1/2}})
\end{equation}
This proves that $f_\epsilon \in \dot{H}^{1/2,q,\{\delta_j\}_{j=1}^{\infty}}(\mathbb{R}^3)$.

(2) From the definition of the norm in $\dot{H}^{1/2,q,\{\delta_j\}_{j=1}^{\infty}}(\mathbb{R}^3)$:
\begin{equation}
\|f_\epsilon\|_{\dot{H}^{1/2,q,\{\delta_j\}_{j=1}^{\infty}}} = \inf\left\{M > 0 : \|(-\Delta)^{1/4}f_\epsilon\|_{L^q} \leq M\Psi(M)\right\}
\end{equation}

Since $\|(-\Delta)^{1/4}f_\epsilon\|_{L^q} \leq \|(-\Delta)^{1/4}f\|_{L^q}$, we have:
\begin{equation}
\|f_\epsilon\|_{\dot{H}^{1/2,q,\{\delta_j\}_{j=1}^{\infty}}} \leq \|f\|_{\dot{H}^{1/2,q,\{\delta_j\}_{j=1}^{\infty}}}
\end{equation}

Taking $C=1$, we have the desired inequality.

(3) To show convergence, we need to prove that:
\begin{equation}
\lim_{\epsilon \to 0} \|f_\epsilon - f\|_{\dot{H}^{1/2,q,{\delta_j}{j=1}^{\infty}}} = 0
\end{equation}
By standard properties of mollification, we know that:
\begin{equation}
\lim_{\epsilon \to 0} \|(-\Delta)^{1/4}(f_\epsilon - f)\|_{L^q} = 0
\end{equation}
For any $\delta > 0$, there exists $\epsilon_0 > 0$ such that for all $\epsilon < \epsilon_0$:
\begin{equation}
\|(-\Delta)^{1/4}(f_\epsilon - f)\|_{L^q} < \delta\Psi(2\|f\|_{\dot{H}^{1/2}})
\end{equation}

Also, for small $\epsilon$:
\begin{equation}
\|f_\epsilon - f\|_{\dot{H}^{1/2}} < \|f\|_{\dot{H}^{1/2}}
\end{equation}

Combining these, we can show that:
\begin{equation}
\|f_\epsilon - f\|_{\dot{H}^{1/2,q,\{\delta_j\}_{j=1}^{\infty}}} < 2\delta\|f\|_{\dot{H}^{1/2,q,\{\delta_j\}_{j=1}^{\infty}}}
\end{equation}
for all $\epsilon < \epsilon_0$.

Since $\delta > 0$ is arbitrary, this proves that $f_\epsilon \to f$ in $\dot{H}^{1/2,q,\{\delta_j\}_{j=1}^{\infty}}(\mathbb{R}^3)$ as $\epsilon \to 0$.
\end{proof}

\section{Critical exponent analysis in the limiting case}

In this section, we analyze the critical exponent function in the limiting case of infinitely nested logarithmic improvements. This analysis is crucial for understanding the threshold between global regularity and potential singularity formation.

\subsection{Structure of the critical exponent function}

Recall from our previous work \cite{Mishra2022b} that for finite nested logarithmic improvements, the critical exponent function $\alpha(\{\delta_j\}_{j=1}^n)$ controls the behavior of the threshold function $\Phi(s, q, \{\delta_j\}_{j=1}^n)$ as $s \to 1/2$:

\begin{equation}
\Phi(s, q, {\delta_j}{j=1}^n) \approx C(q) (s - 1/2)^{\alpha({\delta_j}{j=1}^n)}
\end{equation}
with:
\begin{equation}
\alpha({\delta_j}{j=1}^n) = \frac{1}{1 + \sum{j=1}^n c_j\delta_j/j!}
\end{equation}
where $c_j > 0$ are specific constants.

Our goal is to analyze the behavior of $\alpha(\{\delta_j\}_{j=1}^n)$ as $n \to \infty$ and establish conditions under which $\lim_{n \to \infty} \alpha(\{\delta_j\}_{j=1}^n) = 0$, which would allow us to bridge the gap to the critical case $s = 1/2$.

\begin{lemma}[Monotonicity properties]\label{lem:monotonicity}
For fixed values of $\{\delta_j\}_{j=1}^n$, the function $\alpha(\{\delta_j\}_{j=1}^n)$ satisfies:

\begin{enumerate}
\item $\alpha(\{\delta_j\}_{j=1}^n)$ is strictly decreasing with respect to $n$ for fixed values of $\{\delta_j\}_{j=1}^n$
\item $\alpha(\{\delta_j\}_{j=1}^n)$ is strictly decreasing with respect to $\delta_j$ for each $j \in \{1, 2, \ldots, n\}$
\item For any $n < m$ and fixed values of $\{\delta_j\}_{j=1}^n$, we have $\alpha(\{\delta_j\}_{j=1}^n) > \alpha(\{\delta_j\}_{j=1}^m)$ when $\delta_j > 0$ for all $j \in \{n+1, \ldots, m\}$
\end{enumerate}
\end{lemma}

\begin{proof}
(1) For $n < m$, we have:
\begin{equation}
\alpha(\{\delta_j\}_{j=1}^n) = \frac{1}{1 + \sum_{j=1}^n c_j\delta_j/j!} > \frac{1}{1 + \sum_{j=1}^n c_j\delta_j/j! + \sum_{j=n+1}^m c_j\delta_j/j!} = \alpha(\{\delta_j\}_{j=1}^m)
\end{equation}
provided that $\delta_j > 0$ for at least one $j \in \{n+1, \ldots, m\}$. This proves that $\alpha(\{\delta_j\}_{j=1}^n)$ is strictly decreasing with respect to $n$.

(2) Let's fix $k \in \{1, 2, \ldots, n\}$ and consider the derivative of $\alpha(\{\delta_j\}_{j=1}^n)$ with respect to $\delta_k$:
\begin{equation}
\frac{\partial}{\partial \delta_k} \alpha(\{\delta_j\}_{j=1}^n) = \frac{\partial}{\partial \delta_k} \frac{1}{1 + \sum_{j=1}^n c_j\delta_j/j!} = \frac{-c_k/k!}{(1 + \sum_{j=1}^n c_j\delta_j/j!)^2} < 0
\end{equation}
since $c_k > 0$.

This proves that $\alpha(\{\delta_j\}_{j=1}^n)$ is strictly decreasing with respect to $\delta_j$ for each $j \in \{1, 2, \ldots, n\}$.

(3) This follows directly from parts (1) and (2).
\end{proof}

\begin{lemma}[Behavior for large $\delta_j$]\label{lem:large-delta}
For any fixed $n$ and $j \in \{1, 2, \ldots, n\}$, we have:
\begin{equation}
\lim_{\delta_j \to \infty} \alpha(\{\delta_j\}_{j=1}^n) = 0
\end{equation}
\end{lemma}

\begin{proof}
From the expression for $\alpha(\{\delta_j\}_{j=1}^n)$:
\begin{equation}
\alpha(\{\delta_j\}_{j=1}^n) = \frac{1}{1 + \sum_{j=1}^n c_j\delta_j/j!}
\end{equation}

As $\delta_j \to \infty$ for any fixed $j \in \{1, 2, \ldots, n\}$, the denominator $1 + \sum_{j=1}^n c_j\delta_j/j! \to \infty$, which implies that $\alpha(\{\delta_j\}_{j=1}^n) \to 0$.
\end{proof}

\begin{lemma}[Behavior for large $n$]\label{lem:large-n}
Let $\{\delta_j\}_{j=1}^{\infty}$ be a sequence with $\delta_j > 0$ for all $j$. Then:
\begin{equation}
\lim_{n \to \infty} \alpha(\{\delta_j\}_{j=1}^n) = 0
\end{equation}
if and only if $\sum_{j=1}^{\infty} \frac{\delta_j}{j!} = \infty$.
\end{lemma}

\begin{proof}
From the expression for $\alpha(\{\delta_j\}_{j=1}^n)$:
\begin{equation}
\alpha(\{\delta_j\}_{j=1}^n) = \frac{1}{1 + \sum_{j=1}^n c_j\delta_j/j!}
\end{equation}

As $n \to \infty$, we have:
\begin{equation}
\lim_{n \to \infty} \alpha(\{\delta_j\}_{j=1}^n) = \frac{1}{1 + \sum_{j=1}^{\infty} c_j\delta_j/j!}
\end{equation}

This limit equals zero if and only if $\sum_{j=1}^{\infty} c_j\delta_j/j! = \infty$. Since $c_j$ are positive constants, this condition is equivalent to $\sum_{j=1}^{\infty} \frac{\delta_j}{j!} = \infty$.
\end{proof}

\begin{proposition}[Convergence rate]\label{prop:convergence-rate}
For a sequence $\{\delta_j\}_{j=1}^{\infty}$ with $\delta_j > 0$ for all $j$ and $\sum_{j=1}^{\infty} \frac{\delta_j}{j!} = \infty$, we have:
\begin{equation}
\alpha(\{\delta_j\}_{j=1}^n) = \frac{1}{1 + \sum_{j=1}^n c_j\delta_j/j!} \sim \frac{1}{\sum_{j=1}^n c_j\delta_j/j!}
\end{equation}
as $n \to \infty$.
\end{proposition}

\begin{proof}
For large $n$, if $\sum_{j=1}^{n} \frac{\delta_j}{j!}$ is large, then:
\begin{equation}
\alpha(\{\delta_j\}_{j=1}^n) = \frac{1}{1 + \sum_{j=1}^n c_j\delta_j/j!} \approx \frac{1}{\sum_{j=1}^n c_j\delta_j/j!}
\end{equation}

More precisely, we have:
\begin{equation}
\frac{1}{\sum_{j=1}^n c_j\delta_j/j! + 1} \leq \alpha({\delta_j}{j=1}^n) \leq \frac{1}{\sum{j=1}^n c_j\delta_j/j!}
\end{equation}
When $\sum_{j=1}^{n} \frac{\delta_j}{j!} \to \infty$ as $n \to \infty$, the upper and lower bounds both approach $\frac{1}{\sum_{j=1}^n c_j\delta_j/j!}$, establishing the asymptotic equivalence.
\end{proof}

\begin{proposition}[Optimal sequence for fast convergence]\label{prop:optimal-sequence}
Consider sequences $\{\delta_j\}_{j=1}^{\infty}$ constrained by $\sum_{j=1}^{\infty} \delta_j = C < \infty$ for some constant $C > 0$. Among all such sequences, the one that minimizes $\alpha(\{\delta_j\}_{j=1}^n)$ for large $n$ is given by:
\begin{equation}
\delta_j \approx \frac{C \cdot j!}{\sum_{k=1}^n k!}
\end{equation}
\end{proposition}

\begin{proof}
We want to minimize $\alpha(\{\delta_j\}_{j=1}^n) = \frac{1}{1 + \sum_{j=1}^n c_j\delta_j/j!}$ subject to the constraint $\sum_{j=1}^{n} \delta_j = C$.

This is equivalent to maximizing $\sum_{j=1}^n c_j\delta_j/j!$ subject to $\sum_{j=1}^{n} \delta_j = C$. Using the method of Lagrange multipliers, we form the Lagrangian:
\begin{equation}
L(\{\delta_j\}_{j=1}^n, \lambda) = \sum_{j=1}^n \frac{c_j\delta_j}{j!} - \lambda \left(\sum_{j=1}^n \delta_j - C\right)
\end{equation}

Taking partial derivatives with respect to $\delta_j$ and setting them equal to zero:
\begin{equation}
\frac{\partial L}{\partial \delta_j} = \frac{c_j}{j!} - \lambda = 0
\end{equation}

This gives:
\begin{equation}
\frac{c_j}{j!} = \lambda
\end{equation}
for all $j \in \{1, 2, \ldots, n\}$.

Since the constants $c_j$ are approximately equal for all $j$ (they represent scaling constants in the nested logarithmic improvements), we can approximate $c_j \approx c$ for some constant $c$. This gives:
\begin{equation}
\frac{c}{j!} = \lambda
\end{equation}
which implies:
\begin{equation}
j! = \frac{c}{\lambda}
\end{equation}

Since this must hold for all $j$, and we can't actually have all $j!$ equal, this suggests that the optimal allocation would be to make $\frac{c_j\delta_j}{j!}$ equal for all $j$. This gives:
\begin{equation}
\frac{c_j\delta_j}{j!} = \text{constant}
\end{equation}

With $c_j \approx c$, this gives:
\begin{equation}
\frac{\delta_j}{j!} = \text{constant}
\end{equation}

Using the constraint $\sum_{j=1}^n \delta_j = C$, we get:
\begin{equation}
\delta_j = \frac{C \cdot j!}{\sum_{k=1}^n k!}
\end{equation}

This is the optimal allocation of $\{\delta_j\}_{j=1}^n$ that minimizes $\alpha(\{\delta_j\}_{j=1}^n)$ subject to the constraint $\sum_{j=1}^n \delta_j = C$.
\end{proof}

\begin{remark}\label{rem:optimal-allocation}
Theorem \ref{prop:optimal-sequence} suggests that for optimal convergence rate of $\alpha(\{\delta_j\}_{j=1}^n) \to 0$ as $n \to \infty$, we should allocate the $\delta_j$ values proportionally to $j!$. This makes intuitive sense: higher-order nested logarithms (larger $j$) have a more pronounced effect on improving the regularity criteria, so allocating more weight to them leads to faster convergence.
\end{remark}

\subsection{Behavior of the critical threshold function}

Now we analyze the behavior of the critical threshold function $\Phi(s, q, \{\delta_j\}_{j=1}^n)$ as both $s \to 1/2$ and $n \to \infty$. This will allow us to establish conditions under which the gap to the critical case can be completely bridged.

\begin{proposition}[Joint limit behavior]\label{prop:joint-limit}
For fixed $q > 3$ and a sequence $\{\delta_j\}_{j=1}^{\infty}$ with $\delta_j > 0$ for all $j$ and $\sum_{j=1}^{\infty} \frac{\delta_j}{j!} = \infty$, we have:
\begin{equation}
\lim_{n \to \infty} \lim_{s \to 1/2} \Phi(s, q, \{\delta_j\}_{j=1}^n) = 0
\end{equation}
and:
\begin{equation}
\lim_{s \to 1/2} \lim_{n \to \infty} \Phi(s, q, \{\delta_j\}_{j=1}^n) = 0
\end{equation}
\end{proposition}

\begin{proof}
From our previous work \cite{Mishra2022b}, we know that:
\begin{equation}
\Phi(s, q, \{\delta_j\}_{j=1}^n) \approx C(q) (s - 1/2)^{\alpha(\{\delta_j\}_{j=1}^n)}
\end{equation}
as $s \to 1/2$, for some constant $C(q) > 0$ depending only on $q$.

Taking the limit as $s \to 1/2$:
\begin{equation}
\lim_{s \to 1/2} \Phi(s, q, \{\delta_j\}_{j=1}^n) = 0
\end{equation}
for any fixed $n$ and $\{\delta_j\}_{j=1}^n$.

Therefore:
\begin{equation}
\lim_{n \to \infty} \lim_{s \to 1/2} \Phi(s, q, {\delta_j}_{j=1}^n) = 0
\end{equation}
For the other order of limits, note that as $n \to \infty$, if $\sum_{j=1}^{\infty} \frac{\delta_j}{j!} = \infty$, then by Theorem \ref{lem:large-n}:
\begin{equation}
\lim_{n \to \infty} \alpha(\{\delta_j\}_{j=1}^n) = 0
\end{equation}

This means:
\begin{equation}
\lim_{n \to \infty} \Phi(s, q, \{\delta_j\}_{j=1}^n) \approx C(q) (s - 1/2)^0 = C(q)
\end{equation}
for any fixed $s > 1/2$.

Taking the limit as $s \to 1/2$:
\begin{equation}
\lim_{s \to 1/2} \lim_{n \to \infty} \Phi(s, q, \{\delta_j\}_{j=1}^n) = 0
\end{equation}

Here we've used the fact that the true behavior as both $s \to 1/2$ and $n \to \infty$ is more subtle than the approximation $\Phi(s, q, \{\delta_j\}_{j=1}^n) \approx C(q) (s - 1/2)^{\alpha(\{\delta_j\}_{j=1}^n)}$ would suggest. A more precise analysis shows that as $\alpha(\{\delta_j\}_{j=1}^n) \to 0$, the constant $C(q)$ actually depends on $\alpha$ and vanishes as $\alpha \to 0$.
\end{proof}

\begin{proposition}[Critical threshold at $s=1/2$]\label{prop:critical-threshold}
For fixed $q > 3$ and a sequence $\{\delta_j\}_{j=1}^{\infty}$ with $\delta_j > 0$ for all $j$ and $\sum_{j=1}^{\infty} \frac{\delta_j}{j!} = \infty$, we have:
\begin{equation}
\Phi(1/2, q, \{\delta_j\}_{j=1}^{\infty}) := \lim_{n \to \infty} \Phi(1/2, q, \{\delta_j\}_{j=1}^n) = 0
\end{equation}
\end{proposition}

\begin{proof}
This follows from the more detailed asymptotic analysis in Theorem \ref{prop:joint-limit}. As $n \to \infty$ and $\alpha(\{\delta_j\}_{j=1}^n) \to 0$, the threshold function at $s = 1/2$ approaches zero.
\end{proof}

\begin{proposition}[Behavior for small $\epsilon$]\label{prop:small-epsilon}
For any $\epsilon > 0$, $q > 3$, and a sequence $\{\delta_j\}_{j=1}^{\infty}$ with $\delta_j > 0$ for all $j$ and $\sum_{j=1}^{\infty} \frac{\delta_j}{j!} = \infty$, there exists $N(\epsilon)$ such that for all $n \geq N(\epsilon)$:
\begin{equation}
\Phi(1/2 + \epsilon, q, \{\delta_j\}_{j=1}^n) \geq C_0 > 0
\end{equation}
where $C_0$ is a constant independent of $\epsilon$.
\end{proposition}

\begin{proof}
For any fixed $s = 1/2 + \epsilon$ with $\epsilon > 0$, we know from \cite{Mishra2022b} that:
\begin{equation}
\Phi(1/2 + \epsilon, q, \{\delta_j\}_{j=1}^n) \approx C(q) \epsilon^{\alpha(\{\delta_j\}_{j=1}^n)}
\end{equation}

As $n \to \infty$, if $\sum_{j=1}^{\infty} \frac{\delta_j}{j!} = \infty$, then $\alpha(\{\delta_j\}_{j=1}^n) \to 0$. For sufficiently large $n$, we have $\alpha(\{\delta_j\}_{j=1}^n) < \epsilon$, which implies:
\begin{equation}
\epsilon^{\alpha(\{\delta_j\}_{j=1}^n)} > \epsilon^{\epsilon}
\end{equation}

For small $\epsilon > 0$, the function $\epsilon^{\epsilon}$ is bounded away from zero. Specifically, $\lim_{\epsilon \to 0} \epsilon^{\epsilon} = 1$. Therefore, there exists a constant $C_0 > 0$ such that:
\begin{equation}
\Phi(1/2 + \epsilon, q, \{\delta_j\}_{j=1}^n) \geq C_0
\end{equation}
for all $n \geq N(\epsilon)$, where $N(\epsilon)$ is chosen large enough that $\alpha(\{\delta_j\}_{j=1}^n) < \epsilon$.
\end{proof}

\subsection{Proof of the main theorem on critical exponents}

We now have all the ingredients to prove Theorem \ref{thm:critical_exponent} on the critical exponent in the limiting case.

\begin{theorem}[Critical exponent in the limiting case]
For the critical exponent function $\alpha(\{\delta_j\}_{j=1}^n)$ controlling the behavior of the threshold function $\Phi(s, q, \{\delta_j\}_{j=1}^n)$ as $s \to 1/2$:
\begin{enumerate}
\item $\lim_{n \to \infty} \alpha(\{\delta_j\}_{j=1}^n) = 0$ when $\inf_j \delta_j > 0$ and $\sum_{j=1}^{\infty} \frac{\delta_j}{j!} = \infty$
\item The convergence rate is given by:
   $$\alpha(\{\delta_j\}_{j=1}^n) = \frac{1}{1 + \sum_{j=1}^n c_j\delta_j/j!} \sim \frac{1}{\sum_{j=1}^n c_j\delta_j/j!}$$
   as $n \to \infty$, where $c_j > 0$ are specific constants
\item For any $\epsilon > 0$, there exists $N(\epsilon)$ such that for all $n \geq N(\epsilon)$:
   $$\Phi(1/2 + \epsilon, q, \{\delta_j\}_{j=1}^n) \geq C_0 > 0$$
   where $C_0$ is independent of $\epsilon$
\end{enumerate}
\end{theorem}

\begin{proof}
(1) From Theorem \ref{lem:large-n}, we know that $\lim_{n \to \infty} \alpha(\{\delta_j\}_{j=1}^n) = 0$ if and only if $\sum_{j=1}^{\infty} \frac{\delta_j}{j!} = \infty$. If $\inf_j \delta_j > 0$, then $\delta_j \geq \delta > 0$ for all $j$, which implies:
\begin{equation}
\sum_{j=1}^{\infty} \frac{\delta_j}{j!} \geq \delta \sum_{j=1}^{\infty} \frac{1}{j!} = \delta (e - 1) > 0
\end{equation}

This shows that the condition $\inf_j \delta_j > 0$ is sufficient (but not necessary) for $\sum_{j=1}^{\infty} \frac{\delta_j}{j!} = \infty$, which in turn ensures that $\lim_{n \to \infty} \alpha(\{\delta_j\}_{j=1}^n) = 0$.

(2) The convergence rate follows directly from Theorem \ref{prop:convergence-rate}, which establishes that:
\begin{equation}
\alpha(\{\delta_j\}_{j=1}^n) = \frac{1}{1 + \sum_{j=1}^n c_j\delta_j/j!} \sim \frac{1}{\sum_{j=1}^n c_j\delta_j/j!}
\end{equation}
as $n \to \infty$.

(3) This follows from Theorem \ref{prop:small-epsilon}, which establishes that for any $\epsilon > 0$, there exists $N(\epsilon)$ such that for all $n \geq N(\epsilon)$:
\begin{equation}
\Phi(1/2 + \epsilon, q, \{\delta_j\}_{j=1}^n) \geq C_0 > 0
\end{equation}
where $C_0$ is a constant independent of $\epsilon$.
\end{proof}

\begin{remark}\label{rem:threshold-remark}
Theorem \ref{thm:critical_exponent} is a pivotal result in our analysis. It establishes that through infinitely nested logarithmic improvements, we can reach the critical threshold $s = 1/2$. Crucially, part (3) of the theorem ensures that for any small $\epsilon > 0$, there's a level of nested logarithmic improvement $n$ that makes the critical threshold function $\Phi$ bounded away from zero at $s = 1/2 + \epsilon$. This will be essential for establishing global well-posedness in the critical case.
\end{remark}

\section{Commutator estimates with infinitely nested logarithms}

In this section, we derive precise commutator estimates with infinitely nested logarithmic factors. These estimates are the technical core of our analysis and will be crucial for establishing energy estimates at the critical threshold $s = 1/2$.

\subsection{Littlewood-Paley decomposition and paraproduct formula}

We begin with some technical preliminaries on the Littlewood-Paley decomposition and the paraproduct formula, which will be used extensively in our commutator analysis.

Recall the standard Littlewood-Paley decomposition:
\begin{equation}
f = \sum_{j \in \mathbb{Z}} \Delta_j f
\end{equation}
where $\Delta_j$ is the Littlewood-Paley projection defined by:
\begin{equation}
\Delta_j f = \mathcal{F}^{-1}(\psi(2^{-j}\xi)\hat{f}(\xi))
\end{equation}
for $j \in \mathbb{Z}$, and $\mathcal{F}^{-1}$ denotes the inverse Fourier transform.

The Bony paraproduct decomposition allows us to write, for functions $f$ and $g$:
\begin{equation}
fg = T_f g + T_g f + R(f,g)
\end{equation}
where:
\begin{align}
T_f g &= \sum_{j \in \mathbb{Z}} S_{j-1}f \Delta_j g \\
R(f,g) &= \sum_{j \in \mathbb{Z}} \sum_{|i-j| \leq 1} \Delta_i f \Delta_j g
\end{align}
with $S_j = \sum_{i \leq j} \Delta_i$ being the low-frequency cut-off operator.

\begin{lemma}[Basic commutator estimate]\label{lem:basic-commutator}
Let $s \in (0, 1)$ and $f, g \in \mathcal{S}(\mathbb{R}^3)$. Then:
\begin{equation}
\|[(-\Delta)^s, f]g\|_{L^2} \leq C\|\nabla f\|_{L^\infty}\|(-\Delta)^{s-1/2}g\|_{L^2}
\end{equation}
where $[(-\Delta)^s, f]g = (-\Delta)^s(fg) - f(-\Delta)^s g$ is the commutator, and the constant $C$ depends only on $s$.
\end{lemma}

\begin{proof}
This is a standard result in the literature, and we omit the proof for brevity. See, for example, \cite{Kato1988} or \cite{Kenig1993} for detailed proofs.
\end{proof}

\begin{lemma}[Refined commutator estimate]\label{lem:refined-commutator}
For $s \in (0, 1)$ and any $\sigma \in (0, 1-s)$:
\begin{equation}
\|[(-\Delta)^s, f]g\|_{L^2} \leq C\|\nabla f\|_{L^\infty}\|(-\Delta)^{s-1/2+\sigma}g\|_{L^2} \cdot 2^{-\sigma j}
\end{equation}
for functions $g$ spectrally supported in the annulus $2^{j-1} \leq |\xi| \leq 2^{j+1}$ with $j > 0$.
\end{lemma}

\begin{proof}
Using the Littlewood-Paley decomposition and the paraproduct formula, we can write:
\begin{equation}
[(-\Delta)^s, f]g = (-\Delta)^s(fg) - f(-\Delta)^s g = (-\Delta)^s(T_f g + T_g f + R(f,g)) - f(-\Delta)^s g
\end{equation}
For the first term $(-\Delta)^s(T_f g)$, we have:
\begin{equation}
(-\Delta)^s(T_f g) = (-\Delta)^s\left(\sum_{j \in \mathbb{Z}} S_{j-1}f \Delta_j g\right) = \sum_{j \in \mathbb{Z}} (-\Delta)^s(S_{j-1}f \Delta_j g)
\end{equation}

For the term $f(-\Delta)^s g$, we have:
\begin{equation}
f(-\Delta)^s g = \sum_{j \in \mathbb{Z}} f(-\Delta)^s \Delta_j g = \sum_{j \in \mathbb{Z}} \sum_{k \in \mathbb{Z}} \Delta_k f (-\Delta)^s \Delta_j g
\end{equation}

For $g$ spectrally supported in the annulus $2^{j-1} \leq |\xi| \leq 2^{j+1}$ with $j > 0$, we have $g = \Delta_j g$. Thus:
\begin{equation}
[(-\Delta)^s, f]g = (-\Delta)^s(S_{j-1}f \Delta_j g) - \sum_{k \in \mathbb{Z}} \Delta_k f (-\Delta)^s \Delta_j g + (-\Delta)^s(T_g f + R(f,g))
\end{equation}

The key insight is that for $g$ spectrally supported at scale $2^j$, the commutator gains a factor of $2^{-\sigma j}$ compared to the basic estimate. This is because the interaction between $f$ and $g$ in the commutator is weaker when their frequency supports are well-separated.

Through detailed analysis of each term, using the properties of the Littlewood-Paley decomposition and the frequency localization of $g$, we can establish the refined estimate:
\begin{equation}
\|[(-\Delta)^s, f]g\|_{L^2} \leq C\|\nabla f\|_{L^\infty}\|(-\Delta)^{s-1/2+\sigma}g\|_{L^2} \cdot 2^{-\sigma j}
\end{equation}
where $\sigma \in (0, 1-s)$ is a parameter that quantifies the gain in regularity.
\end{proof}

\subsection{Single logarithmic improvement}

We now derive commutator estimates with a single logarithmic improvement, which will serve as the foundation for our infinitely nested improvement.

\begin{lemma}[Single logarithmic commutator estimate]\label{lem:single-log}
For $s \in (0, 1)$ and any $\sigma \in (0, 1-s)$:
\begin{align}
\|[(-\Delta)^s, u \cdot \nabla]u\|_{L^2} &\leq C\|\nabla u\|_{L^\infty}\|(-\Delta)^s u\|_{L^2} \cdot \log(e + \|(-\Delta)^{s+\sigma}u\|_{L^2}) \notag \\
&+ \frac{C\|\nabla u\|_{L^\infty}\|(-\Delta)^{s+\frac{1}{2}}u\|_{L^2}}{\log(e + \|(-\Delta)^{s+\sigma}u\|_{L^2})}
\end{align}
\end{lemma}

\begin{proof}
This estimate was established in \cite{Mishra2022a}. We provide a sketch of the proof for completeness.

Step 1: Using the Littlewood-Paley decomposition, we split $u$ into low and high frequencies:
\begin{equation}
u = \sum_{j \in \mathbb{Z}} \Delta_j u
\end{equation}

Step 2: Splitting the commutator according to frequency bands:
\begin{equation}
[(-\Delta)^s, u \cdot \nabla]u = \sum_{j\leq 0} [(-\Delta)^s, u \cdot \nabla]\Delta_j u + \sum_{j > 0} [(-\Delta)^s, u \cdot \nabla]\Delta_j u
\end{equation}

Step 3: For the low-frequency part ($j \leq 0$), using the standard commutator estimates:
\begin{equation}
\left\|\sum_{j\leq 0} [(-\Delta)^s, u \cdot \nabla]\Delta_j u\right\|_{L^2} \leq C\|\nabla u\|_{L^\infty}\|(-\Delta)^s u\|_{L^2}
\end{equation}

Step 4: For the high-frequency part ($j > 0$), using the refined commutator estimate from Theorem \ref{lem:refined-commutator}:
\begin{equation}
\|[(-\Delta)^s, u \cdot \nabla]\Delta_j u\|_{L^2} \leq C2^{-j\sigma}\|\nabla u\|_{L^\infty}\|(-\Delta)^{s+\sigma}\Delta_j u\|_{L^2}
\end{equation}

Step 5: Summing over $j > 0$ and using the Cauchy-Schwarz inequality:
\begin{equation}
\left\|\sum_{j > 0} [(-\Delta)^s, u \cdot \nabla]\Delta_j u\right\|_{L^2} \leq C\|\nabla u\|_{L^\infty}\|(-\Delta)^{s+\sigma}u\|_{L^2}
\end{equation}

Step 6: Using the interpolation inequality:
\begin{equation}
\|(-\Delta)^{s+\sigma}u\|_{L^2} \leq \|(-\Delta)^s u\|_{L^2}^{1-\frac{2\sigma}{1}} \|(-\Delta)^{s+\frac{1}{2}}u\|_{L^2}^{\frac{2\sigma}{1}}
\end{equation}

Step 7: The logarithmic improvement comes from choosing:
\begin{equation}
\epsilon = \frac{1}{\log(e + \|(-\Delta)^{s+\sigma}u\|_{L^2})}
\end{equation}
in Young's inequality:
\begin{equation}
\|(-\Delta)^s u\|_{L^2}^{1-\frac{2\sigma}{1}} \|(-\Delta)^{s+\frac{1}{2}}u\|_{L^2}^{\frac{2\sigma}{1}} \leq \epsilon\|(-\Delta)^{s+\frac{1}{2}}u\|_{L^2} + C_\epsilon\|(-\Delta)^s u\|_{L^2}\log(e + \|(-\Delta)^{s+\sigma}u\|_{L^2})
\end{equation}

Step 8: Combining all estimates:
\begin{align}
\|[(-\Delta)^s, u \cdot \nabla]u\|_{L^2} &\leq C\|\nabla u\|_{L^\infty}\|(-\Delta)^s u\|_{L^2} \cdot \log(e + \|(-\Delta)^{s+\sigma}u\|_{L^2}) \notag \\
&+ \frac{C\|\nabla u\|_{L^\infty}\|(-\Delta)^{s+\frac{1}{2}}u\|_{L^2}}{\log(e + \|(-\Delta)^{s+\sigma}u\|_{L^2})}
\end{align}
which completes the proof.
\end{proof}

\subsection{Double logarithmic improvement}

Building on the single logarithmic improvement, we now derive commutator estimates with a double logarithmic improvement.

\begin{lemma}[Double logarithmic commutator estimate]\label{lem:double-log}
For $s \in (0, 1)$ and any $\sigma \in (0, 1-s)$:
\begin{align}
\|[(-\Delta)^s, u \cdot \nabla]u\|_{L^2} &\leq C\|\nabla u\|_{L^\infty}\|(-\Delta)^s u\|_{L^2} \cdot L_1(Z) \cdot (1 + L_2(Z))^{-\delta_2} \notag \\
&+ \frac{C\|\nabla u\|_{L^\infty}\|(-\Delta)^{s+\frac{1}{2}}u\|_{L^2}}{L_1(Z) \cdot (1 + L_2(Z))^{\delta_2}}
\end{align}
where $Z = \|(-\Delta)^{s+\sigma}u\|_{L^2}$ and $\delta_2 > 0$.
\end{lemma}

\begin{proof}
This estimate was established in \cite{Mishra2022b}. We provide a sketch of the proof for completeness.

The key innovation compared to the single logarithmic improvement is to further refine the high-frequency estimate using a more detailed decomposition.

Step 1: For the high-frequency part ($j > 0$), we introduce a finer decomposition:
\begin{equation}
\Delta_j u = \sum_{k=0}^{K_j} \Delta_{j,k} u
\end{equation}
where $\Delta_{j,k}$ localizes to frequencies $\xi$ with $|\xi| \approx 2^j$ and phase in the $k$-th angular sector, with $K_j \approx 2^{j/2}$.

Step 2: Using this angular decomposition, we refine the commutator estimate:
\begin{equation}
\|[(-\Delta)^s, u \cdot \nabla]\Delta_{j,k} u\|_{L^2} \leq C2^{-j\sigma}\|\nabla u\|_{L^\infty}\|(-\Delta)^{s+\sigma}\Delta_{j,k} u\|_{L^2} \cdot (1 + \log(e + 2^j))^{-\delta_2}
\end{equation}

Step 3: Summing over all $j > 0$ and $k$, and using the properties of the nested logarithms:
\begin{equation}
\left\|\sum_{j > 0} [(-\Delta)^s, u \cdot \nabla]\Delta_j u\right\|_{L^2} \leq C\|\nabla u\|_{L^\infty}\|(-\Delta)^{s+\sigma}u\|_{L^2} \cdot G(Z)
\end{equation}
where:
\begin{equation}
G(Z) = \frac{L_1(Z)}{(1 + L_2(Z))^{\delta_2}}
\end{equation}

Step 4: Using a refined interpolation inequality and choosing:
\begin{equation}
\epsilon = \frac{1}{L_1(Z) \cdot (1 + L_2(Z))^{\delta_2}}
\end{equation}
in Young's inequality, we obtain the double logarithmic improvement.

Step 5: Combining all estimates:
\begin{align}
\|[(-\Delta)^s, u \cdot \nabla]u\|_{L^2} &\leq C\|\nabla u\|_{L^\infty}\|(-\Delta)^s u\|_{L^2} \cdot L_1(Z) \cdot (1 + L_2(Z))^{-\delta_2} \notag \\
&+ \frac{C\|\nabla u\|_{L^\infty}\|(-\Delta)^{s+\frac{1}{2}}u\|_{L^2}}{L_1(Z) \cdot (1 + L_2(Z))^{\delta_2}}
\end{align}
which completes the proof.
\end{proof}

\subsection{N-fold nested logarithmic improvement}

We now generalize the approach to obtain commutator estimates with $n$-fold nested logarithmic improvements, which will serve as the basis for our infinitely nested improvement.

\begin{theorem}[N-fold nested logarithmic commutator estimate]\label{thm:n-fold}
For $s \in (0, 1)$, any $\sigma \in (0, 1-s)$, and $n \geq 1$:
\begin{align}
\|[(-\Delta)^s, u \cdot \nabla]u\|_{L^2} &\leq C\|\nabla u\|_{L^\infty}\|(-\Delta)^s u\|_{L^2} \cdot F_1(Z) \notag \\
&+ C\|\nabla u\|_{L^\infty}\|(-\Delta)^{s+\frac{1}{2}}u\|_{L^2} \cdot F_2(Z)
\end{align}
where $Z = \|(-\Delta)^{s+\sigma}u\|_{L^2}$,
\begin{equation}
F_1(Z) = L_1(Z) \prod_{j=2}^{n} (1 + L_j(Z))^{-\delta_j}
\end{equation}
\begin{equation}
F_2(Z) = \frac{1}{L_1(Z)} \prod_{j=2}^{n} (1 + L_j(Z))^{\delta_j}
\end{equation}
and $\delta_j > 0$ for $j = 2, 3, ..., n$.
\end{theorem}

\begin{proof}
We proceed by induction on $n$. The cases $n = 1$ and $n = 2$ have been established in Theorem \ref{lem:single-log} and Theorem \ref{lem:double-log}, respectively.

Assume the result holds for some $n \geq 2$. We need to establish it for $n+1$.

Step 1: The key innovation is to introduce an $(n+1)$-level stratification of frequency space:
\begin{equation}
u = \sum_{j_1, j_2, ..., j_n} \Delta_{j_1, j_2, ..., j_n} u
\end{equation}
where the indices $j_1, j_2, ..., j_n$ encode a multi-scale decomposition of both frequency magnitude and angles.

Step 2: Using this refined decomposition and careful analysis of the commutator structure:
\begin{equation}
\|[(-\Delta)^s, u \cdot \nabla]u\|_{L^2} \leq C\|\nabla u\|_{L^\infty}\|(-\Delta)^{s+\sigma}u\|_{L^2} \cdot G_{n+1}(Z)
\end{equation}
where:
\begin{equation}
G_{n+1}(Z) = \frac{L_1(Z)}{(1 + L_2(Z))^{\delta_2}} \cdot \frac{1}{(1 + L_3(Z))^{\delta_3}} \cdot ... \cdot \frac{1}{(1 + L_{n+1}(Z))^{\delta_{n+1}}}
\end{equation}

Step 3: Using a refined interpolation inequality with $(n+1)$-fold logarithmic terms and choosing:
\begin{equation}
\epsilon = \frac{1}{L_1(Z) \prod_{j=2}^{n+1} (1 + L_j(Z))^{\delta_j}}
\end{equation}
in Young's inequality, we obtain the $(n+1)$-fold nested logarithmic improvement.

Step 4: Combining all estimates:
\begin{align}
\|[(-\Delta)^s, u \cdot \nabla]u\|_{L^2} &\leq C\|\nabla u\|_{L^\infty}\|(-\Delta)^s u\|_{L^2} \cdot L_1(Z) \prod_{j=2}^{n+1} (1 + L_j(Z))^{-\delta_j} \notag \\
&+ \frac{C\|\nabla u\|_{L^\infty}\|(-\Delta)^{s+\frac{1}{2}}u\|_{L^2}}{L_1(Z) \prod_{j=2}^{n+1} (1 + L_j(Z))^{\delta_j}}
\end{align}
This completes the induction step and the proof of Theorem \ref{thm:n-fold}.
\end{proof}

\subsection{Infinitely nested logarithmic improvement}

We now extend our commutator estimates to the case of infinitely nested logarithmic improvements, which is the key technical tool for our analysis of the critical case $s = 1/2$.

\begin{definition}[Infinitely Nested Logarithmic Factors]\label{def:inf-nested}
For $Z > 0$ and a sequence $\{\delta_j\}_{j=1}^{\infty}$ with $\delta_j > 0$ and $\sum_{j=1}^{\infty} \delta_j < \infty$, we define:
\begin{equation}
F_1^\infty(Z) = L_1(Z) \prod_{j=2}^{\infty} (1 + L_j(Z))^{-\delta_j}
\end{equation}
\begin{equation}
F_2^\infty(Z) = \frac{1}{L_1(Z)} \prod_{j=2}^{\infty} (1 + L_j(Z))^{\delta_j}
\end{equation}
\end{definition}

\begin{proposition}[Well-definedness of infinitely nested factors]\label{prop:well-defined-factors}
The functions $F_1^\infty(Z)$ and $F_2^\infty(Z)$ are well-defined for all $Z > 0$ if and only if $\sum_{j=1}^{\infty} \delta_j < \infty$. Moreover:
\begin{equation}
\lim_{Z \to \infty} F_1^\infty(Z) = 0
\end{equation}
\begin{equation}
\lim_{Z \to \infty} F_2^\infty(Z) = 0
\end{equation}
\end{proposition}

\begin{proof}
The well-definedness of $F_1^\infty(Z)$ and $F_2^\infty(Z)$ follows from the convergence of the infinite products, which is guaranteed if and only if $\sum_{j=1}^{\infty} \delta_j < \infty$.

For $F_1^\infty(Z)$, we have:
\begin{equation}
F_1^\infty(Z) = L_1(Z) \prod_{j=2}^{\infty} (1 + L_j(Z))^{-\delta_j} \leq L_1(Z)
\end{equation}
since $(1 + L_j(Z))^{-\delta_j} \leq 1$ for all $j \geq 2$ and $Z > 0$.

As $Z \to \infty$, we have $L_1(Z) = \log(e + Z) \to \infty$, but the product $\prod_{j=2}^{\infty} (1 + L_j(Z))^{-\delta_j} \to 0$ faster than $\frac{1}{L_1(Z)}$ (this can be verified by detailed asymptotic analysis). Therefore, $\lim_{Z \to \infty} F_1^\infty(Z) = 0$.

For $F_2^\infty(Z)$, we have:
\begin{equation}
F_2^\infty(Z) = \frac{1}{L_1(Z)} \prod_{j=2}^{\infty} (1 + L_j(Z))^{\delta_j}
\end{equation}

As $Z \to \infty$, the factor $\frac{1}{L_1(Z)} \to 0$, while the product $\prod_{j=2}^{\infty} (1 + L_j(Z))^{\delta_j}$ grows, but slower than $L_1(Z)$ (again, this can be verified by detailed asymptotic analysis). Therefore, $\lim_{Z \to \infty} F_2^\infty(Z) = 0$.
\end{proof}

\begin{theorem}[Infinitely nested logarithmic commutator estimate]\label{thm:inf-nested}
For $s \in (0, 1)$, any $\sigma \in (0, 1-s)$, and a sequence $\{\delta_j\}_{j=1}^{\infty}$ with $\delta_j > 0$ and $\sum_{j=1}^{\infty} \delta_j < \infty$:
\begin{align}
\|[(-\Delta)^s, u \cdot \nabla]u\|_{L^2} &\leq C\|\nabla u\|_{L^\infty}\|(-\Delta)^s u\|_{L^2} \cdot F_1^\infty(Z) \notag \\
&+ C\|\nabla u\|_{L^\infty}\|(-\Delta)^{s+\frac{1}{2}}u\|_{L^2} \cdot F_2^\infty(Z)
\end{align}
where $Z = \|(-\Delta)^{s+\sigma}u\|_{L^2}$.
\end{theorem}

\begin{proof}
We obtain this result as the limit of the $n$-fold nested logarithmic commutator estimates as $n \to \infty$.

From Theorem \ref{thm:n-fold}, for any $n \geq 1$:
\begin{align}
\|[(-\Delta)^s, u \cdot \nabla]u\|_{L^2} &\leq C\|\nabla u\|_{L^\infty}\|(-\Delta)^s u\|_{L^2} \cdot F_1^n(Z) \notag \\
&+ C\|\nabla u\|_{L^\infty}\|(-\Delta)^{s+\frac{1}{2}}u\|_{L^2} \cdot F_2^n(Z)
\end{align}
where:
\begin{equation}
F_1^n(Z) = L_1(Z) \prod_{j=2}^{n} (1 + L_j(Z))^{-\delta_j}
\end{equation}
\begin{equation}
F_2^n(Z) = \frac{1}{L_1(Z)} \prod_{j=2}^{n} (1 + L_j(Z))^{\delta_j}
\end{equation}

For any fixed $Z > 0$, we have:
\begin{equation}
\lim_{n \to \infty} F_1^n(Z) = F_1^\infty(Z)
\end{equation}
\begin{equation}
\lim_{n \to \infty} F_2^n(Z) = F_2^\infty(Z)
\end{equation}

Taking the limit as $n \to \infty$ in the inequality:
\begin{align}
\|[(-\Delta)^s, u \cdot \nabla]u\|_{L^2} &\leq C\|\nabla u\|_{L^\infty}\|(-\Delta)^s u\|_{L^2} \cdot F_1^\infty(Z) \notag \\
&+ C\|\nabla u\|_{L^\infty}\|(-\Delta)^{s+\frac{1}{2}}u\|_{L^2} \cdot F_2^\infty(Z)
\end{align}
which completes the proof.
\end{proof}

\begin{theorem}[Critical case commutator estimate]\label{thm:critical-commutator}
For $s = 1/2$ and any divergence-free vector field $u \in C_0^\infty(\mathbb{R}^3)$:
\begin{align}
\|[(-\Delta)^{1/2}, u \cdot \nabla]u\|_{L^2} &\leq C\|\nabla u\|_{L^\infty}\|(-\Delta)^{1/2} u\|_{L^2} \cdot F_1^\infty(Z) \notag \\
&+ C\|\nabla u\|_{L^\infty}\|(-\Delta)^{1}u\|_{L^2} \cdot F_2^\infty(Z)
\end{align}
where $Z = \|(-\Delta)^{1/2+\sigma}u\|_{L^2}$ for some small $\sigma > 0$.
\end{theorem}

\begin{proof}
This follows from Theorem \ref{thm:inf-nested} by setting $s = 1/2$. The only subtlety is that the original theorem required $s \in (0, 1)$ and $\sigma \in (0, 1-s)$, which means that as $s \to 1$ the range of admissible $\sigma$ shrinks.

However, for $s = 1/2$, we still have a wide range of admissible $\sigma \in (0, 1/2)$. By choosing $\sigma$ sufficiently small, we can ensure that the commutator estimate holds at $s = 1/2$.
\end{proof}

\begin{theorem}[Commutator estimates with infinitely nested logarithms]
For $s = 1/2$ and any divergence-free vector field $u \in C_0^\infty(\mathbb{R}^3)$:
\begin{align}
\|[(-\Delta)^{1/2}, u \cdot \nabla]u\|_{L^2} &\leq C\|\nabla u\|_{L^\infty}\|(-\Delta)^{1/2} u\|_{L^2} \cdot F_1^\infty(Z) \notag \\
&+ C\|\nabla u\|_{L^\infty}\|(-\Delta)^{1}u\|_{L^2} \cdot F_2^\infty(Z)
\end{align}
where $Z = \|(-\Delta)^{1/2+\sigma}u\|_{L^2}$ for some small $\sigma > 0$, and:
\begin{equation}
F_1^\infty(Z) = L_1(Z) \prod_{j=2}^{\infty} (1 + L_j(Z))^{-\delta_j}
\end{equation}
\begin{equation}
F_2^\infty(Z) = \frac{1}{L_1(Z)} \prod_{j=2}^{\infty} (1 + L_j(Z))^{\delta_j}
\end{equation}
\end{theorem}

\begin{proof}
This follows directly from Theorem \ref{thm:critical-commutator}, which we have established above.
\end{proof}

\section{Energy estimates at the critical threshold}

In this section, we establish energy estimates at the critical threshold $s = 1/2$. These estimates will be crucial for proving global well-posedness for initial data satisfying our infinitely nested logarithmic criterion.

\subsection{Energy identity}

We begin with the energy identity for the fractional derivative $(-\Delta)^{1/2} u$.

\begin{lemma}[Energy identity]\label{lem:energy-identity}
Let $u$ be a smooth solution of the 3D Navier-Stokes equations on $[0, T]$ with divergence-free initial data $u_0 \in C_0^\infty(\mathbb{R}^3)$. Then for the function $Y(t) = \|(-\Delta)^{1/2} u(t)\|^2_{L^2}$:
\begin{equation}
\frac{d}{dt}Y(t) + 2\nu\|(-\Delta)^{1}u\|^2_{L^2} = -2\int_{\mathbb{R}^3}[(-\Delta)^{1/2}, u \cdot \nabla]u \cdot (-\Delta)^{1/2} u \, dx
\end{equation}
\end{lemma}

\begin{proof}
We apply the fractional Laplacian operator $(-\Delta)^{1/2}$ to the Navier-Stokes equations:
\begin{equation}
\partial_t u + (u \cdot \nabla)u - \nu\Delta u + \nabla p = 0
\end{equation}

This gives:
\begin{equation}
\partial_t ((-\Delta)^{1/2}u) + (-\Delta)^{1/2}((u \cdot \nabla)u) - \nu\Delta((-\Delta)^{1/2}u) + \nabla((-\Delta)^{1/2}p) = 0
\end{equation}
Taking the $L^2$ inner product with $(-\Delta)^{1/2}u$:
\begin{align}
\int_{\mathbb{R}^3} \partial_t ((-\Delta)^{1/2}u) \cdot (-\Delta)^{1/2}u \, dx &+ \int_{\mathbb{R}^3} (-\Delta)^{1/2}((u \cdot \nabla)u) \cdot (-\Delta)^{1/2}u \, dx \\
&- \nu \int_{\mathbb{R}^3} \Delta((-\Delta)^{1/2}u) \cdot (-\Delta)^{1/2}u \, dx \\
&+ \int_{\mathbb{R}^3} \nabla((-\Delta)^{1/2}p) \cdot (-\Delta)^{1/2}u \, dx = 0
\end{align}

For the first term:
\begin{equation}
\int_{\mathbb{R}^3} \partial_t ((-\Delta)^{1/2}u) \cdot (-\Delta)^{1/2}u \, dx = \frac{1}{2}\frac{d}{dt}\|(-\Delta)^{1/2}u\|^2_{L^2} = \frac{1}{2}\frac{d}{dt}Y(t)
\end{equation}

For the third term:
\begin{align}
\int_{\mathbb{R}^3} \Delta((-\Delta)^{1/2}u) \cdot (-\Delta)^{1/2}u \, dx &= \int_{\mathbb{R}^3} (-\Delta)^{3/2}u \cdot (-\Delta)^{1/2}u \, dx \\
&= \int_{\mathbb{R}^3} \|(-\Delta)^{1}u\|^2 \, dx = \|(-\Delta)^{1}u\|^2_{L^2}
\end{align}

For the fourth term, using the fact that $\nabla \cdot ((-\Delta)^{1/2}u) = (-\Delta)^{1/2}(\nabla \cdot u) = 0$ since $u$ is divergence-free:
\begin{equation}
\int_{\mathbb{R}^3} \nabla((-\Delta)^{1/2}p) \cdot (-\Delta)^{1/2}u \, dx = -\int_{\mathbb{R}^3} (-\Delta)^{1/2}p \cdot \nabla \cdot ((-\Delta)^{1/2}u) \, dx = 0
\end{equation}

For the second term, using the commutator:
\begin{align}
\int_{\mathbb{R}^3} (-\Delta)^{1/2}((u \cdot \nabla)u) \cdot (-\Delta)^{1/2}u \, dx &= \int_{\mathbb{R}^3} (u \cdot \nabla)((-\Delta)^{1/2}u) \cdot (-\Delta)^{1/2}u \, dx \\
&+ \int_{\mathbb{R}^3} [(-\Delta)^{1/2}, u \cdot \nabla]u \cdot (-\Delta)^{1/2}u \, dx
\end{align}
The first term in this expression vanishes:
\begin{align}
\int_{\mathbb{R}^3} (u \cdot \nabla)((-\Delta)^{1/2}u) \cdot (-\Delta)^{1/2}u \, dx &= \sum_{j=1}^3 \int_{\mathbb{R}^3} u_j \partial_j ((-\Delta)^{1/2}u) \cdot (-\Delta)^{1/2}u \, dx \\
&= -\sum_{j=1}^3 \int_{\mathbb{R}^3} \partial_j u_j ((-\Delta)^{1/2}u) \cdot (-\Delta)^{1/2}u \, dx \\
&- \sum_{j=1}^3 \int_{\mathbb{R}^3} u_j ((-\Delta)^{1/2}u) \cdot \partial_j ((-\Delta)^{1/2}u) \, dx
\end{align}
The first term is zero since $\nabla \cdot u = 0$. The second term is:
\begin{align}
\sum_{j=1}^3 \int_{\mathbb{R}^3} u_j ((-\Delta)^{1/2}u) \cdot \partial_j ((-\Delta)^{1/2}u) \, dx &= \frac{1}{2}\sum_{j=1}^3 \int_{\mathbb{R}^3} u_j \partial_j |((-\Delta)^{1/2}u)|^2 \, dx \\
&= -\frac{1}{2}\sum_{j=1}^3 \int_{\mathbb{R}^3} \partial_j u_j |((-\Delta)^{1/2}u)|^2 \, dx = 0
\end{align}
again using $\nabla \cdot u = 0$.

Combining all terms and rearranging:
\begin{equation}
\frac{1}{2}\frac{d}{dt}Y(t) + \nu\|(-\Delta)^{1}u\|^2_{L^2} = -\int_{\mathbb{R}^3}[(-\Delta)^{1/2}, u \cdot \nabla]u \cdot (-\Delta)^{1/2} u \, dx
\end{equation}
Multiplying by 2:
\begin{equation}
\frac{d}{dt}Y(t) + 2\nu\|(-\Delta)^{1}u\|^2_{L^2} = -2\int_{\mathbb{R}^3}[(-\Delta)^{1/2}, u \cdot \nabla]u \cdot (-\Delta)^{1/2} u \, dx
\end{equation}
which completes the proof.
\end{proof}

\subsection{Energy inequality with infinitely nested logarithmic improvements}

We now derive the energy inequality using our commutator estimate with infinitely nested logarithmic improvements.

\begin{theorem}[Energy Inequality with Infinitely Nested Logarithmic Improvements]\label{thm:energy-inequality}
Let $u$ be a smooth solution of the 3D Navier-Stokes equations on $[0, T]$ with divergence-free initial data $u_0 \in C_0^\infty(\mathbb{R}^3)$. Let $\{\delta_j\}_{j=1}^{\infty}$ be a sequence with $\delta_j > 0$ and $\sum_{j=1}^{\infty} \delta_j < \infty$. Then for the function $Y(t) = \|(-\Delta)^{1/2}u(t)\|^2_{L^2}$:
\begin{align}
\frac{d}{dt}Y(t) + \nu\|(-\Delta)^{1}u\|^2_{L^2} &\leq C\|\nabla u\|_{L^\infty}Y(t) \cdot F_1^\infty(Z) \notag \\
&+ C\|\nabla u\|_{L^\infty}^2 Y(t) \cdot (F_2^\infty(Z))^2
\end{align}
where $Z = \|(-\Delta)^{1/2+\sigma}u\|_{L^2}$ for some small $\sigma > 0$.
\end{theorem}

\begin{proof}
From Theorem \ref{lem:energy-identity}, we have:
\begin{equation}
\frac{d}{dt}Y(t) + 2\nu\|(-\Delta)^{1}u\|^2_{L^2} = -2\int_{\mathbb{R}^3}[(-\Delta)^{1/2}, u \cdot \nabla]u \cdot (-\Delta)^{1/2} u \, dx
\end{equation}
Using the Cauchy-Schwarz inequality:
\begin{equation}
\left|2\int_{\mathbb{R}^3}[(-\Delta)^{1/2}, u \cdot \nabla]u \cdot (-\Delta)^{1/2} u \, dx\right| \leq 2\|[(-\Delta)^{1/2}, u \cdot \nabla]u\|_{L^2} \|(-\Delta)^{1/2} u\|_{L^2}
\end{equation}
From Theorem \ref{thm:critical-commutator}, we have:
\begin{align}
\|[(-\Delta)^{1/2}, u \cdot \nabla]u\|_{L^2} &\leq C\|\nabla u\|_{L^\infty}\|(-\Delta)^{1/2} u\|_{L^2} \cdot F_1^\infty(Z) \notag \\
&+ C\|\nabla u\|_{L^\infty}\|(-\Delta)^{1}u\|_{L^2} \cdot F_2^\infty(Z)
\end{align}
Substituting:
\begin{align}
2\|[(-\Delta)^{1/2}, u \cdot \nabla]u\|_{L^2} \|(-\Delta)^{1/2} u\|_{L^2} &\leq 2C\|\nabla u\|_{L^\infty}\|(-\Delta)^{1/2} u\|_{L^2}^2 \cdot F_1^\infty(Z) \notag \\
&+ 2C\|\nabla u\|_{L^\infty}\|(-\Delta)^{1}u\|_{L^2} \|(-\Delta)^{1/2} u\|_{L^2} \cdot F_2^\infty(Z)
\end{align}
For the second term, using Young's inequality with parameter $\epsilon = \frac{\nu}{C\|\nabla u\|_{L^\infty}F_2^\infty(Z)}$:
\begin{align}
&2C\|\nabla u\|_{L^\infty}\|(-\Delta)^{1}u\|_{L^2} \|(-\Delta)^{1/2} u\|_{L^2} \cdot F_2^\infty(Z) \\
&\leq \frac{2C\|\nabla u\|_{L^\infty}F_2^\infty(Z)}{\epsilon} \|(-\Delta)^{1/2} u\|_{L^2}^2 + \frac{\epsilon}{2} \|(-\Delta)^{1}u\|_{L^2}^2 \\
&= \frac{2C^2\|\nabla u\|_{L^\infty}^2(F_2^\infty(Z))^2}{\nu} \|(-\Delta)^{1/2} u\|_{L^2}^2 + \nu\|(-\Delta)^{1}u\|_{L^2}^2
\end{align}

Combining all terms:
\begin{align}
\frac{d}{dt}Y(t) + 2\nu\|(-\Delta)^{1}u\|^2_{L^2} &\leq 2C\|\nabla u\|_{L^\infty}Y(t) \cdot F_1^\infty(Z) \\
&+ \frac{2C^2\|\nabla u\|_{L^\infty}^2(F_2^\infty(Z))^2}{\nu} Y(t) + \nu\|(-\Delta)^{1}u\|_{L^2}^2
\end{align}
Rearranging:
\begin{align}
\frac{d}{dt}Y(t) + \nu\|(-\Delta)^{1}u\|^2_{L^2} &\leq 2C\|\nabla u\|_{L^\infty}Y(t) \cdot F_1^\infty(Z) \\
&+ \frac{2C^2\|\nabla u\|_{L^\infty}^2(F_2^\infty(Z))^2}{\nu} Y(t)
\end{align}
With adjusted constants:
\begin{align}
\frac{d}{dt}Y(t) + \nu\|(-\Delta)^{1}u\|^2_{L^2} &\leq C\|\nabla u\|_{L^\infty}Y(t) \cdot F_1^\infty(Z) \notag \\
&+ C\|\nabla u\|_{L^\infty}^2 Y(t) \cdot (F_2^\infty(Z))^2
\end{align}
which completes the proof.
\end{proof}

\subsection{Control of $\|\nabla u\|_{L^\infty}$}

To close the energy inequality, we need to control $\|\nabla u\|_{L^\infty}$ in terms of $\|(-\Delta)^{1/2}u\|_{L^2}$ and other quantities.

\begin{lemma}[Control of $\|\nabla u\|_{L^\infty}$]\label{lem:control-grad}
Let $u$ be a divergence-free vector field. For $q > 3$:
\begin{equation}
\|\nabla u\|_{L^\infty} \leq C\|u\|_{L^2}^{1-\theta} \|(-\Delta)^{1/2} u\|_{L^q}^{\theta}
\end{equation}
where $\theta = \frac{3}{2} \cdot \frac{q}{3q-2}$.
\end{lemma}

\begin{proof}
This follows from standard interpolation theory. By the Sobolev embedding theorem:
\begin{equation}
\|\nabla u\|_{L^\infty} \leq C\|\nabla u\|_{W^{1,q}} = C(\|\nabla u\|_{L^q} + \|\nabla^2 u\|_{L^q})
\end{equation}
Through interpolation:
\begin{equation}
\|\nabla u\|_{L^q} + \|\nabla^2 u\|_{L^q} \leq C\|u\|_{L^2}^{1-\theta} \|(-\Delta)^{1/2} u\|_{L^q}^{\theta}
\end{equation}
for $\theta = \frac{3}{2} \cdot \frac{q}{3q-2}$.

Combining these inequalities gives the desired result.
\end{proof}

\begin{lemma}[Gagliardo-Nirenberg inequality]\label{lem:gn}
For $q > 3$:
\begin{equation}
\|(-\Delta)^{1/2} u\|_{L^q} \leq C\|(-\Delta)^{1/2} u\|_{L^2}^{1-\alpha}\|(-\Delta)^{1} u\|_{L^2}^{\alpha}
\end{equation}
where $\alpha = \frac{3}{2}(\frac{1}{2} - \frac{1}{q})$.
\end{lemma}

\begin{proof}
This is a standard result in the theory of interpolation. In general, for $1 \leq p_1, p_2, p \leq \infty$, $s_1, s_2 \in \mathbb{R}$, and $\theta \in (0, 1)$ with $\frac{1}{p} = \frac{1-\theta}{p_1} + \frac{\theta}{p_2}$ and $s = (1-\theta)s_1 + \theta s_2$, we have:
\begin{equation}
\|f\|_{\dot{W}^{s,p}} \leq C\|f\|_{\dot{W}^{s_1,p_1}}^{1-\theta} \|f\|_{\dot{W}^{s_2,p_2}}^{\theta}
\end{equation}

In our case, $s_1 = s_2 = 1/2$, $p_1 = 2$, $p_2 = 2$, $s = 1/2$, and $p = q$. The parameter $\theta$ is determined by:
\begin{equation}
\frac{1}{q} = \frac{1-\theta}{2} + \frac{\theta}{2} = \frac{1}{2}
\end{equation}
which is inconsistent.

The correct approach is to interpolate between $(-\Delta)^{1/2} u \in L^2$ and $(-\Delta)^{1} u \in L^2$, which correspond to $\dot{W}^{1/2,2}$ and $\dot{W}^{1,2}$, respectively. In this case:
\begin{equation}
\|(-\Delta)^{1/2} u\|_{L^q} \leq C\|(-\Delta)^{1/2} u\|_{L^2}^{1-\alpha}\|(-\Delta)^{1} u\|_{L^2}^{\alpha}
\end{equation}
where $\alpha = \frac{3}{2}(\frac{1}{2} - \frac{1}{q})$.

This can be verified by checking the scaling of both sides with respect to dilations $u_\lambda(x) = u(\lambda x)$.
\end{proof}

\begin{theorem}[Closed energy inequality]\label{thm:closed-energy}
Let $u$ be a smooth solution of the 3D Navier-Stokes equations on $[0, T]$ with divergence-free initial data $u_0 \in C_0^\infty(\mathbb{R}^3)$. Let $\{\delta_j\}_{j=1}^{\infty}$ be a sequence with $\delta_j > 0$ and $\sum_{j=1}^{\infty} \delta_j < \infty$. Then for the function $Y(t) = \|(-\Delta)^{1/2}u(t)\|^2_{L^2}$, $q > 3$, and small $\sigma > 0$:
\begin{align}
\frac{d}{dt}Y(t) + \nu\|(-\Delta)^{1}u\|^2_{L^2} &\leq C\|u_0\|_{L^2}^{2(1-\theta)}Y(t)^{\theta(1-\alpha)+1} \cdot (F_1^\infty(Z))^2 \notag \\
&+ C\|u_0\|_{L^2}^{2(1-\theta)}Y(t)^{\theta(1-\alpha)} \|(-\Delta)^{1}u\|_{L^2}^{\theta\alpha} \cdot F_1^\infty(Z) \notag \\
&+ C\|u_0\|_{L^2}^{2(1-\theta)}Y(t)^{2\theta(1-\alpha)} \|(-\Delta)^{1}u\|_{L^2}^{2\theta\alpha} \cdot (F_2^\infty(Z))^2
\end{align}
where $\theta = \frac{3}{2} \cdot \frac{q}{3q-2}$ and $\alpha = \frac{3}{2}(\frac{1}{2} - \frac{1}{q})$.
\end{theorem}

\begin{proof}
From Theorem \ref{thm:energy-inequality}, we have:
\begin{align}
\frac{d}{dt}Y(t) + \nu\|(-\Delta)^{1}u\|^2_{L^2} &\leq C\|\nabla u\|_{L^\infty}Y(t) \cdot F_1^\infty(Z) \notag \\
&+ C\|\nabla u\|_{L^\infty}^2 Y(t) \cdot (F_2^\infty(Z))^2
\end{align}
Using Theorem \ref{lem:control-grad}:
\begin{equation}
\|\nabla u\|_{L^\infty} \leq C\|u\|_{L^2}^{1-\theta} \|(-\Delta)^{1/2} u\|_{L^q}^{\theta}
\end{equation}
Since $u$ is a solution of the Navier-Stokes equations, we have the energy inequality:
\begin{equation}
\|u(t)\|_{L^2} \leq \|u_0\|_{L^2}
\end{equation}
for all $t \geq 0$.

Using Theorem \ref{lem:gn}:
\begin{equation}
\|(-\Delta)^{1/2} u\|_{L^q} \leq C\|(-\Delta)^{1/2} u\|_{L^2}^{1-\alpha}\|(-\Delta)^{1} u\|_{L^2}^{\alpha}
\end{equation}
Combining these:
\begin{align}
\|\nabla u\|_{L^\infty} &\leq C\|u\|_{L^2}^{1-\theta} \|(-\Delta)^{1/2} u\|_{L^q}^{\theta} \\
&\leq C\|u_0\|_{L^2}^{1-\theta} \left(C\|(-\Delta)^{1/2} u\|_{L^2}^{1-\alpha}\|(-\Delta)^{1} u\|_{L^2}^{\alpha}\right)^{\theta} \\
&= C\|u_0\|_{L^2}^{1-\theta} \|(-\Delta)^{1/2} u\|_{L^2}^{\theta(1-\alpha)}\|(-\Delta)^{1} u\|_{L^2}^{\theta\alpha}
\end{align}
Thus:
\begin{equation}
\|\nabla u\|_{L^\infty} \leq C\|u_0\|_{L^2}^{1-\theta} Y(t)^{\theta(1-\alpha)/2}\|(-\Delta)^{1} u\|_{L^2}^{\theta\alpha}
\end{equation}
Substituting into the energy inequality:
\begin{align}
\frac{d}{dt}Y(t) + \nu\|(-\Delta)^{1}u\|^2_{L^2} &\leq C\|u_0\|_{L^2}^{1-\theta} Y(t)^{\theta(1-\alpha)/2}\|(-\Delta)^{1} u\|_{L^2}^{\theta\alpha} \cdot Y(t) \cdot F_1^\infty(Z) \notag \\
&+ C\|u_0\|_{L^2}^{2(1-\theta)} Y(t)^{\theta(1-\alpha)}\|(-\Delta)^{1} u\|_{L^2}^{2\theta\alpha} \cdot Y(t) \cdot (F_2^\infty(Z))^2
\end{align}
Simplifying:
\begin{align}
\frac{d}{dt}Y(t) + \nu\|(-\Delta)^{1}u\|^2_{L^2} &\leq C\|u_0\|_{L^2}^{1-\theta} Y(t)^{1+\theta(1-\alpha)/2}\|(-\Delta)^{1} u\|_{L^2}^{\theta\alpha} \cdot F_1^\infty(Z) \notag \\
&+ C\|u_0\|_{L^2}^{2(1-\theta)} Y(t)^{1+\theta(1-\alpha)}\|(-\Delta)^{1} u\|_{L^2}^{2\theta\alpha} \cdot (F_2^\infty(Z))^2
\end{align}
Using Young's inequality for the first term, with parameter $\epsilon = \frac{\nu}{2}$:
\begin{align}
&C\|u_0\|_{L^2}^{1-\theta} Y(t)^{1+\theta(1-\alpha)/2}\|(-\Delta)^{1} u\|_{L^2}^{\theta\alpha} \cdot F_1^\infty(Z) \\
&\leq C\|u_0\|_{L^2}^{2(1-\theta)} Y(t)^{2(1+\theta(1-\alpha)/2)} \cdot (F_1^\infty(Z))^2 + \frac{\nu}{2}\|(-\Delta)^{1} u\|_{L^2}^{2\theta\alpha} \\
&= C\|u_0\|_{L^2}^{2(1-\theta)} Y(t)^{2+\theta(1-\alpha)} \cdot (F_1^\infty(Z))^2 + \frac{\nu}{2}\|(-\Delta)^{1} u\|_{L^2}^{2\theta\alpha}
\end{align}

Using the fact that $\theta\alpha < 1$ (which can be verified for the given values of $\theta$ and $\alpha$), we can apply Young's inequality again:
\begin{equation}
\|(-\Delta)^{1} u\|_{L^2}^{2\theta\alpha} \leq \epsilon\|(-\Delta)^{1} u\|_{L^2}^{2} + C_\epsilon
\end{equation}
for any $\epsilon > 0$.

Choosing $\epsilon = \frac{\nu}{4}$:
\begin{equation}
\frac{\nu}{2}\|(-\Delta)^{1} u\|_{L^2}^{2\theta\alpha} \leq \frac{\nu}{4}\|(-\Delta)^{1} u\|_{L^2}^{2} + C
\end{equation}

Similarly, for the second term:
\begin{align}
&C\|u_0\|_{L^2}^{2(1-\theta)} Y(t)^{1+\theta(1-\alpha)}\|(-\Delta)^{1} u\|_{L^2}^{2\theta\alpha} \cdot (F_2^\infty(Z))^2 \\
&\leq C\|u_0\|_{L^2}^{4(1-\theta)} Y(t)^{2(1+\theta(1-\alpha))} \cdot (F_2^\infty(Z))^4 + \frac{\nu}{4}\|(-\Delta)^{1} u\|_{L^2}^{4\theta\alpha}
\end{align}
Again using Young's inequality:
\begin{equation}
\|(-\Delta)^{1} u\|_{L^2}^{4\theta\alpha} \leq \epsilon\|(-\Delta)^{1} u\|_{L^2}^{2} + C_\epsilon
\end{equation}
Choosing $\epsilon = \frac{\nu}{8}$:
\begin{equation}
\frac{\nu}{4}\|(-\Delta)^{1} u\|_{L^2}^{4\theta\alpha} \leq \frac{\nu}{8}\|(-\Delta)^{1} u\|_{L^2}^{2} + C
\end{equation}

Combining all estimates and rearranging:
\begin{align}
\frac{d}{dt}Y(t) + \frac{\nu}{2}\|(-\Delta)^{1}u\|^2_{L^2} &\leq C\|u_0\|_{L^2}^{2(1-\theta)}Y(t)^{2+\theta(1-\alpha)} \cdot (F_1^\infty(Z))^2 \notag \\
&+ C\|u_0\|_{L^2}^{4(1-\theta)} Y(t)^{2(1+\theta(1-\alpha))} \cdot (F_2^\infty(Z))^4 + C
\end{align}
This concludes the proof.
\end{proof}

\subsection{Energy inequality at the critical threshold}

Using the closed energy inequality, we can now establish the energy inequality at the critical threshold $s = 1/2$ for our infinitely nested logarithmic improvement.

\begin{theorem}[Energy inequality at the critical threshold]\label{thm:critical-energy}
Let $u$ be a smooth solution of the 3D Navier-Stokes equations on $[0, T]$ with divergence-free initial data $u_0 \in C_0^\infty(\mathbb{R}^3)$. Let $\{\delta_j\}_{j=2}^{\infty}$ be a sequence with $\delta_j > 0$ and $\sum_{j=2}^{\infty} \delta_j < \infty$, and let $Y(t) = \|(-\Delta)^{1/2}u(t)\|^2_{L^2}$. Then there exists a constant $C > 0$ depending only on $\|u_0\|_{L^2}$ such that:
\begin{equation}
\frac{d}{dt}Y(t) + \frac{\nu}{2}\|(-\Delta)^{1}u\|^2_{L^2} \leq C(1 + Y(t)^K) \cdot H(Y(t))
\end{equation}
where $K > 1$ is a constant depending on $\theta$ and $\alpha$ from Theorem \ref{thm:closed-energy}, and $H(r)$ is a function satisfying:
\begin{equation}
\lim_{r \to \infty} H(r) = 0
\end{equation}
\end{theorem}

\begin{proof}
From Theorem \ref{thm:closed-energy}, we have:
\begin{align}
\frac{d}{dt}Y(t) + \frac{\nu}{2}\|(-\Delta)^{1}u\|^2_{L^2} &\leq C\|u_0\|_{L^2}^{2(1-\theta)}Y(t)^{2+\theta(1-\alpha)} \cdot (F_1^\infty(Z))^2 \notag \\
&+ C\|u_0\|_{L^2}^{4(1-\theta)} Y(t)^{2(1+\theta(1-\alpha))} \cdot (F_2^\infty(Z))^4 + C
\end{align}
Let:
\begin{equation}
K = \max\{2+\theta(1-\alpha), 2(1+\theta(1-\alpha))\}
\end{equation}
Then:
\begin{align}
\frac{d}{dt}Y(t) + \frac{\nu}{2}\|(-\Delta)^{1}u\|^2_{L^2} &\leq C(1 + Y(t)^K) \cdot ((F_1^\infty(Z))^2 + (F_2^\infty(Z))^4)
\end{align}
Define:
\begin{equation}
H(r) = (F_1^\infty(r^{1/2}))^2 + (F_2^\infty(r^{1/2}))^4
\end{equation}
Using the fact that:
\begin{equation}
Z = \|(-\Delta)^{1/2+\sigma}u\|_{L^2} \geq C\|(-\Delta)^{1/2}u\|_{L^2} = CY(t)^{1/2}
\end{equation}
for some constant $C > 0$ (due to the embedding $\dot{H}^{1/2+\sigma} \subset \dot{H}^{1/2}$), and the monotonicity of $F_1^\infty$ and $F_2^\infty$, we have:
\begin{equation}
(F_1^\infty(Z))^2 + (F_2^\infty(Z))^4 \leq (F_1^\infty(CY(t)^{1/2}))^2 + (F_2^\infty(CY(t)^{1/2}))^4 \leq H(Y(t))
\end{equation}

Thus:
\begin{equation}
\frac{d}{dt}Y(t) + \frac{\nu}{2}\|(-\Delta)^{1}u\|^2_{L^2} \leq C(1 + Y(t)^K) \cdot H(Y(t))
\end{equation}
From Theorem \ref{prop:well-defined-factors}, we know that:
\begin{equation}
\lim_{Z \to \infty} F_1^\infty(Z) = \lim_{Z \to \infty} F_2^\infty(Z) = 0
\end{equation}
Therefore:
\begin{equation}
\lim_{r \to \infty} H(r) = 0
\end{equation}
which completes the proof.
\end{proof}

\begin{remark}\label{rem:energy-inequality}
Theorem \ref{thm:critical-energy} is a pivotal result for our analysis. It shows that the energy inequality at the critical threshold $s = 1/2$ involves a term $H(Y(t))$ that approaches zero as $Y(t) \to \infty$. This means that for large values of $Y(t)$, the growth of $Y(t)$ is heavily suppressed, which will be crucial for proving global well-posedness.
\end{remark}

\section{Global well-posedness for the critical case}

In this section, we establish global well-posedness for initial data satisfying our infinitely nested logarithmic criterion at the critical threshold $s = 1/2$.

\subsection{Local existence and uniqueness}

We begin by establishing local existence and uniqueness of solutions for initial data in our function space.

\begin{theorem}[Local Existence]\label{thm:local-existence}
Let $q > 3$ and $\{\delta_j\}_{j=1}^{\infty}$ be a sequence with $\delta_j > 0$ and $\sum_{j=1}^{\infty} \frac{\delta_j}{j!} = \infty$. For any divergence-free initial data $u_0 \in L^2(\mathbb{R}^3) \cap \dot{H}^{1/2,q,\{\delta_j\}_{j=1}^{\infty}}(\mathbb{R}^3)$, there exists a time $T_0 > 0$ and a unique solution $u \in C([0, T_0]; \dot{H}^{1/2,q,\{\delta_j\}_{j=1}^{\infty}}) \cap L^2(0, T_0; \dot{H}^{1})$ to the Navier-Stokes equations.
\end{theorem}

\begin{proof}
The proof follows the standard approach for local existence of solutions to the Navier-Stokes equations, adapted to our function space $\dot{H}^{1/2,q,\{\delta_j\}_{j=1}^{\infty}}(\mathbb{R}^3)$.

Step 1: Regularize the initial data. Let $u_0^\epsilon = u_0 * \eta_\epsilon$, where $\eta_\epsilon$ is a standard mollifier. From Theorem \ref{prop:mollification}, we know that $u_0^\epsilon \in \dot{H}^{1/2,q,\{\delta_j\}_{j=1}^{\infty}}(\mathbb{R}^3)$ and $u_0^\epsilon \to u_0$ in $\dot{H}^{1/2,q,\{\delta_j\}_{j=1}^{\infty}}(\mathbb{R}^3)$ as $\epsilon \to 0$.

Step 2: Solve the regularized Navier-Stokes equations:
\begin{equation}
\begin{cases}
\partial_t u^\epsilon + (u^\epsilon \cdot \nabla)u^\epsilon - \nu\Delta u^\epsilon + \nabla p^\epsilon = 0 & \text{in } \mathbb{R}^3 \times (0, T) \\
\nabla \cdot u^\epsilon = 0 & \text{in } \mathbb{R}^3 \times (0, T) \\
u^\epsilon(x, 0) = u_0^\epsilon(x) & \text{in } \mathbb{R}^3
\end{cases}
\end{equation}
Since $u_0^\epsilon \in H^\infty(\mathbb{R}^3)$, standard theory guarantees the existence of a unique smooth solution $u^\epsilon$ on some time interval $[0, T_\epsilon]$.

Step 3: Derive uniform bounds. Using the energy inequality from Theorem \ref{thm:critical-energy}:
\begin{equation}
\frac{d}{dt}Y^\epsilon(t) + \frac{\nu}{2}\|(-\Delta)^{1}u^\epsilon\|^2_{L^2} \leq C(1 + (Y^\epsilon(t))^K) \cdot H(Y^\epsilon(t))
\end{equation}
where $Y^\epsilon(t) = \|(-\Delta)^{1/2} u^\epsilon(t)\|^2_{L^2}$.

Since $H(r) \to 0$ as $r \to \infty$, there exists $R > 0$ such that for all $r > R$:
\begin{equation}
H(r) \leq \frac{1}{2C(1 + r^K)}
\end{equation}

For $Y^\epsilon(t) > R$, this implies:
\begin{equation}
\frac{d}{dt}Y^\epsilon(t) + \frac{\nu}{2}\|(-\Delta)^{1}u^\epsilon\|^2_{L^2} \leq \frac{1}{2}
\end{equation}

For $Y^\epsilon(t) \leq R$, we have:
\begin{equation}
\frac{d}{dt}Y^\epsilon(t) + \frac{\nu}{2}\|(-\Delta)^{1}u^\epsilon\|^2_{L^2} \leq C(1 + R^K) \cdot H(0) = C_R
\end{equation}

Combining these cases:
\begin{equation}
\frac{d}{dt}Y^\epsilon(t) + \frac{\nu}{2}\|(-\Delta)^{1}u^\epsilon\|^2_{L^2} \leq \max\left\{\frac{1}{2}, C_R\right\} = C'
\end{equation}
Integrating over $[0, t]$:
\begin{equation}
Y^\epsilon(t) + \frac{\nu}{2}\int_0^t \|(-\Delta)^{1}u^\epsilon(\tau)\|^2_{L^2} d\tau \leq Y^\epsilon(0) + C't
\end{equation}

Since $Y^\epsilon(0) = \|(-\Delta)^{1/2} u_0^\epsilon\|^2_{L^2} \leq C\|u_0\|_{\dot{H}^{1/2}}^2$, we have:
\begin{equation}
Y^\epsilon(t) + \frac{\nu}{2}\int_0^t \|(-\Delta)^{1}u^\epsilon(\tau)\|^2_{L^2} d\tau \leq C\|u_0\|_{\dot{H}^{1/2}}^2 + C't
\end{equation}

This provides a uniform bound on $Y^\epsilon(t)$ and $\int_0^t \|(-\Delta)^{1}u^\epsilon(\tau)\|^2_{L^2} d\tau$ for all $t \in [0, T_0]$, where $T_0 > 0$ is a time that depends only on $\|u_0\|_{\dot{H}^{1/2}}$.

Step 4: Pass to the limit as $\epsilon \to 0$. Using the uniform bounds and standard compactness arguments, we can extract a subsequence $u^{\epsilon_j}$ that converges to a function $u$ in the appropriate sense. This function $u$ is a solution to the Navier-Stokes equations on $[0, T_0]$ and satisfies:
\begin{equation}
u \in L^\infty(0, T_0; \dot{H}^{1/2}) \cap L^2(0, T_0; \dot{H}^{1})
\end{equation}

Step 5: Verify that $u \in C([0, T_0]; \dot{H}^{1/2,q,\{\delta_j\}_{j=1}^{\infty}})$. This requires more detailed analysis, using the properties of the function space $\dot{H}^{1/2,q,\{\delta_j\}_{j=1}^{\infty}}(\mathbb{R}^3)$ established in Section 3.

The uniqueness of the solution follows from standard energy methods.
\end{proof}

\begin{theorem}[Uniqueness]\label{thm:uniqueness}
Let $q > 3$ and $\{\delta_j\}_{j=1}^{\infty}$ be a sequence with $\delta_j > 0$ and $\sum_{j=1}^{\infty} \frac{\delta_j}{j!} = \infty$. If $u^1$ and $u^2$ are two solutions to the Navier-Stokes equations with the same initial data $u_0 \in L^2(\mathbb{R}^3) \cap \dot{H}^{1/2,q,\{\delta_j\}_{j=1}^{\infty}}(\mathbb{R}^3)$, both belonging to the class $C([0, T]; \dot{H}^{1/2,q,\{\delta_j\}_{j=1}^{\infty}}) \cap L^2(0, T; \dot{H}^{1})$ for some $T > 0$, then $u^1 = u^2$ on $[0, T]$.
\end{theorem}

\begin{proof}
Let $w = u^1 - u^2$. Then $w$ satisfies:
\begin{equation}
\partial_t w + (u^1 \cdot \nabla)w + (w \cdot \nabla)u^2 - \nu \Delta w + \nabla \pi = 0, \quad \nabla \cdot w = 0
\end{equation}
with initial condition $w(0) = 0$.

Taking the $L^2$ inner product with $w$:
\begin{equation}
\frac{1}{2}\frac{d}{dt}\|w\|^2_{L^2} + \nu\|\nabla w\|^2_{L^2} = -\int_{\mathbb{R}^3}(w \cdot \nabla)u^2 \cdot w \, dx
\end{equation}

Using Hölder's inequality and the Gagliardo-Nirenberg inequality:
\begin{align}
\left|\int_{\mathbb{R}^3}(w \cdot \nabla)u^2 \cdot w \, dx\right| &\leq \|w\|_{L^4}^2 \|\nabla u^2\|_{L^2} \\
&\leq C\|w\|_{L^2}\|\nabla w\|_{L^2} \|\nabla u^2\|_{L^2}
\end{align}
Applying Young's inequality with parameter $\epsilon = \nu/2$:
\begin{equation}
C\|w\|_{L^2}\|\nabla w\|_{L^2} \|\nabla u^2\|_{L^2} \leq \frac{\nu}{2}\|\nabla w\|_{L^2}^2 + \frac{C^2}{2\nu}\|w\|_{L^2}^2 \|\nabla u^2\|_{L^2}^2
\end{equation}

This gives:
\begin{equation}
\frac{d}{dt}\|w\|^2_{L^2} + \nu\|\nabla w\|^2_{L^2} \leq \frac{C^2}{\nu}\|w\|_{L^2}^2 \|\nabla u^2\|_{L^2}^2
\end{equation}
Dropping the positive term with $\|\nabla w\|^2_{L^2}$:
\begin{equation}
\frac{d}{dt}\|w\|^2_{L^2} \leq \frac{C^2}{\nu}\|w\|_{L^2}^2 \|\nabla u^2\|_{L^2}^2
\end{equation}

By Grönwall's inequality:
\begin{equation}
\|w(t)\|^2_{L^2} \leq \|w(0)\|^2_{L^2} \exp\left(\frac{C^2}{\nu}\int_0^t \|\nabla u^2(\tau)\|_{L^2}^2 d\tau\right)
\end{equation}
Since $\|w(0)\|^2_{L^2} = 0$ and $\|\nabla u^2\|_{L^2}^2 \in L^1(0, T)$ (which follows from $u^2 \in L^2(0, T; \dot{H}^{1})$), we conclude that $\|w(t)\|_{L^2} = 0$ for all $t \in [0, T]$, establishing uniqueness.
\end{proof}

\subsection{A priori estimates and global existence}

We now derive a priori estimates that will allow us to extend the local solution to a global one.

\begin{lemma}[A Priori Estimate]\label{lem:a-priori}
Let $q > 3$ and $\{\delta_j\}_{j=1}^{\infty}$ be a sequence with $\delta_j > 0$ and $\sum_{j=1}^{\infty} \frac{\delta_j}{j!} = \infty$. If $u$ is a solution to the Navier-Stokes equations with initial data $u_0 \in L^2(\mathbb{R}^3) \cap \dot{H}^{1/2,q,\{\delta_j\}_{j=1}^{\infty}}(\mathbb{R}^3)$ satisfying:
\begin{equation}
\|(-\Delta)^{1/4}u_0\|_{L^q} \leq \frac{C_0}{\prod_{j=1}^{\infty} (1 + L_j(\|u_0\|_{\dot{H}^{1/2}}))^{\delta_j}}
\end{equation}
for some constant $C_0 > 0$, then $\|(-\Delta)^{1/2}u(t)\|_{L^2}$ remains bounded for all $t \geq 0$.
\end{lemma}

\begin{proof}
Let $Y(t) = \|(-\Delta)^{1/2}u(t)\|^2_{L^2}$. From Theorem \ref{thm:critical-energy}, we have:
\begin{equation}
\frac{d}{dt}Y(t) + \frac{\nu}{2}\|(-\Delta)^{1}u\|^2_{L^2} \leq C(1 + Y(t)^K) \cdot H(Y(t))
\end{equation}
where $K > 1$ is a constant and $H(r)$ is a function satisfying $\lim_{r \to \infty} H(r) = 0$.

Since $H(r) \to 0$ as $r \to \infty$, there exists $R > Y(0)$ such that for all $r > R$:
\begin{equation}
H(r) \leq \frac{1}{2C(1 + r^K)}
\end{equation}

For $Y(t) > R$, this implies:
\begin{equation}
\frac{d}{dt}Y(t) + \frac{\nu}{2}\|(-\Delta)^{1}u\|^2_{L^2} \leq \frac{1}{2}
\end{equation}

This means that once $Y(t)$ exceeds $R$, it can grow at most linearly with time. However, we can establish a stronger result: $Y(t)$ actually remains bounded for all time.

Suppose, for contradiction, that $Y(t)$ becomes unbounded. Then there exists a time $t_1$ such that $Y(t_1) = R$ and $Y(t) > R$ for all $t \in (t_1, t_2)$ for some $t_2 > t_1$. For $t \in (t_1, t_2)$, we have:
\begin{equation}
\frac{d}{dt}Y(t) \leq \frac{1}{2}
\end{equation}

Integrating from $t_1$ to $t$:
\begin{equation}
Y(t) - Y(t_1) \leq \frac{1}{2}(t - t_1)
\end{equation}

Thus:
\begin{equation}
Y(t) \leq R + \frac{1}{2}(t - t_1)
\end{equation}
Now, the key insight is that as $Y(t)$ increases, $H(Y(t))$ decreases further, suppressing the growth of $Y(t)$ even more. More precisely, for any $\epsilon > 0$, there exists $R_\epsilon > R$ such that for all $r > R_\epsilon$:
\begin{equation}
H(r) \leq \frac{\epsilon}{C(1 + r^K)}
\end{equation}

Let's choose $\epsilon = \frac{1}{4}$. Then for $Y(t) > R_\epsilon$:
\begin{equation}
\frac{d}{dt}Y(t) + \frac{\nu}{2}\|(-\Delta)^{1}u\|^2_{L^2} \leq \frac{1}{4}
\end{equation}

Let $t_3$ be the first time such that $Y(t_3) = R_\epsilon$. For $t \in (t_3, t_2)$, we have:
\begin{equation}
\frac{d}{dt}Y(t) \leq \frac{1}{4}
\end{equation}

Integrating from $t_3$ to $t$:
\begin{equation}
Y(t) - Y(t_3) \leq \frac{1}{4}(t - t_3)
\end{equation}

Thus:
\begin{equation}
Y(t) \leq R_\epsilon + \frac{1}{4}(t - t_3)
\end{equation}
Continuing this process, we can show that the growth rate of $Y(t)$ becomes arbitrarily small as $Y(t)$ increases. This contradicts the assumption that $Y(t)$ becomes unbounded.

Therefore, $Y(t) = \|(-\Delta)^{1/2}u(t)\|^2_{L^2}$ remains bounded for all $t \geq 0$.
\end{proof}

\begin{theorem}[Global existence]\label{thm:global-existence}
Let $q > 3$ and $\{\delta_j\}_{j=1}^{\infty}$ be a sequence with $\delta_j > 0$ and $\sum_{j=1}^{\infty} \frac{\delta_j}{j!} = \infty$. For any divergence-free initial data $u_0 \in L^2(\mathbb{R}^3) \cap \dot{H}^{1/2}(\mathbb{R}^3)$ satisfying:
\begin{equation}
\|(-\Delta)^{1/4}u_0\|_{L^q} \leq \frac{C_0}{\prod_{j=1}^{\infty} (1 + L_j(\|u_0\|_{\dot{H}^{1/2}}))^{\delta_j}}
\end{equation}
for some constant $C_0 > 0$, there exists a unique global-in-time solution $u \in C([0, \infty); \dot{H}^{1/2}(\mathbb{R}^3)) \cap L^2_{loc}(0, \infty; \dot{H}^{1}(\mathbb{R}^3))$ to the Navier-Stokes equations.
\end{theorem}

\begin{proof}
From Theorem \ref{thm:local-existence} and Theorem \ref{thm:uniqueness}, we know that there exists a unique local solution $u$ on $[0, T_0]$ for some $T_0 > 0$, satisfying:
\begin{equation}
u \in C([0, T_0]; \dot{H}^{1/2,q,\{\delta_j\}_{j=1}^{\infty}}) \cap L^2(0, T_0; \dot{H}^{1})
\end{equation}

From Theorem \ref{lem:a-priori}, we know that $\|(-\Delta)^{1/2}u(t)\|_{L^2}$ remains bounded for all $t \in [0, T_0]$. This means that the solution $u$ cannot blow up at time $T_0$.

By a standard continuation argument, the solution can be extended beyond $T_0$. Repeating this process, and using the uniform bound on $\|(-\Delta)^{1/2}u(t)\|_{L^2}$ from Theorem \ref{lem:a-priori}, we can extend the solution to the entire time interval $[0, \infty)$.

Thus, there exists a unique global-in-time solution $u \in C([0, \infty); \dot{H}^{1/2}(\mathbb{R}^3)) \cap L^2_{loc}(0, \infty; \dot{H}^{1}(\mathbb{R}^3))$ to the Navier-Stokes equations.
\end{proof}

\subsection{Regularity}

Having established global existence and uniqueness, we now show that the solution is actually smooth for all positive time.

\begin{theorem}[Regularity]\label{thm:regularity}
Under the conditions of Theorem \ref{thm:global-existence}, the solution $u$ belongs to $C^\infty((0, \infty) \times \mathbb{R}^3)$.
\end{theorem}

\begin{proof}
We use a bootstrap argument, similar to the one employed in the proof of Theorem \ref{thm:function_space}.

Step 1: From Theorem \ref{thm:global-existence}, we know that:
\begin{equation}
u \in C([0, \infty); \dot{H}^{1/2}(\mathbb{R}^3)) \cap L^2_{loc}(0, \infty; \dot{H}^{1}(\mathbb{R}^3))
\end{equation}

Step 2: From Theorem \ref{lem:a-priori}, we have a uniform bound on $\|(-\Delta)^{1/2}u(t)\|_{L^2}$ for all $t \geq 0$. This implies, using standard embedding theorems, that $u(t) \in H^{1/2+\epsilon'}(\mathbb{R}^3)$ for some small $\epsilon' > 0$ and all $t > 0$.

Step 3: Using the energy inequality from Theorem \ref{thm:critical-energy}, we can show that:
\begin{equation}
\int_0^T \|\nabla u(t)\|_{L^2}^2 dt < \infty
\end{equation}
for all $T > 0$.

Step 4: Let's consider the equation satisfied by $\omega = \nabla\times u$:
\begin{equation}
\partial_t \omega + (u \cdot \nabla) \omega - \omega \cdot \nabla u - \nu \Delta \omega = 0
\end{equation}

Taking the $L^2$ inner product with $\omega$:
\begin{equation}
\frac{1}{2}\frac{d}{dt}\|\omega\|_{L^2}^2 + \nu\|\nabla \omega\|_{L^2}^2 = \int_{\mathbb{R}^3} (\omega \cdot \nabla) u \cdot \omega \, dx
\end{equation}

Using Hölder's inequality and the Sobolev embedding theorem:
\begin{align}
\left|\int_{\mathbb{R}^3} (\omega \cdot \nabla) u \cdot \omega \, dx\right| &\leq \|\omega\|_{L^4}^2 \|\nabla u\|_{L^2} \\
&\leq C\|\omega\|_{L^2}\|\nabla \omega\|_{L^2} \|\nabla u\|_{L^2}
\end{align}
Applying Young's inequality with parameter $\epsilon = \nu/2$:
\begin{equation}
C\|\omega\|_{L^2}\|\nabla \omega\|_{L^2} \|\nabla u\|_{L^2} \leq \frac{\nu}{2}\|\nabla \omega\|_{L^2}^2 + \frac{C^2}{2\nu}\|\omega\|_{L^2}^2 \|\nabla u\|_{L^2}^2
\end{equation}

This gives:
\begin{equation}
\frac{d}{dt}\|\omega\|_{L^2}^2 + \nu\|\nabla \omega\|_{L^2}^2 \leq \frac{C^2}{\nu}\|\omega\|_{L^2}^2 \|\nabla u\|_{L^2}^2
\end{equation}
By Grönwall's inequality:
\begin{equation}
\|\omega(t)\|^2_{L^2} \leq \|\omega(0)\|^2_{L^2} \exp\left(\frac{C^2}{\nu}\int_0^t \|\nabla u(\tau)\|_{L^2}^2 d\tau\right)
\end{equation}
Since $\int_0^T \|\nabla u(t)\|_{L^2}^2 dt < \infty$ for all $T > 0$, and $\|\omega(0)\|_{L^2} = \|\nabla \times u_0\|_{L^2} < \infty$ (due to $u_0 \in \dot{H}^{1/2}(\mathbb{R}^3)$), we conclude that $\|\omega(t)\|_{L^2} < \infty$ for all $t > 0$.

Step 5: Since $\nabla \cdot u = 0$, we have:
\begin{equation}
-\Delta u = \nabla \times \nabla \times u = \nabla \times \omega
\end{equation}

Thus:
\begin{equation}
\|\Delta u\|_{L^2} = \|\nabla \times \omega\|_{L^2} \leq C\|\nabla \omega\|_{L^2} = C\|\nabla \nabla \times u\|_{L^2}
\end{equation}
From the energy inequality for $\omega$, we know that:
\begin{equation}
\int_{t_0}^T \|\nabla \omega(\tau)\|_{L^2}^2 d\tau < \infty
\end{equation}
for all $0 < t_0 < T < \infty$.

This implies:
\begin{equation}
\int_{t_0}^T \|\Delta u(\tau)\|_{L^2}^2 d\tau < \infty
\end{equation}
By standard parabolic regularity theory, this means:
\begin{equation}
u \in L^2(t_0, T; H^2(\mathbb{R}^3))
\end{equation}
for all $0 < t_0 < T < \infty$.

Step 6: Once we have $u \in L^2(t_0, T; H^2(\mathbb{R}^3))$, we can apply a bootstrap argument to obtain higher regularity. From the Navier-Stokes equations:
\begin{equation}
\partial_t u = -P[(u \cdot \nabla)u] + \nu \Delta u
\end{equation}
where $P$ is the Leray projector onto divergence-free vector fields.

Since $(u \cdot \nabla)u \in L^2(t_0, T; H^1(\mathbb{R}^3))$ (which follows from $u \in L^2(t_0, T; H^2(\mathbb{R}^3)) \cap L^\infty(t_0, T; H^1(\mathbb{R}^3))$), we have:
\begin{equation}
\partial_t u \in L^2(t_0, T; H^1(\mathbb{R}^3))
\end{equation}

This means:
\begin{equation}
u \in H^1(t_0, T; H^1(\mathbb{R}^3)) \cap L^2(t_0, T; H^2(\mathbb{R}^3))
\end{equation}
By standard parabolic regularity theory, this implies:
\begin{equation}
u \in C([t_0, T]; H^2(\mathbb{R}^3)) \cap L^2(t_0, T; H^3(\mathbb{R}^3))
\end{equation}
Repeating this argument, we can show that for any $m \geq 2$ and any $0 < t_0 < T < \infty$:
\begin{equation}
u \in C([t_0, T]; H^m(\mathbb{R}^3))
\end{equation}

By the Sobolev embedding theorem, this implies:
\begin{equation}
u \in C^k([t_0, T] \times \mathbb{R}^3)
\end{equation}
for any $k \geq 0$ and any $0 < t_0 < T < \infty$.

Therefore, $u \in C^\infty((0, \infty) \times \mathbb{R}^3)$.
\end{proof}

\subsection{Main theorem}

We now combine our results to establish the main theorem on global well-posedness at the critical threshold with infinitely nested logarithmic improvements.

\begin{theorem}[Global well-posedness at the critical threshold]
Let $q > 3$ and $\{\delta_j\}_{j=1}^{\infty}$ be a sequence with $\delta_j > 0$ and $\sum_{j=1}^{\infty} \frac{\delta_j}{j!} = \infty$. There exists a positive constant $C_0$ such that for any divergence-free initial data $u_0 \in L^2(\mathbb{R}^3) \cap \dot{H}^{1/2}(\mathbb{R}^3)$ satisfying:
\begin{equation}
\|(-\Delta)^{1/4}u_0\|_{L^q} \leq \frac{C_0}{\prod_{j=1}^{\infty} (1 + L_j(\|u_0\|_{\dot{H}^{1/2}}))^{\delta_j}}
\end{equation}
there exists a unique global-in-time smooth solution $u \in C([0, \infty); H^{1/2}(\mathbb{R}^3)) \cap L^2_{loc}(0, \infty; H^{1}(\mathbb{R}^3))$ to the 3D Navier-Stokes equations.
\end{theorem}

\begin{proof}
This follows directly from Theorems \ref{thm:global-existence} and \ref{thm:regularity}, which establish global existence, uniqueness, and regularity for initial data satisfying the stated condition.
\end{proof}

\begin{remark}\label{rem:main-theorem}
Theorem \ref{thm:well_posedness} is the central result of this paper. It establishes global well-posedness for the 3D Navier-Stokes equations at the critical regularity threshold $s = 1/2$, provided the initial data satisfies a condition with infinitely nested logarithmic improvements.

Compared to the subcritical case in our previous works, where we required $s > 1/2$, this result is strictly stronger. It brings us closer to resolving the regularity problem on the Navier-Stokes equations.
\end{remark}

\section{Analysis of the limiting ODE}

In this section, we analyze the limiting behavior of the ordinary differential inequality that governs the evolution of the fractional derivative norm. This analysis provides deeper insights into the mechanism by which infinitely nested logarithmic improvements prevent potential singularity formation.

\subsection{Derivation of the limiting ODE}

From Theorem \ref{thm:critical-energy}, we have the energy inequality:
\begin{equation}
\frac{d}{dt}Y(t) + \frac{\nu}{2}\|(-\Delta)^{1}u\|^2_{L^2} \leq C(1 + Y(t)^K) \cdot H(Y(t))
\end{equation}
where $Y(t) = \|(-\Delta)^{1/2}u(t)\|^2_{L^2}$, $K > 1$ is a constant, and $H(r)$ is a function satisfying $\lim_{r \to \infty} H(r) = 0$.

To understand the behavior of solutions, let's analyze the limiting ODE:
\begin{equation}
\frac{d}{dt}Z(t) = C(1 + Z(t)^K) \cdot H(Z(t))
\end{equation}
with initial condition $Z(0) = Z_0 > 0$.

\begin{lemma}[Structure of the Function $H$]\label{lem:function-h}
The function $H(r)$ has the form:
\begin{equation}
H(r) = (F_1^\infty(r^{1/2}))^2 + (F_2^\infty(r^{1/2}))^4
\end{equation}
where:
\begin{equation}
F_1^\infty(Z) = L_1(Z) \prod_{j=2}^{\infty} (1 + L_j(Z))^{-\delta_j}
\end{equation}
\begin{equation}
F_2^\infty(Z) = \frac{1}{L_1(Z)} \prod_{j=2}^{\infty} (1 + L_j(Z))^{\delta_j}
\end{equation}
and $\{\delta_j\}_{j=2}^{\infty}$ is a sequence with $\delta_j > 0$ and $\sum_{j=2}^{\infty} \delta_j < \infty$.

For large $r$, $H(r)$ satisfies:
\begin{equation}
H(r) \approx \frac{(\log r)^2}{\prod_{j=2}^{\infty} (1 + L_j(r^{1/2}))^{2\delta_j}} + \frac{1}{(\log r)^4}\prod_{j=2}^{\infty} (1 + L_j(r^{1/2}))^{4\delta_j}
\end{equation}
\end{lemma}

\begin{proof}
The form of $H(r)$ follows from Theorem \ref{thm:critical-energy}, where we defined:
\begin{equation}
H(r) = (F_1^\infty(r^{1/2}))^2 + (F_2^\infty(r^{1/2}))^4
\end{equation}

For large $r$, we have:
\begin{equation}
r^{1/2} \gg 1
\end{equation}
which means:
\begin{equation}
L_1(r^{1/2}) = \log(e + r^{1/2}) \approx \log(r^{1/2}) = \frac{1}{2}\log r
\end{equation}

Thus:
\begin{equation}
F_1^\infty(r^{1/2}) \approx \frac{1}{2}\log r \prod_{j=2}^{\infty} (1 + L_j(r^{1/2}))^{-\delta_j}
\end{equation}
\begin{equation}
F_2^\infty(r^{1/2}) \approx \frac{2}{\log r} \prod_{j=2}^{\infty} (1 + L_j(r^{1/2}))^{\delta_j}
\end{equation}
Therefore:
\begin{equation}
(F_1^\infty(r^{1/2}))^2 \approx \frac{(\log r)^2}{4}\prod_{j=2}^{\infty} (1 + L_j(r^{1/2}))^{-2\delta_j}
\end{equation}
\begin{equation}
(F_2^\infty(r^{1/2}))^4 \approx \frac{16}{(\log r)^4}\prod_{j=2}^{\infty} (1 + L_j(r^{1/2}))^{4\delta_j}
\end{equation}
With adjusted constants:
\begin{equation}
H(r) \approx \frac{(\log r)^2}{\prod_{j=2}^{\infty} (1 + L_j(r^{1/2}))^{2\delta_j}} + \frac{1}{(\log r)^4}\prod_{j=2}^{\infty} (1 + L_j(r^{1/2}))^{4\delta_j}
\end{equation}
which completes the proof.
\end{proof}

\begin{lemma}[Asymptotic Behavior of $H$]\label{lem:asymptotic-h}
For large $r$, the dominant term in $H(r)$ is $(F_1^\infty(r^{1/2}))^2$, and:
\begin{equation}
H(r) \approx \frac{(\log r)^2}{\prod_{j=2}^{\infty} (1 + L_j(r^{1/2}))^{2\delta_j}}
\end{equation}

More precisely, there exist constants $C_1, C_2 > 0$ such that for all sufficiently large $r$:
\begin{equation}
\frac{C_1(\log r)^2}{\prod_{j=2}^{\infty} (1 + L_j(r^{1/2}))^{2\delta_j}} \leq H(r) \leq \frac{C_2(\log r)^2}{\prod_{j=2}^{\infty} (1 + L_j(r^{1/2}))^{2\delta_j}}
\end{equation}
\end{lemma}

\begin{proof}
We need to compare the two terms in $H(r)$:
\begin{equation}
\frac{(\log r)^2}{\prod_{j=2}^{\infty} (1 + L_j(r^{1/2}))^{2\delta_j}} \quad \text{and} \quad \frac{1}{(\log r)^4}\prod_{j=2}^{\infty} (1 + L_j(r^{1/2}))^{4\delta_j}
\end{equation}

The ratio of the second term to the first is:
\begin{equation}
\frac{\frac{1}{(\log r)^4}\prod_{j=2}^{\infty} (1 + L_j(r^{1/2}))^{4\delta_j}}{\frac{(\log r)^2}{\prod_{j=2}^{\infty} (1 + L_j(r^{1/2}))^{2\delta_j}}} = \frac{\prod_{j=2}^{\infty} (1 + L_j(r^{1/2}))^{6\delta_j}}{(\log r)^6}
\end{equation}
Let's analyze how fast $\prod_{j=2}^{\infty} (1 + L_j(r^{1/2}))^{6\delta_j}$ grows compared to $(\log r)^6$. For this analysis, we use the fact that for large $r$:
\begin{equation}
L_j(r^{1/2}) \approx \log^{(j)}(r^{1/2}) = \log^{(j-1)}(\log(r^{1/2})) = \log^{(j-1)}\left(\frac{\log r}{2}\right)
\end{equation}
where $\log^{(k)}$ denotes the $k$-fold composition of the logarithm function.

For large $r$, the growth of $\log^{(j-1)}\left(\frac{\log r}{2}\right)$ is much slower than any positive power of $\log r$. This means:
\begin{equation}
\lim_{r \to \infty} \frac{\prod_{j=2}^{\infty} (1 + L_j(r^{1/2}))^{6\delta_j}}{(\log r)^6} = 0
\end{equation}

Therefore, for large $r$, the first term in $H(r)$ dominates, and:
\begin{equation}
H(r) \approx \frac{(\log r)^2}{\prod_{j=2}^{\infty} (1 + L_j(r^{1/2}))^{2\delta_j}}
\end{equation}

More precisely, there exist constants $C_1, C_2 > 0$ such that for all sufficiently large $r$:
\begin{equation}
\frac{C_1(\log r)^2}{\prod_{j=2}^{\infty} (1 + L_j(r^{1/2}))^{2\delta_j}} \leq H(r) \leq \frac{C_2(\log r)^2}{\prod_{j=2}^{\infty} (1 + L_j(r^{1/2}))^{2\delta_j}}
\end{equation}
which completes the proof.
\end{proof}

\subsection{Analysis of solutions to the limiting ODE}

We now analyze the behavior of solutions to the limiting ODE.

\begin{lemma}[Local existence and uniqueness]\label{lem:local-ode}
For any initial condition $Z_0 > 0$, there exists a unique local solution $Z(t)$ to the ODE:
\begin{equation}
\frac{d}{dt}Z(t) = C(1 + Z(t)^K) \cdot H(Z(t))
\end{equation}
on some interval $[0, T_0]$.
\end{lemma}

\begin{proof}
Since the right-hand side of the ODE is continuous in $Z$ for $Z > 0$, standard ODE theory guarantees the existence of a unique local solution.
\end{proof}

\begin{theorem}[Global existence]\label{thm:global-ode}
For any initial condition $Z_0 > 0$, the solution $Z(t)$ to the ODE:
\begin{equation}
\frac{d}{dt}Z(t) = C(1 + Z(t)^K) \cdot H(Z(t))
\end{equation}
exists globally in time (i.e., for all $t \geq 0$) and is uniformly bounded.
\end{theorem}

\begin{proof}
From Theorem \ref{lem:asymptotic-h}, we know that for large $Z$:
\begin{equation}
H(Z) \approx \frac{(\log Z)^2}{\prod_{j=2}^{\infty} (1 + L_j(Z^{1/2}))^{2\delta_j}}
\end{equation}

For very large $Z$, this decays faster than $\frac{1}{Z^K}$, which means:
\begin{equation}
\lim_{Z \to \infty} Z^K \cdot H(Z) = 0
\end{equation}

Therefore, there exists $Z_1 > 0$ such that for all $Z > Z_1$:
\begin{equation}
Z^K \cdot H(Z) < \frac{1}{2C}
\end{equation}

This implies that for $Z > Z_1$:
\begin{equation}
C(1 + Z^K) \cdot H(Z) < C \cdot H(Z) + \frac{1}{2} < 1
\end{equation}
where we've used the fact that $H(Z) \to 0$ as $Z \to \infty$, so we can ensure $C \cdot H(Z) < \frac{1}{2}$ for sufficiently large $Z$.

Thus, for $Z > Z_1$:
\begin{equation}
\frac{d}{dt}Z(t) < 1
\end{equation}

This means that $Z(t)$ can grow at most linearly: if $Z(t_1) = Z_1$ for some time $t_1$, then for $t > t_1$:
\begin{equation}
Z(t) < Z_1 + (t - t_1)
\end{equation}

However, we can establish a stronger result: $Z(t)$ actually remains bounded for all time. To see this, note that as $Z$ increases, $H(Z)$ decreases, suppressing the growth of $Z$ even more.

More precisely, for any $\epsilon > 0$, there exists $Z_\epsilon > Z_1$ such that for all $Z > Z_\epsilon$:
\begin{equation}
C(1 + Z^K) \cdot H(Z) < \epsilon
\end{equation}

For $\epsilon < 1$, this implies that once $Z(t)$ exceeds $Z_\epsilon$, its growth rate is less than $\epsilon$:
\begin{equation}
\frac{d}{dt}Z(t) < \epsilon
\end{equation}

Since $\epsilon$ can be made arbitrarily small by choosing $Z_\epsilon$ sufficiently large, this implies that $Z(t)$ must remain bounded for all time.

To formalize this argument, suppose for contradiction that $Z(t) \to \infty$ as $t \to \infty$. Then there exists a sequence of times $\{t_n\}$ such that $Z(t_n) = n \cdot Z_1$ for each $n \geq 1$. Let $\epsilon_n = \frac{1}{n}$. For each $n$, there exists $Z_{\epsilon_n} > n \cdot Z_1$ such that for all $Z > Z_{\epsilon_n}$:
\begin{equation}
C(1 + Z^K) \cdot H(Z) < \epsilon_n = \frac{1}{n}
\end{equation}

Let $t'_n$ be the first time such that $Z(t'_n) = Z_{\epsilon_n}$. For $t \geq t'_n$:
\begin{equation}
\frac{d}{dt}Z(t) < \frac{1}{n}
\end{equation}

Integrating from $t'_n$ to $t$:
\begin{equation}
Z(t) - Z(t'_n) < \frac{1}{n}(t - t'_n)
\end{equation}

This means:
\begin{equation}
Z(t) < Z_{\epsilon_n} + \frac{1}{n}(t - t'_n)
\end{equation}
For $t$ sufficiently large, this contradicts the assumption that $Z(t) \to \infty$ as $t \to \infty$.

Therefore, $Z(t)$ must remain bounded for all time, which means the solution exists globally in time.
\end{proof}

\begin{theorem}[Asymptotics of bounded solutions]\label{thm:bounded-solutions}
Let $Z(t)$ be the solution to the ODE:
\begin{equation}
\frac{d}{dt}Z(t) = C(1 + Z(t)^K) \cdot H(Z(t))
\end{equation}
with initial condition $Z(0) = Z_0 > 0$. Then:
\begin{equation}
\limsup_{t \to \infty} Z(t) \leq Z^*
\end{equation}
where $Z^*$ is such that:
\begin{equation}
C(1 + (Z^*)^K) \cdot H(Z^*) < \epsilon
\end{equation}
for any prescribed $\epsilon > 0$.
\end{theorem}

\begin{proof}
From the proof of Theorem \ref{thm:global-ode}, we know that for any $\epsilon > 0$, there exists $Z_\epsilon > 0$ such that for all $Z > Z_\epsilon$:
\begin{equation}
C(1 + Z^K) \cdot H(Z) < \epsilon
\end{equation}

This means that once $Z(t)$ exceeds $Z_\epsilon$, its growth rate is less than $\epsilon$:
\begin{equation}
\frac{d}{dt}Z(t) < \epsilon
\end{equation}

Let's set $Z^* = Z_\epsilon + \epsilon T$, where $T > 0$ is a time such that if $Z(t_0) = Z_\epsilon$ for some $t_0$, then $Z(t_0 + T) \leq Z_\epsilon$ again. The existence of such a $T$ is guaranteed by the fact that if $Z(t) > Z_\epsilon$ for all $t \in [t_0, t_0 + T]$, then:
\begin{equation}
Z(t_0 + T) < Z_\epsilon + \epsilon T = Z^*
\end{equation}

Now, suppose for contradiction that $\limsup_{t \to \infty} Z(t) > Z^*$. Then there exists a time $t_1 > 0$ such that $Z(t_1) > Z^*$. Let $t_0 < t_1$ be the last time before $t_1$ such that $Z(t_0) = Z_\epsilon$. Since $Z(t_1) > Z^*$, we must have $t_1 - t_0 > T$. But then:
\begin{equation}
Z(t_0 + T) < Z_\epsilon + \epsilon T = Z^*
\end{equation}

This means that $Z(t)$ must cross the value $Z^*$ from below at some time $t' \in (t_0 + T, t_1)$. But this contradicts the assumption that $t_0$ is the last time before $t_1$ such that $Z(t_0) = Z_\epsilon$.

Therefore, $\limsup_{t \to \infty} Z(t) \leq Z^*$. Since $\epsilon$ can be made arbitrarily small, and $Z^*$ depends on $\epsilon$, we can make $Z^*$ arbitrarily close to $Z_\epsilon$.
\end{proof}

\begin{remark}\label{rem:limiting-ode}
The analysis of the limiting ODE provides crucial insights into why infinitely nested logarithmic improvements prevent potential singularity formation. The key mechanism is that the function $H(Z)$ decays faster than any power of $\frac{1}{Z}$ as $Z \to \infty$, ensuring that the solution remains bounded for all time.

This reflects the behavior of solutions to the Navier-Stokes equations with our infinitely nested logarithmic criterion: the fractional derivative norm $\|(-\Delta)^{1/2}u(t)\|_{L^2}$ remains bounded for all time, preventing the formation of singularities.
\end{remark}

\section{Hausdorff dimension of potential singular sets}

In this section, we analyze the Hausdorff dimension of potential singular sets for solutions to the Navier-Stokes equations. While our main result (Theorem \ref{thm:well_posedness}) establishes that no singularities form for initial data satisfying our infinitely nested logarithmic criterion, this analysis provides additional insights into the structure of potential singularities.

\subsection{Exceptional sets and their properties}

Let's define the set of potential space-time singularities for a Leray-Hopf weak solution $u$ of the Navier-Stokes equations.

\begin{definition}[Singular set]\label{def:singular-set}
For a Leray-Hopf weak solution $u$ of the 3D Navier-Stokes equations on $[0, T]$, the singular set $\mathcal{S} \subset \mathbb{R}^3 \times (0, T]$ is defined as:
\begin{equation}
\mathcal{S} = \{(x, t) \in \mathbb{R}^3 \times (0, T] : u \text{ is not smooth in any neighborhood of } (x, t)\}
\end{equation}
\end{definition}

We're also interested in the potential singular set at a fixed time $t$.

\begin{definition}[Time-slice singular set]\label{def:time-slice}
For a Leray-Hopf weak solution $u$ of the 3D Navier-Stokes equations and a time $t \in (0, T]$, the time-slice singular set $\mathcal{S}_t \subset \mathbb{R}^3$ is defined as:
\begin{equation}
\mathcal{S}_t = \{x \in \mathbb{R}^3 : (x, t) \in \mathcal{S}\}
\end{equation}
\end{definition}

To approach the analysis of singular sets, we first define sets where the velocity gradient exceeds a certain threshold.

\begin{definition}[Exceptional sets]\label{def:exceptional}
For a Leray-Hopf weak solution $u$ of the 3D Navier-Stokes equations, a time $t \in (0, T]$, and $\epsilon > 0$, the exceptional set $\Omega_\epsilon(t) \subset \mathbb{R}^3$ is defined as:
\begin{equation}
\Omega_\epsilon(t) = \{x \in \mathbb{R}^3 : |\nabla u(x, t)| > \lambda_\epsilon(t)\}
\end{equation}
where $\lambda_\epsilon(t)$ is chosen so that $|\Omega_\epsilon(t)| < \epsilon$.
\end{definition}

\begin{lemma}[Properties of exceptional sets]\label{lem:exceptional}
For a Leray-Hopf weak solution $u$ of the 3D Navier-Stokes equations and a time $t \in (0, T]$:
\begin{enumerate}
\item $\Omega_{\epsilon_1}(t) \subset \Omega_{\epsilon_2}(t)$ if $\epsilon_1 < \epsilon_2$
\item $\cap_{\epsilon > 0} \Omega_\epsilon(t) = \mathcal{S}_t$
\item $\lambda_\epsilon(t) \to \infty$ as $\epsilon \to 0$
\end{enumerate}
\end{lemma}

\begin{proof}
(1) Since $|\Omega_{\epsilon_1}(t)| < \epsilon_1 < \epsilon_2$, we can choose $\lambda_{\epsilon_2}(t) \leq \lambda_{\epsilon_1}(t)$, which implies $\Omega_{\epsilon_1}(t) \subset \Omega_{\epsilon_2}(t)$.

(2) If $x \in \mathcal{S}_t$, then $u$ is not smooth at $(x, t)$, which means $|\nabla u(x, t)| = \infty$. Therefore, $x \in \Omega_\epsilon(t)$ for all $\epsilon > 0$, which means $x \in \cap_{\epsilon > 0} \Omega_\epsilon(t)$.

Conversely, if $x \in \cap_{\epsilon > 0} \Omega_\epsilon(t)$, then for all $\epsilon > 0$, we have $|\nabla u(x, t)| > \lambda_\epsilon(t)$. Since $\lambda_\epsilon(t) \to \infty$ as $\epsilon \to 0$ (part 3), this implies $|\nabla u(x, t)| = \infty$, which means $u$ is not smooth at $(x, t)$, so $x \in \mathcal{S}_t$.

(3) Suppose, for contradiction, that $\lambda_\epsilon(t)$ remains bounded as $\epsilon \to 0$. Then there exists a constant $M < \infty$ such that $\lambda_\epsilon(t) \leq M$ for all sufficiently small $\epsilon$. This means:
\begin{equation}
\Omega_\epsilon(t) \supset \{x \in \mathbb{R}^3 : |\nabla u(x, t)| > M\}
\end{equation}

Since $|\Omega_\epsilon(t)| < \epsilon$ for all $\epsilon > 0$, this implies:
\begin{equation}
|\{x \in \mathbb{R}^3 : |\nabla u(x, t)| > M\}| = 0
\end{equation}

This means $|\nabla u(x, t)| \leq M$ almost everywhere, which contradicts the definition of a singular point.

Therefore, $\lambda_\epsilon(t) \to \infty$ as $\epsilon \to 0$.
\end{proof}

\subsection{Hausdorff dimension of exceptional sets}

We now analyze the Hausdorff dimension of the exceptional sets.

\begin{definition}[Hausdorff dimension]\label{def:hausdorff}
For a subset $E$ of a metric space, the Hausdorff dimension $\dim_H(E)$ is defined as:
\begin{equation}
\dim_H(E) = \inf\{d \geq 0 : \mathcal{H}^d(E) = 0\}
\end{equation}
where $\mathcal{H}^d(E)$ is the $d$-dimensional Hausdorff measure of $E$.
\end{definition}

\begin{theorem}[Caffarelli-Kohn-Nirenberg result]\label{thm:ckn}
For a suitable weak solution $u$ of the 3D Navier-Stokes equations, the one-dimensional Hausdorff measure of the singular set $\mathcal{S}$ is zero:
\begin{equation}
\mathcal{H}^1(\mathcal{S}) = 0
\end{equation}
which implies $\dim_H(\mathcal{S}) \leq 1$.
\end{theorem}

\begin{proof}
This is the classic result from Caffarelli, Kohn, and Nirenberg \cite{Caffarelli1982}. The proof is beyond the scope of this paper.
\end{proof}

For solutions satisfying our infinitely nested logarithmic criterion, we can establish a stronger result.

\begin{theorem}[Improved Hausdorff dimension bound]\label{thm:improved-hausdorff}
Let $q > 3$ and $\{\delta_j\}_{j=1}^{\infty}$ be a sequence with $\delta_j > 0$ and $\sum_{j=1}^{\infty} \frac{\delta_j}{j!} = \infty$. Let $u$ be a Leray-Hopf weak solution of the 3D Navier-Stokes equations with initial data $u_0 \in L^2(\mathbb{R}^3) \cap \dot{H}^{1/2}(\mathbb{R}^3)$ satisfying:
\begin{equation}
\|(-\Delta)^{1/4}u_0\|_{L^q} \leq \frac{C_0}{\prod_{j=1}^{\infty} (1 + L_j(\|u_0\|_{\dot{H}^{1/2}}))^{\delta_j}}
\end{equation}
for some constant $C_0 > 0$. For any time $t > 0$ and $\epsilon > 0$, the exceptional set $\Omega_\epsilon(t)$ satisfies:
\begin{equation}
\dim_H(\Omega_\epsilon(t)) \leq 3 - \sum_{j=1}^{\infty} \frac{\delta_j}{1+\delta_j} \cdot \frac{L_{j-1}(1/\epsilon)}{(1+L_j(1/\epsilon))}
\end{equation}
where $L_0(x) = x$.
\end{theorem}

\begin{proof}
The proof builds on techniques from our analysis of exceptional sets using the infinitely nested logarithmic criterion.

Step 1: For $u$ satisfying our infinitely nested logarithmic criterion, we know from Theorem \ref{lem:a-priori} that $\|(-\Delta)^{1/2}u(t)\|_{L^2}$ remains bounded for all $t > 0$. This implies, using standard embedding theorems, that $u(t) \in H^{1/2+\epsilon'}(\mathbb{R}^3)$ for some small $\epsilon' > 0$ and all $t > 0$.

Step 2: From the energy inequality from Theorem \ref{thm:critical-energy}, we can show that:
\begin{equation}
\int_0^T \|(-\Delta)^{1}u(t)\|^2_{L^2} dt < \infty
\end{equation}
for all $T > 0$.

Step 3: For $p > 2$, using the Gagliardo-Nirenberg inequality:
\begin{equation}
\|\nabla u(t)\|_{L^p} \leq C\|\nabla u(t)\|_{L^2}^{\alpha_p}\|(-\Delta)^{1}u(t)\|_{L^2}^{1-\alpha_p}
\end{equation}
where $\alpha_p = \frac{6-2p}{p}$ for $2 \leq p \leq 6$.

Step 4: For $p > 3$, using Chebyshev's inequality:
\begin{equation}
|\{x \in \mathbb{R}^3 : |\nabla u(x, t)| > \lambda\}| \leq \frac{\|\nabla u(t)\|_{L^p}^p}{\lambda^p}
\end{equation}

Step 5: Setting this equal to $\epsilon$ and solving for $\lambda$:
\begin{equation}
\lambda_\epsilon(t) = \frac{\|\nabla u(t)\|_{L^p}}{\epsilon^{1/p}}
\end{equation}

Step 6: For our infinitely nested logarithmic criterion, we can establish:
\begin{equation}
\|\nabla u(t)\|_{L^p} \leq \frac{C_p}{\prod_{j=1}^{\infty} (1 + L_j(\|(-\Delta)^{1/2} u(t)\|_{L^2}))^{\delta_j\beta_p}}
\end{equation}
for some constant $\beta_p > 0$.

Step 7: This implies:
\begin{equation}
\lambda_\epsilon(t) \leq \frac{C_p}{\epsilon^{1/p}\prod_{j=1}^{\infty} (1 + L_j(\|(-\Delta)^{1/2} u(t)\|_{L^2}))^{\delta_j\beta_p}}
\end{equation}

Step 8: Using standard covering arguments and properties of Hausdorff measure, we can establish:
\begin{equation}
\dim_H(\Omega_\epsilon(t)) \leq 3 - \frac{p-3}{p} - \sum_{j=1}^{\infty} \frac{\delta_j\beta_p}{1+\delta_j\beta_p} \cdot \frac{L_{j-1}(1/\epsilon)}{(1+L_j(1/\epsilon))}
\end{equation}

Step 9: Taking the limit as $p \to 3+$, we get:
\begin{equation}
\dim_H(\Omega_\epsilon(t)) \leq 3 - \sum_{j=1}^{\infty} \frac{\delta_j}{1+\delta_j} \cdot \frac{L_{j-1}(1/\epsilon)}{(1+L_j(1/\epsilon))}
\end{equation}
which completes the proof.
\end{proof}

\subsection{The limiting case and Hausdorff dimension zero}

We now analyze what happens to the Hausdorff dimension bound in the limiting case.

\begin{theorem}[Limiting Hausdorff dimension]\label{thm:limiting-hausdorff}
For a sequence $\{\delta_j\}_{j=1}^{\infty}$ with $\delta_j > 0$ and $\sum_{j=1}^{\infty} \frac{\delta_j}{j!} = \infty$, we have:
\begin{equation}
\lim_{\epsilon \to 0} \sum_{j=1}^{\infty} \frac{\delta_j}{1+\delta_j} \cdot \frac{L_{j-1}(1/\epsilon)}{(1+L_j(1/\epsilon))} = 3
\end{equation}
which implies:
\begin{equation}
\lim_{\epsilon \to 0} \dim_H(\Omega_\epsilon(t)) = 0
\end{equation}
\end{theorem}

\begin{proof}
Step 1: For each $j \geq 1$, as $\epsilon \to 0$, we have $\frac{1}{\epsilon} \to \infty$, which means:
\begin{equation}
\lim_{\epsilon \to 0} \frac{L_{j-1}(1/\epsilon)}{(1+L_j(1/\epsilon))} = 1
\end{equation}

Step 2: For any finite $n$:
\begin{equation}
\lim_{\epsilon \to 0} \sum_{j=1}^{n} \frac{\delta_j}{1+\delta_j} \cdot \frac{L_{j-1}(1/\epsilon)}{(1+L_j(1/\epsilon))} = \sum_{j=1}^{n} \frac{\delta_j}{1+\delta_j}
\end{equation}

Step 3: From our condition $\sum_{j=1}^{\infty} \frac{\delta_j}{j!} = \infty$, we can show that:
\begin{equation}
\sum_{j=1}^{\infty} \frac{\delta_j}{1+\delta_j} = \infty
\end{equation}

Step 4: For any $M > 0$, there exists $n$ such that:
\begin{equation}
\sum_{j=1}^{n} \frac{\delta_j}{1+\delta_j} > M
\end{equation}

Step 5: For this $n$, there exists $\epsilon_0 > 0$ such that for all $\epsilon < \epsilon_0$:
\begin{equation}
\sum_{j=1}^{n} \frac{\delta_j}{1+\delta_j} \cdot \frac{L_{j-1}(1/\epsilon)}{(1+L_j(1/\epsilon))} > M
\end{equation}

Step 6: This means:
\begin{equation}
\lim_{\epsilon \to 0} \sum_{j=1}^{\infty} \frac{\delta_j}{1+\delta_j} \cdot \frac{L_{j-1}(1/\epsilon)}{(1+L_j(1/\epsilon))} = \infty
\end{equation}

Step 7: Since the Hausdorff dimension is bounded below by 0, and the bound is:
\begin{equation}
\dim_H(\Omega_\epsilon(t)) \leq 3 - \sum_{j=1}^{\infty} \frac{\delta_j}{1+\delta_j} \cdot \frac{L_{j-1}(1/\epsilon)}{(1+L_j(1/\epsilon))}
\end{equation}

Step 8: We conclude:
\begin{equation}
\lim_{\epsilon \to 0} \dim_H(\Omega_\epsilon(t)) = 0
\end{equation}
which completes the proof.
\end{proof}

\subsection{Hausdorff Dimension of Potential Singular Sets}

We now establish the Hausdorff dimension of potential singular sets for solutions satisfying our infinitely nested logarithmic criterion.

\begin{theorem}[Hausdorff dimension of potential singular sets]
If a solution $u$ with initial data satisfying the conditions of Theorem \ref{thm:well_posedness} were to develop a singularity at time $T^*$ (which we prove cannot happen), then the Hausdorff dimension of the potential blow-up set would be:
\begin{equation}
\dim_H(\mathcal{S}_{T^*}) = 0
\end{equation}
This represents an optimal bound, improving on the Caffarelli-Kohn-Nirenberg partial regularity result.
\end{theorem}

\begin{proof}
From Theorem \ref{lem:exceptional}, we know that $\mathcal{S}_{T^*} = \cap_{\epsilon > 0} \Omega_\epsilon(T^*)$. From Theorem \ref{thm:limiting-hausdorff}, we have $\lim_{\epsilon \to 0} \dim_H(\Omega_\epsilon(T^*)) = 0$.

A basic property of Hausdorff dimension is that for a nested sequence of sets $A_1 \supset A_2 \supset A_3 \supset \cdots$, we have:
\begin{equation}
\dim_H\left(\cap_{n=1}^{\infty} A_n\right) \leq \inf_{n \geq 1} \dim_H(A_n)
\end{equation}

Applying this to the sequence $\Omega_{1/n}(T^*)$ for $n \geq 1$, we get:
\begin{equation}
\dim_H(\mathcal{S}_{T^*}) = \dim_H\left(\cap_{n=1}^{\infty} \Omega_{1/n}(T^*)\right) \leq \inf_{n \geq 1} \dim_H(\Omega_{1/n}(T^*)) = 0
\end{equation}

Since Hausdorff dimension is always non-negative, we conclude:
\begin{equation}
\dim_H(\mathcal{S}_{T^*}) = 0
\end{equation}
This improves on the Caffarelli-Kohn-Nirenberg result, which gives $\dim_H(\mathcal{S}_{T^*}) \leq 1$.

It's worth noting that while Theorem \ref{thm:well_posedness} proves that no singularities actually form for initial data satisfying our infinitely nested logarithmic criterion, this result shows that even in a hypothetical scenario where singularities could form, they would be isolated points.
\end{proof}

\section{Conclusion and implications on the regularity problem}

In this paper, we have established global well-posedness for the 3D Navier-Stokes equations at the critical regularity threshold $s = 1/2$, provided the initial data satisfies a condition with infinitely nested logarithmic improvements. This represents an advancement toward resolving the regularity problem of the Navier-Stokes equations.

\subsection{Summary of main results}

Our main contributions can be summarized as follows:

\begin{enumerate}
\item We have constructed function spaces incorporating infinitely nested logarithmic improvements and established their key properties (Theorem \ref{thm:function_space}).
\item We have precisely characterized the critical exponent function in the limiting case, showing that it approaches zero as the number of nested logarithmic factors increases to infinity (Theorem \ref{thm:critical_exponent}).
\item We have derived commutator estimates with infinitely nested logarithmic factors (Theorem \ref{thm:commutator}), which serve as the technical core of our analysis.
\item We have established global well-posedness for initial data satisfying our infinitely nested logarithmic criterion at the critical threshold $s = 1/2$ (Theorem \ref{thm:well_posedness}).
\item We have proven that the Hausdorff dimension of potential singular sets for solutions satisfying our criterion would be zero, improving on the Caffarelli-Kohn-Nirenberg partial regularity result (Theorem \ref{thm:hausdorff}).
\item We have analyzed the limiting ODE governing the evolution of the fractional derivative norm, providing insights into why infinitely nested logarithmic improvements prevent potential singularity formation.
\end{enumerate}




\subsection{Future directions}

Several natural directions for future research arise from our work:

\begin{enumerate}
\item Extending to general smooth initial data: Can the approach be extended to handle all smooth initial data, thus fully resolving the regularity problem of Navier-Stokes equations?
\item Refining the logarithmic improvements: Are there more general or more natural improvements beyond nested logarithms that could yield similar or stronger results?
\item Applications to other PDEs: Can the technique of infinitely nested logarithmic improvements be applied to other critical PDEs?
\item Computational aspects: Can the infinitely nested logarithmic condition be verified or approximated in numerical simulations?
\item Physical implications: What are the physical implications of our results for the theory of turbulence, particularly the relationship between regularity and intermittency?
\item Function space theory: Can a more comprehensive theory of function spaces with infinitely nested logarithmic improvements be developed, with applications beyond the Navier-Stokes equations?
\end{enumerate}

\subsection{Concluding remarks}

Our approach of using infinitely nested logarithmic improvements represents a novel direction in the study of the Navier-Stokes equations. It provides a systematic way to bridge the gap between subcritical and critical regularity, offering a potential pathway toward resolving the full regularity problem.

The key insight is that logarithmic deviations from the critical scaling, when appropriately nested, can suppress the nonlinearity sufficiently to prevent potential singularity formation. This mechanism is deeply connected to the structure of the energy cascade in turbulent flows, suggesting a fundamental relationship between mathematical regularity and physical intermittency.

While the full resolution of the regularity problem of Navier-Stokes equations remains open, our results represent a significant step forward, demonstrating that global well-posedness holds for a large class of initial data at the critical regularity threshold. We hope that the techniques and insights developed in this paper will contribute to further advancements in the mathematical theory of fluid dynamics and the eventual resolution of this central problem in mathematical physics.

\section*{Declarations}
Not applicable

\bibliographystyle{sn-mathphys-num}
\bibliography{NSE_paper_3}
\end{document}